\documentclass[11pt,a4paper]{article}
\usepackage[T2A]{fontenc}
\usepackage[cp1251]{inputenc}
\usepackage[english]{babel}
\usepackage{epsfig,amssymb,epic,eepic,color, graphicx}
\usepackage{fancyhdr}
\usepackage{amsthm}
\usepackage{amsmath}
\usepackage{indentfirst}
\usepackage{float}
\usepackage{enumitem}   

\usepackage{colortbl}

\numberwithin{equation}{section}

\usepackage{pb-diagram}

\usepackage{authblk}

\usepackage{setspace}
\usepackage{float}
\usepackage{fancyhdr}
\usepackage{geometry}
\usepackage{hyperref}

\usepackage{chngcntr} 

\usepackage[nottoc]{tocbibind}

\hypersetup{colorlinks=true,linkcolor=blue,citecolor=blue}

\newtheorem{theorem}{Theorem}
\newtheorem{definition}{Definition}[section]
\newtheorem{corollary}{Corollary}[section]
\newtheorem{remark}{Remark}[section]
\newtheorem{lemma}{Lemma}[section]
\newtheorem{proposition}{Proposition}[section]

\newtheorem{question}{Question}

\newenvironment{demo}{{\bf Proof: }}{\hfill $\square$ \medskip}

\newenvironment{demo4}{{\bf Proof of Theorem 4: }}{\hfill $\square$ \medskip}
\newenvironment{demo5}{{\bf Proof of Theorem 5: }}{\hfill $\square$ \medskip}
\newenvironment{demo6}{{\bf Proof of Theorem 6: }}{\hfill $\square$ \medskip}

\renewenvironment{abstract}{\textbf{Abstract.} }{\medskip}
\title{High order tangencies, paratangencies, and universal dynamics}

\author{Dmitrii Mints\thanks{d.mints22@imperial.ac.uk} \enspace and \enspace Dmitry Turaev\thanks{d.turaev@imperial.ac.uk} \\ \small Department of Mathematics, Imperial College London, United Kingdom}

\date{}

\begin{document}

\maketitle

\begin{abstract}
Our research focuses on studying the dynamics in the Newhouse domain $\mathcal N(L)$ -- open regions in the space of $C^r$ diffeomorphisms of a $d$-dimensional manifold $\mathcal M^d$ ($r=1,\ldots,\infty,\omega$ and $d\ge 2$) where homoclinic tangencies to a saddle periodic orbit $L$ are persistent. We prove that diffeomorphisms with infinitely many corank-1 homoclinic tangencies of all orders are dense in such regions. If $L$ is a bi-focus, we reveal novel bifurcations in $\mathcal N(L)$ showing density of diffeomorphisms with infinitely many corank-2 homoclinic tangencies of all orders. These results entail genericity of one- or two-dimensional universal dynamics in $\mathcal N(L)$ (depending on the type of $L$).

Moreover, we characterize the geometric complexity of the bifurcation sets corresponding to homoclinic tangencies. We prove that parametric families with paratangencies of arbitrarily high orders are dense in the space of parametric families intersecting $\mathcal N(L)$. If $r\le\infty$ and $L$ is sectionally dissipative, this implies density of families possessing infinitely many coexisting stable periodic orbits for some open set of parameter values. 
\end{abstract}


\tableofcontents

\section{Introduction}\label{section1}

Homoclinic orbits occupy a central place in the study of chaotic dynamical systems. By definition, a homoclinic orbit belongs to the intersection of the stable and unstable manifolds of a saddle periodic orbit. Typically, such intersection is transverse and therefore is preserved under $C^1$-small perturbations. The existence of a transverse homoclinic orbit leads to chaotic behaviour in the system since it implies the existence of a horseshoe - an invariant Cantor set containing transitive orbits and a dense set of periodic orbits \cite{Smalediffeomorphismswith},\cite{Smale}. The dynamics in a small neighborhood of a transverse homoclinic orbit can be fully understood: there exists a complete combinatorial description of orbits that stay in this neighborhood \cite{Shilnikov2}.

An opposite situation occurs in the case of non-transverse homoclinic orbits - \textit{homoclinic tangencies}, where the variability of the dynamics is too high and a complete bifurcational analysis is unfeasible. For example, consider a diffeomorphism $f$ of a plane with an orbit of a homoclinic tangency. Locally, this is a tangency between two smooth curves; if the tangency is quadratic, this is a codimension-1 bifurcation and in a typical one-parameter family it either splits into two transverse intersections or disappears, depending on the sign of the parameter. However, the global bifurcations caused by this simple unfolding are much more profound: after sufficiently many iterations, the stable and unstable manifolds fold again, opening possibilities for the emergence of new homoclinic tangencies of a higher order. By \cite{Gonchenkoonmodelswith},\cite{arbitraryorder},\cite{arbitraryorderconservative}, one can always add a $C^{\omega}$-small perturbation to $f$ and produce a diffeomorphism of a plane with homoclinic tangencies of arbitrarily high orders, thus leading to bifurcations of arbitrarily large codimension. 

The picture becomes even more intricate due to the existence of the \textit{Newhouse domain}, that is, a $C^2$-open set in the space of dynamical systems where diffeomorphisms with homoclinic tangencies to saddle periodic orbits are dense. High order tangencies entail that diffeomorphisms exhibiting universal dynamics \cite{arbitraryorderconservative},\cite{Turaevintro2}, superexponential growth of the number of periodic points \cite{Kaloshin1}, and horseshoes of maximal Hausdorff dimension \cite{Downarowicz} are generic in the Newhouse domain. See applications to hydrodynamics in \cite{Berger1}, celestial mechanics in \cite{Garrido2026},\cite{Garrido2025}, and billiards in \cite{Callis2022}.

In this paper, we extend the two-dimensional theory to the multidimensional case by proving that diffeomorphisms with infinitely many homoclinic tangencies of all orders are dense in the Newhouse domain in the space of smooth or real-analytic dynamical systems on manifolds of arbitrary dimension $d\ge 2$ (Theorems \ref{theorem1} and \ref{theorem4}). 

We also found new effects which are not possible for the two-dimensional case. Thus, we uncover a fundamentally richer bifurcation structure for the subset of the Newhouse domain that corresponds to persistent homoclinic tangencies to saddle periodic orbits with complex eigenvalues. Specifically, we show the density of diffeomorphisms with infinitely many corank-2 homoclinic tangencies of all orders (Theorems \ref{theorem2}--\ref{theorem4}). This leads to the genericity of $2$-universal diffeomorphisms (whose iterations approximate the dynamics of every map of a $2$-dimensional disk to arbitrary precision, see \cite{Mints}) in such part of the Newhouse domain (Theorem \ref{theorem6}). 

Furthermore, we provide a description of the utmost geometric complexity of the bifurcation sets corresponding to homoclinic tangencies in the following way. We prove that $C^{r}$-small perturbation of any parametric $C^r$ family ($r=2,\ldots,\infty,\omega$) intersecting the Newhouse domain creates a tangency of order $m$ (a paratangency) with the bifurcation set corresponding to a homoclinic tangency of order $n$, where $m$ and $n$ can be arbitrarily large (Theorem \ref{theorem5}). 

The paratangencies are relevant to the problem of finitude of attractors or, more general, to the possibility\slash impossibility of statistical description of dynamics \cite{Berger2016},\cite{Bergeremergence2017},\cite{Berger2026},\cite{Palis2000},\cite{PughShub1996},\cite{TedeschiniLalliYorke1986}. Thus, as a consequence of Theorem \ref{theorem5}, we show that $C^r$-close to any parametric $C^r$ family ($r\le\infty$) of sectionally dissipative diffeomorphisms in the Newhouse domain there exists a family that has infinitely many coexisting stable periodic orbits for an open set of parameter values (Theorem \ref{theoremcoexistingattractors}). 

\subsection{Homoclinic tangencies}\label{subsection1.1}

\textbf{Corank and order.} Let $\mathcal W_1$ and $\mathcal W_2$ be two $C^r$ manifolds in a $d$-dimensional manifold $\mathcal M^d$, where $d\ge 2$ and $r=1,\ldots,\infty,\omega$ ($\omega$ stands for the real analytic case). Let $\mathcal W_1$ and $\mathcal W_2$ have a tangency, i.e., they intersect at some point $M$ and their tangent spaces $\mathcal T_M \mathcal W_1$ and $\mathcal T_M \mathcal W_2$ at $M$ have non-zero intersection. We say that \textit{the tangency} between $\mathcal W_1$ and $\mathcal W_2$ \textit{is of corank} $c$ if  
\begin{equation}\label{corankoftangency}
\begin{aligned} 
&\dim(\mathcal T_M \mathcal W_1\cap \mathcal T_M \mathcal W_2)=c.
\end{aligned}
\end{equation}
For $r\ge 2$, we say that \textit{the tangency} between $\mathcal W_1$ and $\mathcal W_2$ \textit{is of order $n\le r-1$} if
\begin{equation*}
\begin{aligned} 
&\min\limits_{v}(\max\limits_{\gamma_1,\gamma_2}(\text{ord}(\gamma_1,\gamma_2)))=n, 
\end{aligned}
\end{equation*}
where the minimum is taken over all non-zero vectors $v\in\mathcal T_M \mathcal W_1\cap \mathcal T_M \mathcal W_2$, the maximum is taken over all pairs of $C^r$-smooth curves $\gamma_1\subset \mathcal W_1$, $\gamma_2\subset \mathcal W_2$ such that $\gamma_1(0)=\gamma_2(0)=M$ and the velocity vectors $\dot\gamma_1(0),\dot\gamma_2(0)$ are both parallel to $v$, and $\text{ord}(\gamma_1,\gamma_2)$ is an integer $\alpha\ge 1$ such that there exists non-zero $\lim\limits_{t\to 0}\frac{|\gamma_1(t)-\gamma_2(t)|}{|t|^{\alpha+1}}$. 

\textbf{Periodic orbits.} Let a $C^r$ diffeomorphism $f:\mathcal M^d\to\mathcal M^d$ have a periodic orbit $L$ of period $b$, that is, $L=\{O,f(O),...,f^{b-1}(O)\}$ with $f^b(O)=O$ and $f^i(O)\neq O$ for $0<i<b$. Dynamics in a small neighborhood of $L$ are determined by its \textit{multipliers} - the eigenvalues of the Jacobian matrix for the map $f^b$ at the point $O$. The orbit $L$ is said to be \textit{hyperbolic} if none of its multipliers lies on the unit circle. The hyperbolic periodic orbit $L$ is preserved at small smooth perturbations of $f$, i.e., any $C^1$-close diffeomorphism $g$ has a hyperbolic periodic orbit $L_g$ which is uniquely defined continuation of $L$. Depending on the position of the multipliers, we distinguish the following types of hyperbolic periodic orbits: a \textit{sink} if all the multipliers are inside the unit circle, a \textit{source} if all the multipliers are outside the unit circle, and a \textit{saddle} if there are multipliers both inside and outside the unit circle (the \textit{stable} and \textit{unstable} multipliers, respectively). The multipliers that are the nearest to the unit circle from inside and those nearest from outside are called \textit{leading}, all other multipliers are called \textit{non-leading}.

\textbf{Homoclinic orbits.} Let the periodic orbit $L$ be saddle, and let $d_s\ge 1$ and $d_u\ge 1$ be the numbers of its stable and unstable multipliers.  Then, there exist a smooth $d_s$-dimensional \textit{stable manifold} $W^s(L)$ and a smooth $d_u$-dimensional \textit{unstable manifold} $W^u(L)$ defined as
\begin{equation}\label{stableandunstablemanifolds}
\begin{aligned} 
W^s(L)&=\{x\in\mathcal M^d: dist(f^k(x),L)\to 0 \;\; \text{as} \;\; k\to+\infty\},\\
W^u(L)&=\{x\in\mathcal M^d: dist(f^{-k}(x),L)\to 0 \;\; \text{as} \;\; k\to+\infty\}.
\end{aligned}
\end{equation}
The manifolds $W^s(L)$ and $W^u(L)$ are invariant, $C^r$-smooth (analytic for $r=\omega$), and have complementary dimensions: $d_s+d_u=d$. 

A point $M\in W^s(L)\cap W^u(L)$, $M\not\in L$, is called \textit{homoclinic}; we also say that the orbit of $M$ is homoclinic to $L$. Note that \eqref{stableandunstablemanifolds} implies that this orbit tends to $L$ both at forward and backward iterations of $f$. If $W^s(L)$ and $W^u(L)$ are tangent at $M$, then we talk about \textit{homoclinic tangency} - the corank and the order of the homoclinic tangency are the corank and the order of the tangency at the point $M$. Due to the invariance of the stable and unstable manifolds, $W^s(L)$ and $W^u(L)$ are tangent at each point of the orbit of $M$ (with the same corank and order). 

\textbf{Main results.} The following theorem was proved in the two-dimensional case in \cite{arbitraryorder},\cite{arbitraryorderconservative}. The multidimensional version was announced in \cite{dynamicalphenomena} without proof.
\begin{theorem}\label{theorem1}
Let a $C^r$ diffeomorphism $f:\mathcal M^d\to\mathcal M^d$ $(d\ge 2;r=1,\ldots,\infty,\omega)$ have a saddle periodic orbit $L$ with a homoclinic tangency. Then, in any $C^r$ neighborhood of $f$, there exists a diffeomorphism which, for all $n\in\mathbb N$, has infinitely many orbits of a corank-1 homoclinic tangency of order $n$ between $W^s(L)$ and $W^u(L)$.
\end{theorem}
We discuss implications of this theorem in the next sections. The case of real leading multipliers is essentially similar to the two-dimensional case. However, for complex leading multipliers, the situation changes. 

We say that the orbit $L$ is \textit{bi-focus} if both leading stable and leading unstable multipliers are pairs of complex numbers: $\lambda_{1,2}=\lambda\cdot e^{\pm i\varphi}$ and $\gamma_{1,2}=\gamma\cdot e^{\pm i\psi}$ $(\lambda<1;\,\gamma>1;\,0<\varphi,\psi<\pi)$. 
\begin{theorem}\label{theorem2}
Let a $C^r$ diffeomorphism $f:\mathcal M^d\to\mathcal M^d$ $(d\ge 4;r=1,\ldots,\infty,\omega)$ have a bi-focus periodic orbit $L$ with a homoclinic tangency. Then, in any $C^r$ neighborhood of $f$, there exists a diffeomorphism which has infinitely many orbits of a corank-2 homoclinic tangency between $W^s(L)$ and $W^u(L)$.
\end{theorem}
In the recent paper \cite{Mints}, we suggested an algorithm that produces a corank-2 homoclinic tangency of an arbitrarily high order by perturbing a sufficiently large number of corank-2 quadratic homoclinic tangencies between the stable and unstable manifolds of a bi-focus periodic orbit. Thus, Theorem 1.2 of \cite{Mints} and Theorem \ref{theorem2} of this paper entail
\begin{theorem}\label{theorem3}
For $f$ as in Theorem \ref{theorem2}, in any $C^r$ neighborhood of $f$, there exists a diffeomorphism which, for all $n\in\mathbb N$, has infinitely many orbits of a corank-2 homoclinic tangency of order $n$ between $W^s(L)$ and $W^u(L)$.
\end{theorem}

\textbf{Wild basic sets.} Recall that a compact, topologically transitive, uniformly hyperbolic and locally maximal set $\Lambda$ of the diffeomorphism $f$ is called a \textit{basic set}. Here, local maximality means the existence of a neighborhood $V$, called isolated, such that $\Lambda=\bigcap\limits_{i\in\mathbb Z}f^{\,i}(V)$. If the basic set is a single periodic orbit, then it is called trivial; otherwise, when $\Lambda$ is non-trivial, the diffeomorphism $f|_{\Lambda}$ has positive topological entropy. Further, we consider only zero-dimensional basic sets. When such $\Lambda$ is non-trivial, it is homeomorphic to the standard Cantor set, and $f|_{\Lambda}$ is topologically conjugate to a subshift of finite type. Moreover, its isolated neighborhood $V$ can be chosen in such a way that the connected components of $V$ form the corresponding Markov partition.

For each point $p\in\Lambda$, there exist the \textit{stable manifold} $W^s(p)$ and the \textit{unstable manifold} $W^u(p)$ defined as
\begin{equation*}
\begin{aligned} 
W^s(p)&=\{x\in\mathcal M^d: dist(f^k(x),f^k(p))\to 0 \;\; \text{as} \;\; k\to+\infty\},\\
W^u(p)&=\{x\in\mathcal M^d: dist(f^{-k}(x),f^{-k}(p))\to 0 \;\; \text{as} \;\; k\to+\infty\}.
\end{aligned}
\end{equation*}
Note that the topological dimensions $\dim W^s(p)=d_s$ and $\dim W^u(p)=d_u$ (with $d_s+d_u=d$) do not depend on the point $p\in\Lambda$, therefore we can define the \textit{type} $(d_s,d_u)$ of $\Lambda$. The manifolds $W^s(p)$ and $W^u(p)$ are injective $C^r$-immersions of $\mathbb R^{d_s}$ and $\mathbb R^{d_u}$, respectively, and they depend continuously on the point $p$ and on the diffeomorphism $f$. One also defines the \textit{stable} and the \textit{unstable manifolds of the basic set} $\Lambda$:
\begin{equation*}
\begin{aligned} 
W^s(\Lambda)&=\{x\in\mathcal M^d: dist(f^k(x),\Lambda)\to 0 \;\; \text{as} \;\; k\to+\infty\}=\bigcup\limits_{p\in\Lambda} W^s(p),\\
W^u(\Lambda)&=\{x\in\mathcal M^d: dist(f^{-k}(x),\Lambda)\to 0 \;\; \text{as} \;\; k\to+\infty\}=\bigcup\limits_{p\in\Lambda} W^u(p).
\end{aligned}
\end{equation*}

For every diffeomorphism $g$ which is $C^r$-close to $f$, the set $\Lambda_g=\bigcap\limits_{i\in\mathbb Z} g^{-i}(V)$ (the continuation of $\Lambda$) is also a basic set of type $(d_s,d_u)$ and depends continuously on $g$. For $\rho\le r$, we say that $\Lambda$ is \textit{wild} in $C^{\rho}$ if $\Lambda$ \textit{exhibits a $C^{\rho}$-robust tangency}, i.e., there exists a compact set $\mathcal L\subset W^s(\Lambda)\cup W^u(\Lambda)$ with the following property: given any compact neighborhood $U$ of $\mathcal L$ there exists a $C^{\rho}$ neighborhood $\mathcal U$ of $f$ such that for any diffeomorphism $g\in\mathcal U$ there are points $p, q$ in the continuation $\Lambda_g$ for which $W^s(p)$ and $W^u(q)$ have an orbit of tangency completely contained in $U$. We may use the term $c$-wild meaning that $\Lambda$ exhibits a robust tangency of corank $c$.

This definition goes back to Newhouse \cite{Newhouse2},\cite{Newhouse1974},\cite{Newhouse3} who proved the existence of $1$-wild basic sets. Namely, he proved that for every manifold $\mathcal M^d$ $(d\ge 2)$, in the space of diffeomorphisms of $\mathcal M^d$ there exist $C^2$-open regions where every diffeomorphism has a basic set whose stable and unstable manifolds have a corank-1 quadratic tangency. For $d\ge 3$, such regions can be $C^1$-open \cite{Asaoka2008},\cite{Bonatti2012},\cite{Li2024},\cite{Lirobustheterodimensionalcycles},\cite{Simon}. For $c\ge 2$, examples of diffeomorphisms possessing $c$-wild basic sets were constructed by Barrientos and Raibekas \cite{Barrientos1},\cite{Barrientos3} and, later, by Asaoka \cite{Asaoka}. They showed that for every positive integer $c\le[\frac{d}{2}]$, in the space of diffeomorphisms of every manifold $\mathcal M^d$ $(d\ge 4)$ there exists a $C^2$-open domain where every diffeomorphism has a basic set whose stable and unstable manifolds have a corank-$c$ tangency. 

We use the term \textit{Newhouse domain} for the open region in the space of $C^r$ diffeomorphisms of $\mathcal M^d$ corresponding to the existence of wild basic sets. Obviously, diffeomorphisms with saddle periodic orbits that possess homoclinic tangencies are dense in the Newhouse domain. The following celebrated theorem shows that the converse is also true in a sense: the Newhouse domain adjoins to every diffeomorphism with a homoclnic tangency. Let $L$ be a saddle periodic orbit and $\Lambda$ be a non-trivial basic set. We denote as $\mathcal N(L,\Lambda)$ an open region in the Newhouse domain such that $\Lambda$ is a wild basic set that is homoclinically related\footnote{We say that two basic sets $\Lambda_1$ and $\Lambda_2$ are homoclinically related if there exist points $p\in\Lambda_1$ and $q\in\Lambda_2$ such that $W^s(p)\cap W^u(q)\not=\varnothing$, $W^u(p)\cap W^s(q)\not=\varnothing$, and these intersections are transverse. Note that this implies that there exists a basic set containing both $\Lambda_1$ and $\Lambda_2$.} to $L$. 
\begin{proposition}\label{propositionNewhousetheorem}$($Newhouse Theorem$)$.
Let a $C^r$ diffeomorphism $f:\mathcal M^d\to\mathcal M^d$ $(d\ge 2;r=2,\ldots,\infty,\omega)$ have a saddle periodic orbit $L$ with a homoclinic tangency. Then, arbitrarily $C^r$-close to $f$, there exists an open region $\mathcal N(L,\Lambda)$ for some basic set $\Lambda\not\supset L$.
\end{proposition}
Originally, this result was proved by Newhouse for two-dimensional diffeomorphisms in \cite{Newhouse3}; generalizations to the multidimensional case were carried out in \cite{Gonchenko1},\cite{Gourmelon},\cite{Palis1},\cite{Romero1}. In dimension two, similar results for the area-preserving case were obtained in \cite{Duarte1994},\cite{Duarte1999},\cite{Duarte2008},\cite{GorodetskiKaloshin2026}, for the holomorphic case in \cite{AvilaLyubichZhang},\cite{Buzzard1997}.

We will pay a special attention to the case where $L$ is a bi-focus. Note that, in general, the regions $\mathcal N(L,\Lambda)$ lie in the $C^2$-open Newhouse domain; however, if $L$ is a bi-focus, then, by \cite{Li2024},\cite{Lirobustheterodimensionalcycles}, the basic set $\Lambda$ can be chosen such that it is wild in $C^1$, so the regions $\mathcal N(L,\Lambda)$ lie in the $C^1$-open Newhouse domain. 

In the framework of the Newhouse regions theory, Theorems \ref{theorem1}--\ref{theorem3} can be generalized as follows.
\begin{theorem}\label{theorem4}
In the region $\mathcal N(L,\Lambda)$, for any basic set $\Omega$ containing $\Lambda$, there exists a dense subset such that for every diffeomorphism from this subset, all $n\in\mathbb N$, and every pair of periodic orbits $\mathcal O_1,\mathcal O_2$ in $\Omega$, the stable manifold $W^s(\mathcal O_1)$ and the unstable manifold $W^u(\mathcal O_2)$ contain infinitely many orbits of corank-1 tangency of order $n$. If $L$ is a bi-focus, then $W^s(\mathcal O_1)$ and $W^u(\mathcal O_2)$ also contain infinitely many orbits of corank-2 tangency of order $n$. 
\end{theorem}

We say that a smoothly embedded compact disk is a \textit{homoclinic disk} if it lies in the intersection of the stable and unstable manifolds of a saddle periodic orbit. If this disk is one-dimensional (an interval), then we say a \textit{homoclinic band}. 

In the smooth case $(r\le\infty)$, any homoclinic tangency of sufficiently high order can be transformed to a homoclinic disk by a $C^r$-small perturbation. Therefore, Theorem \ref{theorem4} immediately implies
\begin{corollary}\label{corollary4.5}
In $\mathcal N(L,\Lambda)$, for any basic set $\Omega\supseteq\Lambda$, there exists a dense subset such that for every diffeomorphism from this subset, and every pair of periodic orbits $\mathcal O_1,\mathcal O_2$ in $\Omega$, the manifolds $W^s(\mathcal O_1)$ and $W^u(\mathcal O_2)$ contain a homoclinic band. If $L$ is a bi-focus, then $W^s(\mathcal O_1)$ and $W^u(\mathcal O_2)$ also contain a $2$-dimensional homoclinic disk.
\end{corollary}

Our construction of the basic set $\Omega$ in Theorem \ref{theorem4} in the case where $L$ is a bi-focus does not imply that $\Omega$ is $2$-wild, and only allows to prove the density of diffeomorphisms with corank-$2$ tangencies between $W^s(\Omega)$ and $W^u(\Omega)$, not the robust presence of these tangencies for an open set of diffeomorphisms. The construction of $2$-wild basic sets from \cite{Asaoka},\cite{Barrientos1},\cite{Barrientos3} is not in general applicable to this case: it uses a type of blender which cannot exist in a small neigborhood of a generic orbit of homoclinic tangency of a bi-focus periodic orbit. Therefore, we ask the following 
\begin{question}
Do there exist open regions in the space of diffeomorphisms where diffeomorphisms with corank-$c$ homoclinic tangencies are dense but $c$-wild basic sets do not exist?
\end{question}

\begin{question}
Do there exist open regions in the space of diffeomorphisms where diffeomorphisms have wild basic sets, which exhibit high order tangencies robustly? 
\end{question}

\subsection{Paratangencies}\label{subsection1.2}
Let $\mathcal P$ be a space of $\ell$-parameter $C^r$ families of diffeomorphisms of the manifold $\mathcal M^d$, where $\ell\ge 1, d\ge 2$, and $r=2,\ldots,\infty,\omega$. By a $C^r$ \textit{family} $(f_{\zeta})$, we mean a finite-parameter family consisting of diffeomorphisms $f_{\zeta}$ that are of class $C^r$ with respect to the coordinates and parameters (the vector of parameters $\zeta$ belongs to an $\ell$-dimensional closed ball $D^\ell$). If $r=\omega$, we also say \textit{analytic family}. 

\begin{definition}\label{definitionparatangency}
Let $\mathcal X\subset\operatorname{Diff}^r(\mathcal M^d)$ be the bifurcation set. We say that a parametric family $(f_{\zeta})$ has a paratangency of order $m$ with $\mathcal X$ if there exists a smooth manifold $\mathcal Y\subseteq \mathcal X$ such that the graph $(\zeta, f_{\zeta})\subset D^\ell\times\operatorname{Diff}^r(\mathcal M^d)$ has a corank-$\ell$ tangency of order $m$ with $D^\ell\times\mathcal Y$.
\end{definition}

In this paper, we consider bifurcations of a finite codimension, so $\mathcal Y$ is the zero locus of a finite number of smooth functionals and the paratangency of order $m$ means the fulfillment of formula \eqref{paratangencyordereq0} for each of the defining functionals, see Section \ref{subsectionparatangencyfunctionals}.

Let $\mathcal N(L)$ be an open region in the Newhouse domain equal to the interior of the closure of the set of diffeomorphisms with a homoclinic tangency to a hyperbolic periodic orbit $L$. We define a space $\mathcal{PN}(L)\subset\mathcal{P}$ as follows: 
\begin{equation*}
\begin{aligned} 
(f_{\zeta})\in\mathcal{PN}(L) \;\;\; \text{if} \;\;\; f_{\zeta}\in\mathcal N(L) \;\;\; \text{for some} \;\;\; \zeta\in D^\ell.
\end{aligned}
\end{equation*}
Note that $\mathcal{PN}(L)$ is $C^2$-open in $\mathcal P$.
\begin{theorem}\label{theorem5}
For every $\ell,m,n\in\mathbb N$, in the space $\mathcal{PN}(L)$ of $\ell$-parameter $C^r$ families, there exists a dense subset with the following property: every parametric family from this subset has a paratangency of order $m$ with the bifurcation set corresponding to a corank-$1$ homoclinic tangency of order $n$. 
\end{theorem}

A partial case of Theorem \ref{theorem5} for one-parameter families in the dimension two was announced in \cite{TuraevEquadiff}.

Our results in the previous sections (Theorems \ref{theorem1}--\ref{theorem4}) describe the complexity of the dynamics due to bifurcations of homoclinic tangencies. Theorem \ref{theorem5} indicates the geometric complexity of the corresponding bifurcation sets in the Newhouse domain. Indeed, take any finite-dimensional analytic surface (in the space of dynamical systems) that intersects the Newhouse domain. We may consider this surface as the finite-parameter family of diffeomorphisms, so Theorem \ref{theorem5} gives us that a $C^r$-small perturbation of this surface creates a high-order tangency with an infinite-dimensional bifurcation set that corresponds to some homoclinic tangency. In \cite{Gonchenkocomplexity}, a similar statement for one-parameter families was proven for bifurcation sets corresponding to homoclinic loops of a saddle-focus equilibrium. 

In the smooth case $(r\le\infty)$, we have the following immediate corollary of Theorem~\ref{theorem5}: 
\begin{corollary}\label{corollary5.5}
In the space $\mathcal{PN}(L)$, there exists a dense subset with the following property: for every parametric family from this subset there is an $\ell$-dimensional ball $\widehat D^\ell\subset D^\ell$ such that $L$ has a homoclinic band for all $\zeta\in\widehat D^\ell$, which smoothly depends on $\zeta$. 
\end{corollary}

In a celebrated paper \cite{Berger2016}, Berger considers a particular open set $\mathcal {\mathring P}$ in the space of parametric families of diffeomorphisms (of arbitrarily high, but not infinite regularity) of any compact $d$-dimensional manifold with $d\ge 3$. He proves that the open set $\mathcal {\mathring P}$ contains a dense subset such that for every parametric family from this subset and for some finite covering of $D^\ell$ by open balls $\{\widehat D^\ell_i\}$ (independent of the family), for each ball $\widehat D^\ell_i$ and every $\zeta\in\widehat D^\ell_i$ there is an orbit of a corank-1 quadratic homoclinic tangency which smoothly depends on $\zeta$. The fact that non-splitting tangency happens for all $\zeta\in D^\ell$ is crucial for showing that a \textit{generic} family in $\mathcal {\mathring P}$ has infinitely many stable periodic orbits for all $\zeta\in D^\ell$ \cite{Berger2016}. 

The result of Corollary \ref{corollary5.5} that gives the existence of non-splitting tangencies only for a subset $\widehat D^\ell\subset D^\ell$ entails a weaker statement (we prove density, not genericity) yet for a wider class of parametric families, see Theorem \ref{theoremcoexistingattractors} below.\

\subsection{Universal dynamics}\label{subsection1.3}

\textbf{Universal diffeomorphisms.} Let $f:\mathcal M^d\rightarrow \mathcal M^d$ be a $C^r$ diffeomorphism, where $d\ge 2$ and $r=1,.\ldots,\infty$. Take a positive integer $q\le d$ and denote by $\mathbb D^{q}$ the closed unit disk in $\mathbb R^{q}$. Consider any small closed $q$-dimensional disk $B^{q}=\psi(\mathbb D^q)$, where $\psi$ is a $C^r$ embedding of $\mathbb R^q$ to $\mathcal M^d$. Suppose $f^k(B^q)\subset \psi(\mathbb R^q)$ for some positive $k$, so $B^q\cup f^k(B^q)$ is a $C^r$ submanifold of $\mathcal M^d$. By construction, the return map $f^k|_{B^q}$ is $C^r$-conjugate with the diffeomorphism $f_{k,\psi}\equiv\psi^{-1}\circ f^k\circ \psi: \mathbb D^q\rightarrow \mathbb R^q$. The diffeomorphisms $f_{k,\psi}$ are called \textit{renormalized iterations} of $f$; to stress the dependence on $q$, we also call them $q$-\textit{restricted renormalized iterations}. The set of all possible $q$-restricted renormalized iterations is called $q$-\textit{restricted renormalization class} of $f$. It characterizes the complexity of local dynamics of $f$ and is a generalisation of dynamical conjugacy class (see \cite{Turaev2}).
\begin{definition}\label{definitionuniversaldynamics}
Let $f:\mathcal M^d\rightarrow\mathcal M^d$ be a $C^r$ diffeomorphism and let $1\le p< d$. We say that $f$ is $(p,r)$-universal if there exists $q$ such that $p<q\le d$ and the $C^r$-closure of the $q$-restricted renormalization class of $f$ contains all maps $\mathbb D^q\to\mathbb R^q$ of the form
\begin{equation*}
X\mapsto \overline X: \quad
\begin{aligned} 
(\overline X_1,\ldots,\overline X_{p})&= \Phi(X_{1},\ldots,X_p), \\
(\overline X_{p+1},\ldots,\overline X_q)&=0,
\end{aligned}
\end{equation*}
where $\Phi$ is an arbitrary $C^r$ map of the closed unit ball $\mathbb D^p$ into $\mathbb R^p$. 
\end{definition}
The main novelty with respect to the original definitions \cite{Bonatti2003},\cite{Turaev2},\cite{Turaevintro2} is that now we allow the iterations to be restricted to the invariant submanifolds of lower dimension. This is done in order to incorporate the partially hyperbolic case where the diffeomorphism $f$ is uniformly expanding in certain directions. In this case, the forward invariant $d$-dimensional balls do not exist, so the restriction of iterations to $q$-dimensional disks becomes necessary. 

This creates a problem in the real-analytic case since normally-hyperbolic invariant manifolds of dimension $q<d$ typically have only finite regularity. Therefore, we also define the notion of "cross-universality" which essentially expresses the same idea without a reference to invariant submanifolds. This definition is more technical, and we formulate it in the Appendix, see Definition \ref{definitioncrossuniversaldynamics}. 

Let $\lambda$ and $\gamma$ stand for the absolute value of the leading stable and the leading unstable multipliers of a saddle periodic orbit $L$, respectively. Based on Theorems \ref{theorem1}--\ref{theorem4}, we prove the following
\begin{theorem}\label{theorem6}
Let $\lambda\gamma<1$. A $C^r$-generic diffeomorphism from $\mathcal N(L)$ is $(1,r)$-universal $(r=2,\ldots,\infty,\omega)$. If $L$ is a bi-focus, then a $C^r$-generic diffeomorphism from $\mathcal N(L)$ is $(2,r)$-universal $(r=1,\ldots,\infty,\omega)$.
\end{theorem}

For the smooth case, including $C^{\infty}$, this theorem understands universality in the sense of Definition \ref{definitionuniversaldynamics}. In the $C^{\omega}$ case, we, in general, talk about cross-universality. However, for the subclass of $C^{\omega}$ Newhouse regions $\mathcal N(L)$ for which $\dim W^u(L)=p$ (where $p=1$, or, in the bi-focus case, $p=2$), we still have universality with $q=d$, see Remark \ref{remarkcrossuniversaldynamics}. Note that, in general, for Theorem \ref{theorem6}, we take $q=p+\dim W^s(L)$ in the definition of universality.

For the case of two-dimensional diffeomorphisms, Theorem \ref{theorem6} was proved in \cite{arbitraryorder},\cite{arbitraryorderconservative}. For an open subset of the Barrientos-Raibekas-Asaoka domain (see \cite{Asaoka},\cite{Barrientos1},\cite{Barrientos3}), where a $2$-wild basic set is homoclinically related to a bi-focus, the $C^r$-genericity of $(2,r)$-universal $(r=2,\ldots,\infty,\omega)$ diffeomorphisms was proved in \cite{Mints}.

Let us mention other results on universal dynamics (see also reviews in \cite{Berger2026},\cite{Turaevintro1}). By \cite{Berger2026},\cite{Berger2025},\cite{Turaevintro1}, in the class of symplectic $C^{\infty}$ diffeomorphisms of the closed unit ball $\mathbb D^{d}$ into $\mathbb R^{d}$ $(d\in 2\mathbb N)$ with a non-degenerate elliptic periodic orbit, a generic diffeomorphism is universal (with a modification of the definition of the universality to the symplectic case \cite{Turaev2}). If $d=2$, the same holds true for symplectic $C^{\omega}$ diffeomorphisms \cite{Gelfreich2010},\cite{arbitraryorderconservative}. A non-symplectic reversible counterpart of this result (in every dimension $d\ge 2$) is in \cite{Gonchenkoonthreetypes}. In the space of two-dimensional diffeomorphisms, the $C^r$-genericity $(r=2,\ldots,\infty)$ of universal diffeomorphisms in the absolute Newhouse domain (the subset of the Newhouse domain where each diffeomorphism contains a wild basic set with area-contracting and area-expanding saddle periodic orbits) was proved in \cite{Turaevintro2}. Earlier, a similar result for the space of three-dimensional $C^1$-diffeomorphisms was obtained in \cite{Bonatti2003}. 

In these works, the universality property is weaker than in Theorem \ref{theorem6}. The renormalized iterations of the universal maps constructed previously approximate only diffeomorphisms isotopic to identity, while the universal maps of Theorem \ref{theorem6} approximate all the dynamics in the given dimension $p$, including non-invertible ones. The price is that we have to take $p<d$, the dimension of the ambient manifold\footnote{The reduction of the dimension can lead to a different type of strengthening of the notion of universality. Thus, by \cite{Asaokadegeneratebehavior},\cite{Asaokafastgrowth}, in an open subset of the Bonatti-D{\'\i}az domain (an open region corresponding to diffeomorphisms with heterodimensional cycles), $C^r$ diffeomorphisms $(r=3,\ldots,\infty)$ supporting universal semi-groups of interval diffeomorphisms are $C^r$-generic.}. Note that, in general, one cannot take $p>1$ (or, $p>2$ if $L$ is bi-focus) in Theorem \ref{theorem6}, as it follows from the dimension estimates of \cite{Turaev1}, i.e., this theorem is optimal in a sense. 

Another difference is that these dimension estimates prevent from the reduction to the near-identity case which was one of the necessary ingredients in the quoted works. Namely, the previous results often rely on the fact that (weakly) $C^r$-universal $d$-dimensional diffeomorphisms ($r=1,\ldots,\infty$ and $d\ge 2$) exist arbitrarily $C^r$-close to the identity map \cite{Turaevintro2}; for $r\le d$ this follows from \cite{Ruelleonthenature}, and for $p=d-1$ and any $r\ge 1$ - from \cite{Newhouseoccurrence}. A stronger statement was proven in \cite{Bergereverydiffeomorphism}: for any orientation preserving $C^r$ diffeomorphism $\phi$ of the unit ball $\mathbb D^d$ into itself, there exists a diffeomorphism $f$ such that it is $C^r$-close to the identity map and has a renormalized iteration that coincides with $\phi$ (not just approximates it). As we have mentioned, we cannot directly use this technique because of the dimension constraints. However, making use of Corollary \ref{corollary4.5}, we prove a similar in spirit result in Theorem \ref{theorem7} below (without the orientation preservation restriction). 

In this theorem, we take $p=1$ and $r=2,\ldots,\infty$, or, if $L$ is a bi-focus, we may also take $p=2$ and $r=1,\ldots,\infty$.
\begin{theorem}\label{theorem7}
For every $C^r$ diffeomorphism $\phi:\mathbb D^p\to \mathbb D^p$ there exists a dense subset $\mathcal A_{\phi}$ of a Newhouse region $\mathcal N(L)$ such that for every $f\in\mathcal A_{\phi}$, its $p$-restricted renormalized iteration $f_{k,\psi}$, for some $k>0$ and $\psi\in C^{\infty}$, equals to $\phi$.
\end{theorem}

\textbf{Universal parametric families.} Let $r=1,\ldots,\infty$. A family $(f_{\zeta})\in \mathcal P$ can be viewed as a $C^r$ diffeomorphism of $\mathcal M^{d+\ell}=\mathcal M^d\times D^\ell$ that acts as $(x,\zeta)\mapsto(f_{\zeta}(x),\zeta)$. Thus, given a positive integer $q\le d$, one defines $q$-restricted renormalized iterations by restriction to small $(q+\ell)$-dimensional disks and scaling the coordinates and parameters. Namely, we will scale $(x,\zeta)=\mathfrak h(X,Z)=(\psi_{Z}(X),\Psi(Z))$, where $\Psi$ is a $C^r$ embedding $\mathbb D^\ell \hookrightarrow D^\ell$, and $(\psi_Z)$ is a $C^r$ family of embeddings $\mathbb R^q\hookrightarrow\mathcal M^d$. Then, whenever there exists $k\ge 1$ such that $f_{\Psi(Z)}^k\circ \psi_{Z}(\mathbb D^q)\subset\psi_Z(\mathbb R^q)$ for all $Z\in \mathbb D^{\ell}$, the diffeomorphism 
\begin{equation}\label{universalparametricfamiliesrestrictedrenormalizediteration}
\begin{aligned} 
f_{k,\mathfrak h}:(X,Z)\mapsto(\psi_{Z}^{-1}\circ f_{\Psi(Z)}^k\circ \psi_{Z}|_{\mathbb D^q},Z)
\end{aligned}
\end{equation}
is called a $q$-restricted renormalized iteration $f_{k,\mathfrak h}$ for the family $(f_{\zeta})$. Note that $f_{k,\mathfrak h}$ acts from $\mathbb D^q\times \mathbb D^\ell$ to $\mathbb R^q\times \mathbb D^\ell$. One defines the $q$-restricted renormalization class as a set of all $q$-restricted renormalized iterations.   

\begin{definition}\label{definitionuniversaldynamicsparameters}
Let $(f_{\zeta})\in\mathcal P$ and let $p$ be a natural number such that $1\le p< d$. We say that $(f_{\zeta})$ is $(p,r)$-universal if the $C^r$-closure of its $q$-restricted renormalization class $($with $q>p)$ contains all $\ell$-parameter families of maps $\mathbb D^{q}\to\mathbb R^q$ of the form
\begin{equation*}
(X,Z)\mapsto \overline X: \quad
\begin{aligned} 
(\overline X_1,\ldots,\overline X_{p})&=\Phi(X_{1},\ldots,X_p,Z_1,\ldots,Z_\ell), \\
(\overline X_{p+1},\ldots,\overline X_q)&=0,\\
\end{aligned}
\end{equation*}
where $\Phi$ is an arbitrary $C^r$ map $\mathbb D^{p}\times \mathbb D^\ell\to\mathbb R^p$.
\end{definition}

From Theorems \ref{theorem5} and \ref{theorem6}, we infer the following
\begin{theorem}\label{theorem8}
Let $\lambda\gamma<1$. A $C^r$-generic family from $\mathcal{PN}(L)$ is $(1,r)$-universal $(r=2,\ldots,\infty,\omega)$. 
\end{theorem}
Here, as before, we distinguish between smooth and real-analytic cases. For $r\le\infty$, Theorem \ref{theorem8} addresses universality in accordance with Definition \ref{definitionuniversaldynamicsparameters}. For the real-analytic case, it refers to cross-universality (cf. Appendix \ref{sectionappendix}). 

We also obtain a parametric analogue of Theorem \ref{theorem7} for the case $p=1$. We take $r=2,\ldots,\infty$.
\begin{theorem}\label{theorem9}
For every $\ell$-parameter $C^r$ family $\varphi=(\phi_Z)$ of diffeomorphisms $\mathbb D^1\to\mathbb D^1$ there exists a dense subset $\mathcal A_{\varphi}$ of $\mathcal{PN}(L)$ such that for every family $(f_{\zeta})\in\mathcal A_{\varphi}$, its $1$-restricted renormalized iteration $f_{k,\mathfrak h}$, for some $k>0$ and $\mathfrak h\in C^{\infty}$, equals to $\varphi$.
\end{theorem}

\subsection{Coexistence of infinitely many attractors}\label{subsection1.4}
We say that a saddle periodic orbit $L$ is \textit{sectionally dissipative} if it has only one unstable multiplier and the product of any two different multipliers is less than one in the absolute value. This means that the derivative of the diffeomorphism at the points of $L$ contracts two-dimensional areas. Let $\mathcal{DN}(L)\subset\mathcal{N}(L)$ be an open region such that $L$ is sectionally dissipative. By Proposition \ref{propositionNewhousetheorem}, the regions $\mathcal{DN}(L)$ adjoin to any $C^r$ diffeomorphism $(r\ge 2)$ for which the sectionally dissipative orbit $L$ has a homoclinic tangency. 

One of the primary applications of the Newhouse Theorem is the genericity of the coexistence of infinitely many attractors. Thus, according to \cite{dynamicalphenomena},\cite{Newhouse1974},\cite{Palis1}, a generic diffeomorphism in $\mathcal{DN}(L)$ has infinitely many stable periodic orbits (we also mention papers \cite{BenedicksMartensPalmisano2018},\cite{GonchenkoShilnikov1986},\cite{GonchenkoShilnikov1987},\cite{MartensPalmisano2026}, where approaches different from that of \cite{Newhouse1974} were used for the analysis of this phenomenon). A generalization to the non-sectionally dissipative case was obtained in \cite{Gonchenko3}. The generic coexistence of infinitely many non-trivial attractors was shown in \cite{Gonchenko3} (closed invariant curves), \cite{Colli1998},\cite{Leal2008} (H{\'e}non-like attractors), \cite{GonchenkoOvsyannikov2010},\cite{GonchenkoShilnikovTuraev2009} (Lorenz-like attractors), \cite{Turaevintro2} (hyperbolic attractors), and infinitely many attractors and repellers in \cite{Gonchenkoonthreetypes},\cite{GonchenkoTuraevShilnikov1997}. 

The generic coexistence of infinitely many attractors puts in jeopardy the perspective of describing chaotic dynamics in finite terms. As a response to the challenge, the notion of measure-theoretical typicality was proposed to replace topological genericity, as a more adequate reflection of physical reality \cite{Hunt1992},\cite{Kaloshin1996},\cite{Palis2000},\cite{TedeschiniLalliYorke1986}. Indeed, establishing the typicality of infinitely many attractors is harder than establishing the genericity \cite{GorodetskiKaloshin2007},\cite{TedeschiniLalliYorke1986}. In \cite{Berger2016},\cite{BergerCrovisierPujals2022}, the measure-theoretical (Kolmogorov) typicality of the infinitude was proven for the smooth, not sectionally dissipative case (e.g., this excludes two-dimensional diffeomorphisms, and leaves the problem completely open in the analytic case). Theorem~\ref{theoremcoexistingattractors} below provides a partial advancement in the smooth, sectionally-dissipative case - this excludes strong prevalence of the finitude of attractors but does not preclude from the prevalence in the sense of \cite{Kaloshin1996}. 

We define a space $\mathcal{PDN}(L)\subset\mathcal P$ as follows:
\begin{equation*}
\begin{aligned} 
(f_{\zeta})\in\mathcal{PDN}(L)\;\;\; \text{if} \;\;\; f_{\zeta}\in\mathcal {DN}(L) \;\;\; \text{for some} \;\;\; \zeta\in D^\ell.
\end{aligned}
\end{equation*}
In the sectionally dissipative case, the normally-hyperbolic invariant interval $\psi_Z(\mathbb D^1)$ in formula \eqref{universalparametricfamiliesrestrictedrenormalizediteration} for the $1$-restricted renormalized iteration $f_{k,\mathfrak h}$ from Theorem \ref{theorem9} is normally-attracting. Therefore, by taking the diffeomorphism $\phi_Z$ such that it has infinitely many stable fixed points for all $Z$, e.g.,
\begin{equation*}
\phi_{Z}: \;\;\;
\begin{aligned} 
\overline X=X-\frac{1}{2}\cdot e^{-\frac{1}{X^2}}\cdot\sin\left(\frac{\pi}{X}\right) \;\; \text{for} \;\; X\in[-1;1], \\
\end{aligned}
\end{equation*}
Theorem \ref{theorem9} immediately gives us
\begin{theorem}\label{theoremcoexistingattractors}
In the space $\mathcal{PDN}(L)$, there exists a dense subset with the following property: for every parametric family $(f_{\zeta})$ from this subset there is an $\ell$-dimensional ball $\widehat D^\ell\subset D^\ell$ such that $f_{\zeta}$ has infinitely many stable periodic orbits for all $\zeta\in\widehat D^\ell$. 
\end{theorem}

\section{Partially hyperbolic structure}\label{section2}

\subsection{Local map}\label{subsectionlocalmap}
Let $f:\mathcal M^d\to\mathcal M^d$ be a $C^r$ diffeomorphism possessing a hyperbolic saddle periodic orbit $L=\{O,f(O),\ldots,f^{b-1}(O)\}$ of period $b$. Let $\lambda_1,\ldots,\lambda_{d_s}$ and $\gamma_1,\ldots,\gamma_{d_u}$ be the stable and unstable multipliers of $L$, respectively. We distinguish between leading and non-leading multipliers as follows:
\begin{itemize}
    \item $\lambda_1, \dots, \lambda_{m_s}$ are the leading stable and $\lambda_{m_s+1},\dots,\lambda_{d_s}$ are the non-leading stable multipliers: $1>|\lambda_1|=\ldots=|\lambda_{m_s}|>|\lambda_{m_s+1}|\ge\ldots\ge|\lambda_{d_s}|$.
    \item $\gamma_1, \dots, \gamma_{m_u}$ are the leading unstable and $\gamma_{m_u+1},\dots,\gamma_{d_u}$ are the non-leading unstable multipliers: $1<|\gamma_1|=\ldots=|\gamma_{m_u}|<|\gamma_{m_u+1}|\le\ldots\le|\gamma_{d_u}|$.
\end{itemize}
Let $\lambda\equiv|\lambda_1|$ and $\gamma\equiv|\gamma_1|$. Let $\text{alg}(\cdot)$ and $\text{geom}(\cdot)$ denote the algebraic and geometric multiplicities of the multipliers, respectively. We assume 
\textbf{\hypertarget{C1.}}

\begin{doublespace}
\textbf{\hypertarget{C1.1.}{C1.1.}} $\text{alg}(\lambda_i)=\text{geom}(\lambda_i)$ for all $i\in\{1,\ldots,m_s\}$;

\textbf{\hypertarget{C1.2.}{C1.2.}} $\text{alg}(\gamma_j)=\text{geom}(\gamma_j)$ for all $j\in\{1,\ldots,m_u\}$.
\end{doublespace}

Condition \textbf{\hyperlink{C1.}{C1}} is dense in the set of diffeomorphisms that are $C^r$-close to $f$. Indeed, by applying an arbitrarily $C^r$-small perturbation that splits the equal multipliers, we arrive to one of the four main cases:
\begin{itemize}
\item $(1,1):$ \; $\lambda_1$ and $\gamma_1$ are real, and $\lambda>|\lambda_2|$, $\gamma<|\gamma_2|$;

\item $(2,1):$ \; $\lambda_{1,2}=\lambda\cdot e^{\pm i\varphi}$, $\gamma_1$ is real, and $\lambda>|\lambda_3|$, $\gamma<|\gamma_2|$;

\item $(1,2):$ \; $\lambda_1$ is real, $\gamma_{1,2}=\gamma\cdot e^{\pm i\psi}$, and $\lambda>|\lambda_2|$, $\gamma<|\gamma_3|$;

\item $(2,2):$ \; $\lambda_{1,2}=\lambda\cdot e^{\pm i\varphi}$, $\gamma_{1,2}=\gamma\cdot e^{\pm i\psi}$, and $\lambda>|\lambda_3|$, $\gamma<|\gamma_3|$,
\end{itemize}
where $\varphi,\psi\in(0,\pi)$. Note that the union of these four classes is  open and dense among $C^r$ diffeomorphisms possessing a periodic orbit $L$.

Let $f$ be embedded into an $\ell$-parameter $(\ell\ge 1)$ $C^r$ family of diffeomorphisms $(f_{\varepsilon})$, where the parameters $\varepsilon$ vary in the $\ell$-dimensional closed ball $D^\ell$ centered at zero, with $f_0\equiv f$. Let $T_{O}$ be the restriction of the period map $f_{\varepsilon}^b$ to a small neighborhood $U(O)$. We refer to $T_{O}$ as the \textit{local map}. By centering the origin at $O$ and straightening the local stable and local unstable manifolds, we represent $T_{O}$ for all small $\varepsilon$ as
\begin{equation}\label{eqlocal1}
\begin{aligned} 
\overline x&= \mathcal A_1(\varepsilon) \cdot x+\vartheta_1(x,u,y,v,\varepsilon), \quad &\overline y&=\mathcal B_1(\varepsilon) \cdot y+\vartheta_3(x,u,y,v,\varepsilon),\\ 
\overline u&= \mathcal A_2(\varepsilon) \cdot u+\vartheta_2(x,u,y,v,\varepsilon), \quad &\overline v&=\mathcal B_2(\varepsilon) \cdot v+\vartheta_4(x,u,y,v,\varepsilon),\\
\end{aligned}
\end{equation}
where $(x,u,y,v)\in\mathbb R^{m_s}\times\mathbb R^{d_s-m_s}\times\mathbb R^{m_u}\times\mathbb R^{d_u-m_u}$. The matrices $\mathcal A_1(0)$ and $\mathcal A_2(0)$ have eigenvalues corresponding to the leading stable and non-leading stable multipliers of $L$ at $\varepsilon=0$, respectively. Similarly, $\mathcal A_3(0)$ and $\mathcal A_4(0)$ have eigenvalues corresponding to the leading unstable and non-leading unstable multipliers of $L$ at $\varepsilon=0$, respectively. The nonlinearities $\vartheta$ satisfy 
\begin{equation}\label{eqlocal2.0}
\begin{aligned} 
&\vartheta_{1,2}(0,0,y,v,\varepsilon)\equiv0, \qquad \vartheta_{3,4}(x,u,0,0,\varepsilon)\equiv0.\\
\end{aligned} 
\end{equation}
Now, let $r=2,\ldots,\infty,\omega$. Then, following \cite{Gonchenko3}, we apply $C^r$ coordinate transformation to eliminate nonlinear terms in \eqref{eqlocal1} that cause the iterations $T_{O}^k$ to deviate significantly from the iterations of the linearization. After this transformation, the nonlinearities $\vartheta$ satisfy the additional conditions:
\begin{equation}\label{eqlocal2.1}
\begin{aligned} 
\vartheta_{1}(x,u,0,0,\varepsilon)&\equiv0, \qquad \vartheta_{3}(0,0,y,v,\varepsilon)\equiv0,\\
\end{aligned} 
\end{equation}
\begin{equation}\label{eqlocal2.2}
\begin{aligned} 
\frac{\partial \vartheta_{1,2}}{\partial x}(0,0,y,v,\varepsilon)&\equiv0, \qquad \frac{\partial \vartheta_{3,4}}{\partial y}(x,u,0,0,\varepsilon)\equiv0.\\
\end{aligned} 
\end{equation}
When expressed in these coordinates, the map $T_{O}$ is said to be in the \textit{main normal form}. Let us clarify that the smoothness of the coordinate transformation that brings $T_O$ to the main normal form is, in general, only $C^{r-2}$ with respect to parameters. Namely, the transformation and its derivatives up to the second order are $C^{r-2}$ jointly with respect to coordinates and parameters (see \cite{Gonchenko3}).
\begin{remark}
In general, we do not assume that the diagonalization condition \textbf{\hyperlink{C1.}{C1}} stays for $\varepsilon\not=0$. However, in the four main cases mentioned above, the diagonalizability holds for all small $\varepsilon$. Thus, one can apply a $C^r$ change of coordinates to reduce the matrices $\mathcal A_1(\varepsilon)$ and $\mathcal B_1(\varepsilon)$ to the real Jordan form while keeping $T_{O}$ in the main normal form. Then, for all small $\varepsilon$, these matrices are given by: 
\begin{equation*}
\begin{aligned} 
&\mathcal A_1(\varepsilon)=\pm\lambda_{\varepsilon} \;\; (\text{cases} \;\; (1,1), (1,2)), \;\; \text{and} \;\; \mathcal A_1(\varepsilon)=\lambda_{\varepsilon}\cdot\begin{pmatrix}\cos\varphi_{\varepsilon}& -\sin\varphi_{\varepsilon}\\ \sin\varphi_{\varepsilon}& \cos\varphi_{\varepsilon} \end{pmatrix} \;\; (\text{cases} \;\;  (2,1), (2,2));
\end{aligned}
\end{equation*}
\begin{equation*}
\begin{aligned} 
&\mathcal B_1(\varepsilon)=\pm\gamma_{\varepsilon} \;\; (\text{cases} \;\; (1,1), (2,1)), \;\; \text{and} \;\; \mathcal B_1(\varepsilon)=\gamma_{\varepsilon}\cdot\begin{pmatrix}\cos\psi_{\varepsilon}& -\sin\psi_{\varepsilon}\\ \sin\psi_{\varepsilon}& \cos\psi_{\varepsilon} \end{pmatrix} \;\; (\text{cases} \;\;  (1,2), (2,2)).
\end{aligned}
\end{equation*}
In this coordinate system, $T_{O}$ is said to be in the \textit{special main normal form}.
\end{remark}

Let $(\tilde x,\tilde u,\tilde y,\tilde v)=T_{O}^k(x,u,y,v)$. According to \cite{Shilnikov2}, the coordinates $(\tilde x, \tilde u, y, v)$ can be uniquely expressed in terms of $(x,u,\tilde y,\tilde v)$; that is, the map $T_{O}^k$ can be written in \textit{cross form}. We will employ this form in all subsequent proofs. When $T_{O}$ is in the main normal form, the relation between $(\tilde x,\tilde u,y,v)$ and $(x,u,\tilde y,\tilde v)$ is given by  
\begin{proposition}\label{propositionlocalmap}(\cite{Gonchenko3}, Lemma 7).
Let the local map $T_{O}$ be brought to the main normal form. Then, for all small $\varepsilon$ and all large $k$, the following relations hold:
\begin{equation}\label{eqlocal3}
\begin{aligned} 
\tilde x&= \mathcal A_1^k(\varepsilon) \cdot x+\hat\lambda^k\cdot\hat\xi_{1,k}\left(x,u,\tilde y,\tilde v,\varepsilon\right), &y&= \mathcal B_1^{-k}(\varepsilon) \cdot \tilde y+\hat\gamma^{-k}\cdot\hat\xi_{3,k}\left(x,u,\tilde y,\tilde v,\varepsilon\right),\\
\tilde u&=\hat\lambda^k\cdot\hat\xi_{2,k}\left(x,u,\tilde y,\tilde v,\varepsilon\right), &v&=\hat\gamma^{-k}\cdot\hat\xi_{4,k}\left(x,u,\tilde y,\tilde v,\varepsilon\right),
\end{aligned}
\end{equation}
where $\hat\lambda$ and $\hat\gamma$ are some constants such that $0<\hat\lambda<\lambda$, $\hat\gamma>\gamma$. The functions $\hat\xi_{k}$ and their derivatives up to order $r-2$ tend to zero uniformly as $k\to+\infty$.
\end{proposition}

For most of the paper, we will work with particular sequences of $k$ values and also restrict the range of the values of $\varepsilon$ such that the formula \eqref{eqlocal3} gets simplified as in Corollary \ref{corollarylocalmap} below.   
\begin{definition}\label{definitiontamesequence}
Given a sequence $\langle k_\alpha\rangle_{\alpha\in\mathbb N}$ of natural numbers with $k_\alpha\xrightarrow{}+\infty$, a sequence $\langle\mathfrak{D}_{\alpha}\rangle_{\alpha\in\mathbb N}$ of open balls centered at $\varepsilon=0$ in the $\varepsilon$-parameter space is called $\langle k_\alpha\rangle$-tame if the radii of the balls satisfy
\begin{equation}\label{lemmasequenceepseq1}
\begin{aligned} 
\rho\left(\mathfrak{D}_{\alpha}\right)&=o\left(\frac{1}{k_\alpha}\right).
\end{aligned}
\end{equation}
\end{definition}

\begin{lemma}\label{lemmasequence}
Let the periodic orbit $L$ satisfy condition \textbf{\hyperlink{C1.}{C1}} at $\varepsilon=0$. Then there exists a sequence $\langle k_\alpha\rangle_{\alpha\in\mathbb N}$ such that for any $\langle k_{\alpha}\rangle$-tame sequence $\langle\mathfrak{D}_{\alpha}\rangle_{\alpha\in\mathbb N}$ and all $\varepsilon\in\mathfrak{D}_{\alpha}$, one has
\begin{equation}\label{lemmasequenceepseq2}
\begin{aligned} 
&\mathcal A_1^{k_{\alpha}}(\varepsilon)=\lambda^{k_{\alpha}}\cdot I+o_{r-2}(\lambda^{k_{\alpha}}) \quad \text{and} \quad \mathcal B_1^{-k_{\alpha}}(\varepsilon)=\gamma^{-k_{\alpha}}\cdot I+o_{r-2}(\gamma^{-k_{\alpha}}),
\end{aligned}
\end{equation}
where $o_{r-2}$ means that the limit is taken with $r-2$ derivatives with respect to $\varepsilon$. 
\end{lemma}

\begin{demo}
\textit{Step 1.} By condition \textbf{\hyperlink{C1.}{C1}}, we can assume that the matrices $\mathcal A_1(0)$ and $\mathcal B_1(0)$ are in the real Jordan form:   
\begin{equation}\label{lemmasequenceeq1}
\begin{aligned} 
\mathcal{A}_1(0)&=\lambda\cdot\text{diag}\left(S_1,\ldots,S_{n_s}\right)\quad \text{and} \quad \mathcal{B}_1(0)=\gamma\cdot\text{diag}\left(U_{1},\ldots,U_{n_u}\right),
\end{aligned}
\end{equation}
where each block $S_{i}$ and $U_{j}$ is either a number $\pm 1$ or a $2\times2$ rotation matrix. To each block $S_i$, we associate a complex argument $\varphi_i\in[0,\pi]$. Specifically, $S_i=1$ corresponds to $\varphi_i=0$ and $S_i=-1$ corresponds to $\varphi_i=\pi$. In the case of a $2\times2$ rotation matrix, the argument $\varphi_i$ is the angle of rotation. Similarly, to each block $U_j$, we associate a complex argument $\psi_j\in[0,\pi]$. Denote by $\overline\varphi$ the vector $(\varphi_1,\ldots,\varphi_{n_s})$ and by $\overline\psi$ the vector $(\psi_1,\ldots,\psi_{n_u})$.

Let $\mathbb T^{\,n_s+n_u}$ be the $(n_s+n_u)$-dimensional torus and let $\Phi: \mathbb T^{\,n_s+n_u}\to\mathbb T^{\,n_s+n_u}$ be the translation 
\begin{equation*}
\begin{aligned} 
w_1&=w_1+\overline\varphi, \quad w_2=w_2-\overline\psi \quad \left(mod\; 2\pi\right). 
\end{aligned}
\end{equation*}
Since the torus $\mathbb T^{\,n_s+n_u}$ is compact, there is a point $(w_1^*,w_2^*)\in\mathbb T^{\,n_s+n_u}$ which is recurrent for $\Phi$. This means there exists a sequence $\langle k_\alpha\rangle_{\alpha\in\mathbb N}$ such that
\begin{equation}\label{lemmasequenceeq2}
\begin{aligned} 
\Phi^{k_\alpha}(w_1^*,w_2^*)= (w_1^*,w_2^*)+o_{r-1}(1)_{\alpha\rightarrow+\infty} \quad \left(mod\; 2\pi\right).
\end{aligned}
\end{equation}
Let $\Psi$ be the translation of the torus $\mathbb T^{\,n_s+n_u}$ which takes $(0,0)$ to $(w_1^*,w_2^*)$. The maps $\Phi$ and $\Psi$ commute; therefore, $\Phi^{k_{\alpha}}=\Psi^{-1}\circ\Phi^{k_{\alpha}}\circ\Psi$. This implies 
\begin{equation}\label{lemmasequenceeq3}
\begin{aligned} 
\Phi^{k_{\alpha}}\left(0,0\right)=\Psi^{-1}\circ\Phi^{k_{\alpha}}\circ\Psi\left(0,0\right)=\Psi^{-1}\circ\Phi^{k_{\alpha}}(w_1^*,w_2^*).
\end{aligned}
\end{equation}
By \eqref{lemmasequenceeq2} and \eqref{lemmasequenceeq3},
\begin{equation}\label{lemmasequenceeq4}
\begin{aligned} 
\Phi^{k_{\alpha}}(0,0)=(0,0)+o_{r-1}(1)_{\alpha\rightarrow+\infty} \quad \left(mod\; 2\pi\right). 
\end{aligned}
\end{equation}

By \eqref{lemmasequenceeq4}, each $S^{k_{\alpha}}_{i}$ and $U^{-k_{\alpha}}_{j}$ approach the identity matrix. By \eqref{lemmasequenceeq1}, 
\begin{equation}\label{lemmasequenceeq4.5}
\begin{aligned} 
\mathcal{A}^{k_{\alpha}}_1(0)&=\lambda^{k_{\alpha}}\cdot\text{diag}\left(S^{k_{\alpha}}_1,\ldots,S^{k_{\alpha}}_{n_s}\right),\\ 
\mathcal{B}^{-k_{\alpha}}_1(0)&=\gamma^{-k_{\alpha}}\cdot\text{diag}\left(U^{-k_{\alpha}}_{1},\ldots,U^{-k_{\alpha}}_{n_u}\right).
\end{aligned}
\end{equation}
Thus, 
\begin{equation}\label{lemmasequenceeq5}
\begin{aligned} 
&\mathcal{A}^{k_{\alpha}}_1(0)=\lambda^{k_{\alpha}}\cdot I+o_{r-1}(\lambda^{k_{\alpha}}) \quad \text{and} \quad \mathcal{B}^{-k_{\alpha}}_1(0)=\gamma^{-k_{\alpha}}\cdot I+o_{r-1}(\gamma^{-k_{\alpha}}).
\end{aligned}
\end{equation}

\textit{Step 2.} Let $\langle\mathfrak{D}_{\alpha}\rangle_{\alpha\in\mathbb N}$ be a $\langle k_{\alpha}\rangle$-tame sequence. For $\varepsilon\in\mathfrak{D}_{\alpha}$, we write $\mathcal{A}_1(\varepsilon)$ as
\begin{equation*}
\begin{aligned} 
&\mathcal{A}_1(\varepsilon)=\mathcal{A}_1(0)+\mathcal R(\varepsilon),
\end{aligned}
\end{equation*}
where the remainder $\mathcal R(\varepsilon)$ satisfies $||\mathcal R(\varepsilon)||_{r-2}=O_{r-2}(||\varepsilon||)$. It is directly verified that the following relation holds:
\begin{equation}\label{lemmasequenceeq7}
\begin{aligned} 
&\mathcal{A}^{k_{\alpha}}_1(\varepsilon)-\mathcal{A}^{k_{\alpha}}_1(0)=\sum\limits_{i=0}^{k_{\alpha}-1}\mathcal A_1^{k_{\alpha}-i-1}(\varepsilon)\cdot\mathcal R(\varepsilon)\cdot\mathcal A_1^{i}(0).
\end{aligned}
\end{equation}
Let us estimate the $(r-2)$-norm of the left-hand side in \eqref{lemmasequenceeq7}. By \eqref{lemmasequenceeq4.5}, for all $i\in\{0,\ldots,k_{\alpha}-1\}$,
\begin{equation}\label{lemmasequenceeq8}
\begin{aligned} 
&||\mathcal A_1^{i}(0)||_{r-1}=\lambda^{i}.
\end{aligned}
\end{equation}
Condition \eqref{lemmasequenceepseq1} implies that for $\varepsilon\in\mathfrak{D}_{\alpha}$,
\begin{equation}\label{lemmasequenceeq9}
\begin{aligned} 
&||\mathcal R(\varepsilon)||_{r-2}=O_{r-2}(||\varepsilon||)=o_{r-2}\left(\frac{1}{k_{\alpha}}\right).
\end{aligned}
\end{equation}
The $C^{r-2}$-norm of $\mathcal{A}_1(\varepsilon)$ can be estimated as
\begin{equation*}
\begin{aligned} 
&||\mathcal{A}_1(\varepsilon)||_{r-2}=||\mathcal{A}_1(0)+\mathcal R(\varepsilon)||_{r-2}\le||\mathcal{A}_1(0)||_{r-2}+||\mathcal R(\varepsilon)||_{r-2}=\lambda+o_{r-2}\left(\frac{1}{k_{\alpha}}\right),
\end{aligned}
\end{equation*}
therefore, for all $i\in\{0,\ldots,k_{\alpha}-1\}$,
\begin{equation}\label{lemmasequenceeq10}
\begin{aligned} 
||\mathcal A_1(\varepsilon)||^{k_{\alpha}-i-1}_{r-2}\le&\left(\lambda+o_{r-2}\left(\frac{1}{k_{\alpha}}\right)\right)^{k_{\alpha}-i-1}=\lambda^{k_{\alpha}-i-1}\cdot\left(1+\frac{1}{\lambda}\cdot o_{r-2}\left(\frac{1}{k_{\alpha}}\right)\right)^{k_{\alpha}-i-1}=\\
=&\lambda^{k_{\alpha}-i-1}+o_{r-2}(\lambda^{k_{\alpha}-i-1}).
\end{aligned}
\end{equation}
Combining \eqref{lemmasequenceeq8}--\eqref{lemmasequenceeq10}, we estimate 
\eqref{lemmasequenceeq7} as 
\begin{equation*}
\begin{aligned} 
||\mathcal{A}^{k_{\alpha}}_1(\varepsilon)-\mathcal{A}^{k_{\alpha}}_1(0)||_{r-2}\le&\sum\limits_{i=0}^{k_{\alpha}-1}||\mathcal A_1^{k_{\alpha}-i-1}(\varepsilon)||_{r-2}\cdot||\mathcal R(\varepsilon)||_{r-2}\cdot||\mathcal A_1^{i}(0)||_{r-2}=\\
=&\sum\limits_{i=0}^{k_{\alpha}-1}(\lambda^{k_{\alpha}-i-1}+o_{r-2}(\lambda^{k_{\alpha}-i-1}))\cdot o_{r-2}\left(\frac{1}{k_{\alpha}}\right)\cdot\lambda^{i}=\\
=&(\lambda^{k_{\alpha}-1}+o_{r-2}(\lambda^{k_{\alpha}-1}))\cdot o_{r-2}(1).
\end{aligned}
\end{equation*}
By \eqref{lemmasequenceeq5}, this gives the result \eqref{lemmasequenceepseq2} for $\mathcal{A}^{k_{\alpha}}_1(\varepsilon)$. The estimate for $\mathcal{B}^{-k_{\alpha}}_1(\varepsilon)$ is done in the same way.
\end{demo}

\begin{corollary}\label{corollarylocalmap}
Let the local map $T_{O}$ be brought to the main normal form. Let the periodic orbit $L$ satisfy condition \textbf{\hyperlink{C1.}{C1}} at $\varepsilon=0$. Then there exists a sequence $\langle k_\alpha\rangle_{\alpha\in\mathbb N}$ such that for any $\langle k_{\alpha}\rangle$-tame sequence $\langle\mathfrak{D}_{\alpha} \rangle_{\alpha\in\mathbb N}$ and all $\varepsilon\in\mathfrak{D}_{\alpha}$, one has
\begin{equation}\label{eqlocal3tamesequence}
\begin{aligned} 
\tilde x&= \lambda^{k_{\alpha}}\cdot (x+\xi_{1,k_\alpha}(x,u,\tilde y,\tilde v,\varepsilon)), &y&= \gamma^{-k_{\alpha}} \cdot (\tilde y+\xi_{3,k_\alpha}(x,u,\tilde y,\tilde v,\varepsilon)),\\
\tilde u&=\lambda^{k_{\alpha}}\cdot\xi_{2,k_\alpha}(x,u,\tilde y,\tilde v,\varepsilon), &v&=\gamma^{-k_{\alpha}}\cdot\xi_{4,k_\alpha}(x,u,\tilde y,\tilde v,\varepsilon),
\end{aligned}
\end{equation}
where the functions $\xi_{k_\alpha}$ and their derivatives up to order $r-2$ tend to zero uniformly as $\alpha\rightarrow+\infty$.
\end{corollary}

\subsection{Global map}\label{subsectionglobalmap}
Let the diffeomorphism $f$ have a homoclinic orbit $\Gamma$ which is contained in the intersection of the stable manifold $W^s(L)$ and the unstable manifold $W^u(L)$. Note that the periodic orbit $L$ can be viewed as a partially hyperbolic set with $(m_s+m_u)$-dimensional center. Below we describe conditions which ensure that the corresponding partially hyperbolic structure is extended to $cl(\Gamma)=L\cup\Gamma$.

Since the orbit $\Gamma$ is homoclinic, we can locate points $M^+$ and $M^-$ within this orbit such that $M^+\in W^s_{loc}(O)$, $M^-\in W^u_{loc}(O)$, and $M^+=f^{\kappa}(M^-)$ for some $\kappa>0$. We take small neighborhoods $\Pi^+$ of $M^+$ and $\Pi^-$ of $M^-$. We define the \textit{global map} $T_{\Gamma}$ as $T_{\Gamma}\equiv f^{\kappa}:\Pi^-\to\Pi^+$.

We denote the invariant subspaces of the linearization of the local map $T_{O}$ corresponding to the non-leading stable and non-leading unstable multipliers as $\mathcal E^{ss}$ and $\mathcal E^{uu}$, respectively. It is well known (see, e.g., \cite{Shilnikov1}) that there exist $C^r$-smooth invariant manifolds $W^{ss}_{loc}(O)\subset W^s_{loc}(O)$ and $W^{uu}_{loc}(O)\subset W^u_{loc}(O)$ tangent to $\mathcal E^{ss}$ and $\mathcal E^{uu}$ at $O$. These manifolds are embedded in invariant $C^r$-smooth foliations $\mathcal L^{ss}$ and $\mathcal L^{uu}$, which are uniquely defined on $W^s_{loc}(O)$ and $W^u_{loc}(O)$. 

Furthermore, let $\mathcal E^{se}$ be the eigensubspace of the linear part of $T_{O}$ corresponding to the stable and leading unstable multipliers, and let $\mathcal E^{ue}$ be the eigensubspace corresponding to the unstable and leading stable multipliers. There are invariant manifolds $W^{se}_{loc}(O)$ and $W^{ue}_{loc}(O)$ that are tangent to $\mathcal E^{se}$ and $\mathcal E^{ue}$ at $O$ and are at least $C^1$-smooth. While these manifolds are not unique, the tangent spaces $\mathcal T_{M^+}(W^{se}_{loc}(O))$ and $\mathcal T_{M^-}(W^{ue}_{loc}(O))$ are uniquely defined. Now we can formulate the partial-hyperbolicity conditions: \textbf{\hypertarget{C2.}}

\begin{doublespace}

\textbf{\hypertarget{C2.1.}{C2.1.}} $T_{\Gamma}(W^{ue}_{loc}(O))$ is transverse to the leaf $l^{ss}$ of $\mathcal L^{ss}$ passing through $M^+$;

\textbf{\hypertarget{C2.2.}{C2.2.}} $W^{se}_{loc}(O)$ is transverse to $T_{\Gamma}(l^{uu})$, where $l^{uu}$ is the leaf of $\mathcal L^{uu}$ passing through $M^-$.
\end{doublespace}

Conditions \textbf{\hyperlink{C2.1.}{C2.1}} and \textbf{\hyperlink{C2.2.}{C2.2}} are independent of the choice of homoclinic points $M^+$ and $M^-$. They are open in the class of systems exhibiting homoclinic tangency and, if not satisfied, can be achieved via an arbitrarily small $C^r$-perturbation of $f$ (see \cite{Turaev1}). 

Let $(x,u,y,v)$ be coordinates such that the local map $T_{O}$ is in the main normal form. In these coordinates, by \eqref{eqlocal2.0}, the manifolds $W^s_{loc}(O)$ and $W^u_{loc}(O)$ are straightened and given by 
\begin{equation}\label{foliationseq1}
\begin{aligned} 
&W^s_{loc}(O):\;(y,v)=0, \quad\quad W^u_{loc}(O):\;(x,u)=0.
\end{aligned}
\end{equation}
By \eqref{eqlocal2.1}, the manifolds $W^{ss}_{loc}(O)$, $W^{uu}_{loc}(O)$ are also straightened and can be written as  
\begin{equation}\label{foliationseq2}
\begin{aligned} 
&W^{ss}_{loc}(O):\;(x,y,v)=0, \quad\quad W^{uu}_{loc}(O):\;(x,y,u)=0.
\end{aligned}
\end{equation}
The corresponding foliations $\mathcal L^{ss}$ and $\mathcal L^{uu}$ are given by
\begin{equation}\label{foliationseq3}
\begin{aligned} 
&\mathcal L^{ss}:\; x=const,\;(y,v)=0, \quad\quad \mathcal L^{uu}:\;y=const,\;(x,u)=0.
\end{aligned}
\end{equation}
Furthermore, the equalities \eqref{eqlocal2.1} imply  
\begin{equation}\label{foliationseq4}
\begin{aligned} 
&\mathcal T_{M^+}(W^{se}_{loc}(O)):\;v=0, \quad\quad \mathcal T_{M^-}(W^{ue}_{loc}(O)):\;u=0.
\end{aligned}
\end{equation}

If the diffeomorphism $f$ is embedded into a finite-parameter family $(f_{\varepsilon})$ (with $f_0\equiv f$), then for all small $\varepsilon$, formulas \eqref{foliationseq1}--\eqref{foliationseq4} remain valid (we take the points $M^+$ and $M^-$ independent of $\varepsilon$). The global map $T_{\Gamma}:\Pi^-\to\Pi^+$ stay well-defined for all small $\varepsilon$, even though the homoclinic orbit $\Gamma$ may disappear if it is non-transverse. We provide explicit formulas for the global map in Section \ref{subsectionformulasforglobalmap} (see Lemma \ref{lemmaglobalmap}).

\section{Choice of coordinates}\label{sectionchoiceofcoordinates}

\subsection{Coordinate definitions of corank and order of tangency}\label{subsectioncorankandorder}

Let $\mathcal M^d$ be a $C^r$-smooth $d$-dimensional manifold, where $d\ge 2$ and $r=1,\ldots,\infty,\omega$. Let $\mathcal W_1$ and $\mathcal W_2$ be submanifolds of $\mathcal M^d$ with complementary dimensions $d_1$ and $d_2$, which are tangent at a point $M$. Near this point, we introduce coordinates $z\in\mathbb R^d$ such that the manifolds $\mathcal W_1$ and $\mathcal W_2$ are given by 
\begin{equation*}
\begin{aligned} 
&\mathcal W_1:\quad 0=F_1(z) \quad\quad \text{and} \quad\quad \mathcal W_2:\quad 0=F_2(z),
\end{aligned}
\end{equation*}
where $F_1:\mathbb R^{d}\to\mathbb R^{d_2}$ and $F_2:\mathbb R^{d}\to\mathbb R^{d_1}$ are smooth functions such that $\text{rank}\left.\frac{\partial F_1}{\partial z}\right|_{M}=d_2$ and $\text{rank}\left.\frac{\partial F_2}{\partial z}\right|_{M}=d_1$. The tangent spaces $\mathcal T_M \mathcal W_1$ and $\mathcal T_M \mathcal W_2$ are given by 
\begin{equation*}
\begin{aligned} 
&\mathcal T_M \mathcal W_1:\quad 0=\left.\frac{\partial F_1}{\partial z}\right|_{M}\cdot z \quad\quad \text{and} \quad\quad \mathcal T_M \mathcal W_2:\quad 0=\left.\frac{\partial F_2}{\partial z}\right|_{M}\cdot z.
\end{aligned}
\end{equation*}
Define the $d\times d$ matrix $\mathcal C$ as
\begin{equation*}
\begin{aligned} 
&\mathcal C=\begin{pmatrix} \frac{\partial F_1}{\partial z}|_{M} \\ \frac{\partial F_2}{\partial z}|_{M} \end{pmatrix}.
\end{aligned}
\end{equation*}
We can define the corank of tangency as follows:
\begin{definition}\label{coordinatedefinitioncorank}
The tangency between $\mathcal W_1$ and $\mathcal W_2$ is of corank $c$ if
\begin{equation}\label{corankoftangencyeq1}
\begin{aligned} 
&d-\text{rank}\;\mathcal C=\text{corank}\;\mathcal C=c.
\end{aligned}
\end{equation}
\end{definition}
\begin{remark}
When \eqref{corankoftangencyeq1} holds, we say that the system of equations $\{0=F_1(z), 0=F_2(z)\}$ has a corank-$c$ solution $($the point $M$$)$.
\end{remark}
The equivalence of Definition \ref{coordinatedefinitioncorank} to the definition of corank in Section \ref{subsection1.1} follows directly from \eqref{corankoftangency}. 

Now, suppose the tangency between $\mathcal W_1$ and $\mathcal W_2$ at $M$ is of corank $c$. Following \cite{bifurcationsandstability} (section II.6.), there exist a coordinate system $z=(z_1,z_2,z_3,z_4)\in\mathbb R^c\times\mathbb R^{d_1-c}\times\mathbb R^c\times\mathbb R^{d_2-c}$ near $M$ such that $\mathcal W_1$ and $\mathcal W_2$ are represented as 
\begin{equation}\label{orderoftangencyeq1}
\begin{aligned} 
&\mathcal W_1:\;(z_3,z_4)=(0,0), \quad\quad \mathcal W_2:\; z_3=F(z_1)+\theta_{1}(z_1,z_4), \;\; z_2=\theta_{2}(z_1,z_4),
\end{aligned}
\end{equation}
where $F$ and $\theta$ are $C^r$-functions satisfying 
\begin{equation}\label{orderoftangencyeq2}
\begin{aligned} 
F(0)&=0, &\left.\frac{\partial F}{\partial z_1}\right|_{z_1=0}&=0,\quad \text{and} \quad \theta_1(z_1,0)=0.
\end{aligned}
\end{equation}

Using this representation, we can formulate the coordinate-based definition of the order of tangency. 
\begin{definition}\label{coordinatedefinitionorder}
A corank-$c$ tangency is of order $n\le r-1$ $(r=2,\ldots,\infty,\omega)$ if
\begin{equation}\label{coordinatedefinitionordereq}
\begin{aligned} 
&\left.\frac{\partial^{i} F}{\partial z_1^{i}}\right|_{z_1=0}=0 \;\;\; \text{for} \;\;\; 0\le |i|\le n, \;\; \text{and} \;\;\; \left.\frac{\partial^{i} F}{\partial z_1^{i}}\right|_{z_1=0}\not=0 \;\; \text{for some} \;\; i \;\; \text{with} \;\; |i|=n+1.
\end{aligned}
\end{equation}
\end{definition}

\begin{remark}
The manifolds $\mathcal W_1$ and $\mathcal W_2$ coincide along a $c$-dimensional, smoothly embedded, compact disk if $F\equiv0$ in some neighborhood of $z_1=0$.
\end{remark}
It can be verified (see \cite{Domitrz}) that Definition \ref{coordinatedefinitionorder} is equivalent to the coordinate-independent definition provided in Section \ref{subsection1.1}. Consequently, the order of tangency is invariant under choices of coordinates that preserve $\mathcal W_1$ and $\mathcal W_2$ in the form \eqref{orderoftangencyeq1}.

\subsection{Adapted coordinates for a pair of homoclinic tangencies}\label{subsectionadaptedcoordinates}

\subsubsection{General case}\label{subsubsectionadaptedcoordinatesgeneralcase}

Let $f_0:\mathcal M^d\rightarrow \mathcal M^d$ be a $C^r$ diffeomorphism, where $d\ge 2$ and $r=1,.\ldots,\infty,\omega$. Let $f_0$ have a saddle periodic orbit $L$ whose stable and unstable manifolds $W^s(L)$ and $W^u(L)$ intersect over two orbits $\Gamma$, $\widetilde\Gamma$ of corank-$c$ homoclinic tangency. Let $L$ satisfy condition \textbf{\hyperlink{C1.}{C1}}, and $\Gamma,\widetilde\Gamma$ satisfy condition \textbf{\hyperlink{C2.}{C2}}. Partial-hyperbolicity condition \textbf{\hyperlink{C2.}{C2}} entails that $c$ satisfies $c\le\min\{m_s,m_u\}$, where $m_s$ and $m_u$ are the numbers of leading stable and leading unstable multipliers of $L$, respectively. Let $O$ be a point on the orbit $L$. We choose homoclinic points $M^{+},M^{-}\in\Gamma$ and $\widetilde M^{+},\widetilde M^{-}\in\widetilde\Gamma$ in such a way that $M^{+}, \widetilde M^{+}\in W^s_{loc}(O)$ and $M^{-}, \widetilde M^{-}\in W^u_{loc}(O)$. 

Assume $f_0$ is embedded into a finite-parameter $C^r$ family $(f_{\varepsilon})$. We fix a small neighborhood $U(O)$ of the point $O$, where the local map $T_{O}$ is defined for all small $\varepsilon$. We bring $T_O$ to the main normal form (see Section \ref{subsectionlocalmap}). We fix small neighborhoods $\Pi^{+},\Pi^{-}$ and $\widetilde \Pi^+,\widetilde\Pi^-$ of $M^{+},M^{-}$ and $\widetilde M^{+}, \widetilde M^{-}$, respectively, where the global maps $T_{\Gamma}$ and $T_{\widetilde\Gamma}$ are defined for all small $\varepsilon$. We denote
\begin{equation*}
\begin{aligned} 
&\mathcal W^s(O)=W^s_{loc}(O)\cap\Pi^+, \quad &\mathcal W^u(O)&=T_{\Gamma}(W^u_{loc}(O))\cap\Pi^+,\\
&\widetilde{\mathcal W}^s(O)=T^{-1}_{\widetilde\Gamma}(W^s_{loc}(O))\cap\tilde\Pi^-, \quad &\widetilde{\mathcal W}^u(O)&=W^u_{loc}(O)\cap\tilde\Pi^-.
\end{aligned}
\end{equation*}
The corresponding coordinates restricted to $\Pi^+$ and $\widetilde\Pi^-$ will be denoted by $(x,u,y,v)$ and $(\tilde x,\tilde u,\tilde y,\tilde v)$, respectively. We denote the tangent space coordinates at $M^+$ and $\widetilde M^-$ by $\left(\Delta x,\Delta u,\Delta y,\Delta v\right)$ and $(\Delta \tilde x,\Delta \tilde u,\Delta \tilde y,\Delta \tilde v)$. 

Our goal is to make linear transformations of the leading coordinates in the neighborhood $U(O)$ (note that they keep $T_O$ in the main normal form) such that equations for the manifolds $\mathcal W^u(O)$ and $\widetilde{\mathcal W}^s(O)$ take a particular form as summarised in Lemma \ref{lemmachoiceofcoordinates}. We do this in several steps. 

\textit{Step 1.} Since $T_O$ is in the main normal form, the manifolds $W^s_{loc}(O)$ and $W^u_{loc}(O)$ are straightened, i.e., $\mathcal W^s(O)$ and $\widetilde{\mathcal W}^u(O)$ are given by
\begin{equation}\label{lemmachoiceofcoordinateseq0}
\begin{aligned} 
&\mathcal W^s(O):\;(y,v)=(0,0), \quad\quad\quad \widetilde{\mathcal W}^u(O):\;(\tilde x,\tilde u)=(0,0).
\end{aligned}
\end{equation}
The points $M^+$ and $\widetilde M^-$ are given by $M^+=(x^+,u^+,0,0)$ and $\widetilde M^-=(0,0,\tilde y^-,\tilde v^-)$. The manifolds $\mathcal W^u(O)$ and $\widetilde{\mathcal W}^s(O)$ are smoothly embedded, so, for all small $\varepsilon$, they can be written as 
\begin{equation}\label{lemmachoiceofcoordinateseq3}
\mathcal W^u(O): \;\;
\begin{aligned} 
0&= P\left(x- x^+,u- u^+,y,v,\varepsilon\right),\\
\end{aligned}
\end{equation}
and
\begin{equation*}
\widetilde{\mathcal W}^s(O): \;\;
\begin{aligned} 
0&=\widetilde P\left(\tilde x,\tilde u,\tilde y-\tilde y^-,\tilde v - \tilde v^-,\varepsilon\right),
\end{aligned}
\end{equation*}
where $P:\mathbb R^{d+l}\to\mathbb R^{d_s}$ and $\widetilde P:\mathbb R^{d+l}\to\mathbb R^{d_u}$ are $C^r$-smooth functions such that
\begin{equation*}
\begin{aligned} 
\text{rank}\,\left.\frac{\partial P}{\partial (x,u,y,v)}\right|_{M^+}&=d_s &\text{and}& &\text{rank}\,\left.\frac{\partial \widetilde P}{\partial (\tilde x,\tilde u,\tilde y,\tilde v)}\right|_{\widetilde M^-}&=d_u.
\end{aligned}
\end{equation*}

\textit{Step 2.} According to condition \textbf{\hyperlink{C2.1.}{C2.1}}, the manifold $T_{\Gamma}(W^{ue}_{loc}(O))$ and the leaf $l^{ss}$ are transverse at the point $M^+$ when $\varepsilon=0$. Since $\mathcal W^u(O)\subset T_{\Gamma}(W^{ue}_{loc}(O))$, condition \textbf{\hyperlink{C2.1.}{C2.1}} entails
\begin{equation}\label{lemmachoiceofcoordinateseq5}
\begin{aligned} 
&\mathcal T_{M^+}(\mathcal W^u(O))\cap\mathcal T_{M^+}(l^{ss})=0.
\end{aligned}
\end{equation}
The equation for $\mathcal T_{M^+}(\mathcal W^u(O))$ is 
\begin{equation*}
\begin{aligned} 
&0=\left.\frac{\partial P}{\partial x}\right|_{M^+} \cdot \Delta x+\left.\frac{\partial P}{\partial u}\right|_{M^+} \cdot \Delta u+\left.\frac{\partial P}{\partial y}\right|_{M^+} \cdot \Delta y+\left.\frac{\partial P}{\partial v}\right|_{M^+} \cdot \Delta v,
\end{aligned}
\end{equation*}
while, by \eqref{foliationseq3}, the equation for $\mathcal T_{M^+}(l^{ss})$ is given by 
\begin{equation*}
\begin{aligned} 
&0=(\Delta x,\Delta y, \Delta v).
\end{aligned}
\end{equation*}
Consequently, condition \eqref{lemmachoiceofcoordinateseq5} implies that the equation $\left.\frac{\partial P}{\partial u}\right|_{M^+}\cdot \Delta u=0$ has only the zero solution. Therefore, $\text{rank} \left.\frac{\partial P}{\partial u}\right|_{M^+}=d_s$, and, applying the Implicit Function Theorem, we rewrite \eqref{lemmachoiceofcoordinateseq3} for all small $\varepsilon$ as
\begin{equation}\label{lemmachoiceofcoordinateseq8}
\mathcal W^u(O): \;\;
\begin{aligned} 
0&= Q\left(x- x^+, y, v,\varepsilon\right),\\
u-u^+&= Q^*\left(x- x^+, y, v,\varepsilon\right),\\
\end{aligned}
\end{equation}
where $Q,Q^*$ are $C^r$-smooth functions. Similarly, we rewrite $\widetilde{\mathcal W}^s(O)$ as
\begin{equation}\label{lemmachoiceofcoordinateseq9}
\widetilde{\mathcal W}^s(O): \;\;
\begin{aligned} 
0&= \widetilde Q\left(\tilde x,\tilde u,\tilde y-\tilde y^-,\varepsilon\right),\\
\tilde v-\tilde v^-&=\widetilde Q^*\left(\tilde x,\tilde u,\tilde y-\tilde y^-,\varepsilon\right),\\
\end{aligned}
\end{equation}
where $\widetilde Q,\widetilde Q^*$ are $C^r$-smooth functions.

\textit{Step 3.} Recall that $\Gamma$ and $\widetilde\Gamma$ are orbits of corank-$c$ tangency. By Definition \ref{coordinatedefinitioncorank}, considering that $\mathcal W^s(O),\widetilde{\mathcal W}^u(O)$ are given by \eqref{lemmachoiceofcoordinateseq0} and $\mathcal W^u(O), \widetilde{\mathcal W}^s\left(O\right)$ are given by \eqref{lemmachoiceofcoordinateseq8}, \eqref{lemmachoiceofcoordinateseq9}, we get
\begin{equation}\label{lemmachoiceofcoordinateseq10}
\begin{aligned} 
\text{corank}\,\left.\frac{\partial Q}{\partial x}\right|_{M^+}=c \quad \text{and} \quad \text{corank}\,\left.\frac{\partial \widetilde Q}{\partial \tilde y}\right|_{\widetilde M^-}=c \quad  \text{at}\quad \varepsilon=0.
\end{aligned}
\end{equation}
By the Singular Value Decomposition, for all small $\varepsilon$, there exist special orthogonal matrices $C_1$ and $C_2$ such that the matrix $C_1\cdot \left.\frac{\partial Q}{\partial x}\right|_{M^+} \cdot C_2$ is block-diagonal. By \eqref{lemmachoiceofcoordinateseq10}, this matrix consists of a $c\times c$ block (which vanishes at $\varepsilon=0$) and an $(m_s-c)\times(m_s-c)$ block with a non-zero determinant. We apply a coordinate transformation $x=C_2\cdot x^{new}$ and multiply the first equation in \eqref{lemmachoiceofcoordinateseq8} by $C_1$ from the left. As a result, the equations \eqref{lemmachoiceofcoordinateseq8}, for all small $\varepsilon$, become
\begin{equation}\label{lemmachoiceofcoordinateseq11}
\mathcal W^u(O): \;\;
\begin{aligned} 
0&= R_1\left(x- x^+, y, v,\varepsilon\right),\\
0&= R_2\left(x- x^+, y, v,\varepsilon\right),\\
u-u^+&= R_3\left(x- x^+, y, v,\varepsilon\right),\\
\end{aligned}
\end{equation}
where $R$ are $C^r$-smooth functions whose derivatives at $M^+$ satisfy
\begin{equation*}
\begin{aligned} 
\left.\frac{\partial R_1}{\partial x_1}\right|_{M^+}&=0 \quad \text{at}\quad \varepsilon=0,\\
\end{aligned}
\end{equation*}
\begin{equation}\label{lemmachoiceofcoordinateseq13}
\begin{aligned} 
\left.\frac{\partial R_1}{\partial x_2}\right|_{M^+}=0, \quad \left.\frac{\partial R_2}{\partial x_1}\right|_{M^+}=0, \quad \det\left.\frac{\partial R_2}{\partial x_2}\right|_{M^+}\not=0, \quad \text{for all small} \; \varepsilon.
\end{aligned}
\end{equation}
Using the third relation in \eqref{lemmachoiceofcoordinateseq13} and the Implicit Function Theorem, we recast \eqref{lemmachoiceofcoordinateseq11} as
\begin{equation}\label{lemmachoiceofcoordinateseq14}
\mathcal W^u(O): \;\;
\begin{aligned} 
0&= S_1\left(x_1- x_1^+, y, v,\varepsilon\right),\\
x_2- x_2^+&= S_2\left(x_1- x_1^+, y, v,\varepsilon\right),\\
u-u^+&= S_3\left(x_1- x_1^+, y, v,\varepsilon\right),\\
\end{aligned}
\end{equation}
where $S$ are $C^r$-smooth functions satisfying
\begin{equation}\label{lemmachoiceofcoordinateseq14.5}
\begin{aligned} 
&\left.\frac{\partial S_1}{\partial x_1}\right|_{M^+}=0 \quad  \text{at}\quad \varepsilon=0, \quad \text{and} \quad \left.\frac{\partial S_2}{\partial x_1}\right|_{M^+}=0 \quad \text{for all small $\varepsilon$}.\\
\end{aligned}
\end{equation}
Similarly, we perform a special orthogonal change of $\tilde y$-coordinates and multiply the first equation in \eqref{lemmachoiceofcoordinateseq9} by a special orthogonal matrix from the left. In this way, we rewrite \eqref{lemmachoiceofcoordinateseq9} as
\begin{equation}\label{lemmachoiceofcoordinateseq15}
\widetilde{\mathcal W}^s(O): \;\;
\begin{aligned} 
0&= \widetilde S_1\left(\tilde x,\tilde u,\tilde y_1-\tilde y_1^-,\varepsilon\right),\\
\tilde y_2-\tilde y_2^-&=\widetilde S_2\left(\tilde x,\tilde u,\tilde y_1-\tilde y_1^-,\varepsilon\right),\\
\tilde v-\tilde v^-&=\widetilde S_3\left(\tilde x,\tilde u,\tilde y_1-\tilde y_1^-,\varepsilon\right),\\
\end{aligned}
\end{equation}
where $\widetilde S$ are $C^r$-smooth functions satisfying
\begin{equation}\label{lemmachoiceofcoordinateseq15.5}
\begin{aligned} 
&\left.\frac{\partial \widetilde S_1}{\partial \tilde y_1}\right|_{\widetilde M^-}=0 \quad  \text{at}\quad \varepsilon=0, \quad \text{and} \quad \left.\frac{\partial \widetilde S_2}{\partial \tilde y_1}\right|_{\widetilde M^-}=0 \quad \text{for all small $\varepsilon$}.\\
\end{aligned}
\end{equation}

\textit{Step 4.} According to condition \textbf{\hyperlink{C2.2.}{C2.2}}, the manifold $W^{se}_{loc}(O)$ and the leaf $T_{\Gamma}(l^{uu})$ are transverse at the point $M^+$ at $\varepsilon=0$. Considering that $l^{uu}\subset W^{u}_{loc}(O)$, we can write this condition as the transversality of $\mathcal W^u(O)\cap W^{se}_{loc}(O)$ and $T_{\Gamma}(l^{uu})$ inside $\mathcal W^u(O)$. This implies that 
\begin{equation}\label{lemmachoiceofcoordinateseq16}
\begin{aligned} 
\dim\left(\mathcal T_{M^+}\left(\mathcal W^u(O)\right)\cap\mathcal T_{M^+}(W^{se}_{loc}(O))\right)\ge c.
\end{aligned}
\end{equation}
By \eqref{lemmachoiceofcoordinateseq14}, the equation for $\mathcal T_{M^+}(\mathcal W^u(O))$ has the form 
\begin{equation*}
\begin{aligned} 
&0=\left.\frac{\partial S_{1}}{\partial y}\right|_{M^+} \cdot \Delta y+\left.\frac{\partial S_{1}}{\partial v}\right|_{M^+} \cdot \Delta v, \quad (\Delta x_2,\Delta u)=O(\Delta x_1,\Delta y,\Delta v),
\end{aligned}
\end{equation*}
and, by \eqref{foliationseq4}, the equation for $\mathcal T_{M^+}(W^{se}_{loc}(O))$ is given by $\Delta v=0$. Consequently, condition \eqref{lemmachoiceofcoordinateseq16} implies 
\begin{equation}\label{lemmachoiceofcoordinateseq18}
\begin{aligned} 
&\text{rank}\, \left.\frac{\partial S_{1}}{\partial y}\right|_{M^+}=c.
\end{aligned}
\end{equation}
Applying similar reasoning to $\widetilde{\mathcal W}^s(O)$, we obtain 
\begin{equation}\label{lemmachoiceofcoordinateseq19}
\begin{aligned} 
&\text{rank}\, \left.\frac{\partial \widetilde S_{1}}{\partial \tilde x}\right|_{\widetilde M^-}=c.
\end{aligned}
\end{equation}
Note that the relations \eqref{lemmachoiceofcoordinateseq18} and \eqref{lemmachoiceofcoordinateseq19} mean that the matrices $\left.\frac{\partial S_{1}}{\partial y}\right|_{M^+}$ and $\left.\frac{\partial \widetilde S_{1}}{\partial \tilde x}\right|_{\widetilde M^-}$ have full rank. Therefore, these relations hold not only at $\varepsilon=0$, but for all small $\varepsilon$.

\vspace{0.5cm}

To etablish the lemma below, we impose the following genericity conditions: \textbf{\hypertarget{C3.}}

\begin{doublespace}
\textbf{\hypertarget{C3.1.}{C3.1.}} $\det\left.\frac{\partial S_{1}}{\partial y_1}\right|_{M^+}\not=0$ at $\varepsilon=0$;

\textbf{\hypertarget{C3.2.}{C3.2.}} $\det\left.\frac{\partial \widetilde S_{1}}{\partial \tilde x_1}\right|_{\widetilde M^-}\not=0$ at $\varepsilon=0$.
\end{doublespace}

Note that condition \textbf{\hyperlink{C3.}{C3}} holds for an open and dense subset of $C^r$ diffeomorphisms for which the periodic orbit $L$ satisfies \textbf{\hyperlink{C1.}{C1}} and the orbits $\Gamma,\widetilde\Gamma$ of a homoclinic tangency satisfy \textbf{\hyperlink{C2.}{C2}}. Openness is obvious and density is given by Remark \ref{remarkperturbationofcoefficients}. Also note that condition \textbf{\hyperlink{C3.1.}{C3.1}} always holds if $m_u=c$, and condition \textbf{\hyperlink{C3.2.}{C3.2}} always holds if $m_s=c$.

\begin{lemma}\label{lemmachoiceofcoordinates}
Under conditions \textbf{\hyperlink{C1.}{C1}}--\textbf{\hyperlink{C3.}{C3}}, the system of coordinates in $U(O)$ can be chosen in such a way that, for all small $\varepsilon$, the local map $T_{O}$ is in the main normal form, the equations of the manifolds $\mathcal W^s(O),\widetilde{\mathcal W}^u(O)$ are given by \eqref{lemmachoiceofcoordinateseq0}, and the equations of the manifolds $\mathcal W^u(O), \widetilde{\mathcal W}^s(O)$ are given by
\begin{equation}\label{lemmachoiceofcoordinateseq1}
\mathcal W^u(O): \;\;
\begin{aligned} 
y_1&=F(x_1-x_1^+,\varepsilon)+b\cdot y_2+\theta_1(x_1-x_1^+,y_2,v,\varepsilon),\\
x_2- x_2^+&=\theta_2(x_1-x_1^+,y_2,v,\varepsilon),\\
u-u^+&=\theta_3(x_1-x_1^+,y_2,v,\varepsilon),
\end{aligned}
\end{equation}
and
\begin{equation}\label{lemmachoiceofcoordinateseq2}
\widetilde{\mathcal W}^s(O): \;\;
\begin{aligned} 
\tilde x_1&=\widetilde F(\tilde y_1-\tilde y_1^-,\varepsilon)+\tilde b\cdot \tilde x_2+\tilde\theta_1(\tilde x_2,\tilde u,\tilde y_1-\tilde y_1^-,\varepsilon),\\
\tilde y_2-\tilde y_2^-&=\tilde\theta_2(\tilde x_2,\tilde u,\tilde y_1-\tilde y_1^-,\varepsilon),\\
\tilde v-\tilde v^-&=\tilde\theta_3(\tilde x_2,\tilde u,\tilde y_1-\tilde y_1^-,\varepsilon),
\end{aligned}
\end{equation}
where $\left(x_1,x_2,u,y_1,y_2,v\right),\left(\tilde x_1,\tilde x_2,\tilde u,\tilde y_1,\tilde y_2,\tilde v\right)\in\mathbb R^{c}\times\mathbb R^{m_s-c}\times\mathbb R^{d_s-m_s}\times\mathbb R^{c}\times\mathbb R^{m_u-c}\times\mathbb R^{d_u-m_u}$. At $\varepsilon=0$, one has
\begin{equation}\label{lemmachoiceofcoordinatesestimateseq1}
\begin{aligned} 
F(x_1=x_1^+,0)&=0, &\left.\frac{\partial F}{\partial x_1}\right|_{(x_1=x_1^+,0)}&=0, &\widetilde F(\tilde y_1=\tilde y_1^-,0)&=0, &\left.\frac{\partial \widetilde F}{\partial\tilde y_1}\right|_{(\tilde y_1=\tilde y_1^-,0)}&=0.
\end{aligned}
\end{equation}

For all small $\varepsilon$, one has
\begin{subequations}\label{lemmachoiceofcoordinatesestimateseq2}
\begin{align} 
\theta(x_1=x_1^+,0,0,\varepsilon)&=0, &\theta_1(x_1-x_1^+,0,0,\varepsilon)&=0, &\left.\frac{\partial\theta_1}{\partial y_2}\right|_{M^+}&=0, &\left.\frac{\partial\theta_2}{\partial x_1}\right|_{M^+}&=0;\\
\tilde\theta(0,0,\tilde y_1=\tilde y_1^-,\varepsilon)&=0, &\tilde\theta_1(0,0,\tilde y_1-\tilde y_1^-,\varepsilon)&=0, &\left.\frac{\partial\tilde\theta_1}{\partial \tilde x_2}\right|_{\widetilde M^-}&=0, &\left.\frac{\partial\tilde\theta_2}{\partial \tilde y_1}\right|_{\widetilde M^-}&=0.
\end{align}
\end{subequations}
\end{lemma}

\begin{demo}
By condition \textbf{\hyperlink{C3.1.}{C3.1}}, we express $y_1$ in the first equation of \eqref{lemmachoiceofcoordinateseq14} and plug it into other equations in \eqref{lemmachoiceofcoordinateseq14}. In the same way, by condition \textbf{\hyperlink{C3.2.}{C3.2}}, we express $\tilde x_1$ in the first equation of \eqref{lemmachoiceofcoordinateseq15} and plug it into other equations in \eqref{lemmachoiceofcoordinateseq15}. Thus, considering \eqref{lemmachoiceofcoordinateseq14.5},\eqref{lemmachoiceofcoordinateseq15.5}, the equations of the manifolds $\mathcal W^u(O), \widetilde{\mathcal W}^s(O)$ take the forms \eqref{lemmachoiceofcoordinateseq1},\eqref{lemmachoiceofcoordinateseq2} with the relations \eqref{lemmachoiceofcoordinatesestimateseq1},\eqref{lemmachoiceofcoordinatesestimateseq2}.
\end{demo}

In our subsequent proofs, we will use the coordinates constructed in the lemma above. Accordingly, we introduce the following
\begin{definition}\label{adaptedcoordinates}
The coordinates constructed in Lemma \ref{lemmachoiceofcoordinates} are called $(T_{O},\Gamma,\widetilde\Gamma)$-adapted.
\end{definition}

\begin{remark}
In $(T_{O},\Gamma,\widetilde\Gamma)$-adapted coordinates, the tangencies $\Gamma$ and $\widetilde\Gamma$ are of orders $n$ and $\tilde n$, respectively, if and only if the functions $F$ and $\widetilde F$ of \eqref{lemmachoiceofcoordinateseq1} and \eqref{lemmachoiceofcoordinateseq2} satisfy condition \eqref{coordinatedefinitionordereq}.
\end{remark}

\subsubsection{Pair of corank-1 homoclinic tangencies to bi-focus}

In this section, we consider a special case when the periodic orbit $L$ is a bi-focus and $c=1$. Let $\varphi_{\varepsilon}$ and $\psi_{\varepsilon}$ be the arguments of its leading stable and leading unstable multipliers, respectively. We assume the following genericity condition:

\textbf{\hypertarget{C3*.}{C3*.}} $(\varphi_0;\psi_0)\not=\left(\frac{\pi}{2};\frac{\pi}{2}\right)$. 

\begin{lemma}\label{lemmachoiceofcoordinatesandhomoclinicpoints}
Under conditions \textbf{\hyperlink{C2.}{C2}},\textbf{\hyperlink{C3*.}{C3*}}, the homoclinic point $M^+$ and the system of coordinates in $U(O)$ can be chosen in such a way that, for all small $\varepsilon$, the local map $T_{O}$ is in the special main normal form, the equations of the manifolds $\mathcal W^s(O), \widetilde{\mathcal W}^u(O)$ are given by \eqref{lemmachoiceofcoordinateseq0}, and the equations of the manifolds $\mathcal W^u(O), \widetilde{\mathcal W}^s(O)$ are given by \eqref{lemmachoiceofcoordinateseq1},\eqref{lemmachoiceofcoordinateseq2} with the relations \eqref{lemmachoiceofcoordinatesestimateseq1},\eqref{lemmachoiceofcoordinatesestimateseq2}.
\end{lemma}

\begin{demo}
We perform steps 1--4 from Section \ref{subsubsectionadaptedcoordinatesgeneralcase}, making the following modifications. In Step 1, we fix the coordinates in $U(O)$ such that the local map $T_{O}$ is not just in the main normal form, but in the special main normal form. Steps 2--4 remain unchanged (note that, in step 3, the change of coordinates is a rotation, so the local map $T_O$ remains in the special main normal form). In order to fulfil conditions \textbf{\hyperlink{C3.1.}{C3.1}} and \textbf{\hyperlink{C3.2.}{C3.2}}, we will iterate forward the manifold $\mathcal W^u(O)$ (this means that we choose a different homoclinic point $M^+$) and simultaneously rotate the $x$-coordinates as described below. Replicating proof of Lemma \ref{lemmachoiceofcoordinates} then will yield the result. 

Iterating  \eqref{eqlocal1}, we obtain the following formula for $T_O^k$:
\begin{equation}\label{lemmachoiceofcoordinatesandhomoclinicpointseq1}
\begin{aligned} 
x_k&= \mathcal A^k_1(\varepsilon) \cdot x+\vartheta_{1,k}(x,u,y,v,\varepsilon), \quad &y_k&=\mathcal B^k_1(\varepsilon) \cdot y+\vartheta_{3,k}(x,u,y,v,\varepsilon),\\ 
u_k&= \mathcal A^k_2(\varepsilon) \cdot u+\vartheta_{2,k}(x,u,y,v,\varepsilon), \quad &v_k&=\mathcal B^k_2(\varepsilon) \cdot v+\vartheta_{4,k}(x,u,y,v,\varepsilon),\\
\end{aligned}
\end{equation}
where
\begin{equation*}
\begin{aligned} 
&\mathcal A^k_1(\varepsilon)=\lambda_{\varepsilon}\cdot\begin{pmatrix}\cos(k\varphi_{\varepsilon})& -\sin (k\varphi_{\varepsilon})\\ \sin(k\varphi_{\varepsilon})& \cos(k\varphi_{\varepsilon}) \end{pmatrix}, \quad &\mathcal B^k_1(\varepsilon)&=\gamma_{\varepsilon}\cdot\begin{pmatrix}\cos(k\psi_{\varepsilon})& -\sin(k\psi_{\varepsilon})\\ \sin(k\psi_{\varepsilon})& \cos(k\psi_{\varepsilon})\end{pmatrix},
\end{aligned}
\end{equation*}
and the functions $\vartheta_{k}$ satisfy \eqref{eqlocal2.0}--\eqref{eqlocal2.2}. We iterate by $T_O^k$ the manifold $\mathcal W^u(O)$ given by \eqref{lemmachoiceofcoordinateseq14} and perform the rotation:
\begin{equation}\label{lemmachoiceofcoordinatesandhomoclinicpointseq3}
\begin{aligned} 
x_k&= \mathcal A^k_1(\varepsilon) \cdot x_k^{new}.
\end{aligned}
\end{equation}
Applying \eqref{lemmachoiceofcoordinatesandhomoclinicpointseq3} and viewing $T_O^k$ as the coordinate transformation \eqref{lemmachoiceofcoordinatesandhomoclinicpointseq1} applied to the equation \eqref{lemmachoiceofcoordinateseq14}, we write the equation of the manifold $T^k_O(\mathcal W^u(O))$ near $T^k_O(M^+)$ as
\begin{equation*}
T^k_O(\mathcal W^u(O)): \;\;
\begin{aligned} 
0&= S^*_{1}\left(x_{1,k}- x_1^+, y_k, v_k,\varepsilon\right),\\
x_{2,k}- x_2^+&= S^*_{2}\left(x_{1,k}- x_1^+, y_k, v_k,\varepsilon\right),\\
u_k-u^+&= S^*_{3}\left(x_{1,k}- x_1^+, y_k, v_k,\varepsilon\right),\\
\end{aligned}
\end{equation*}
where $S^*$ are $C^r$-smooth functions satisfying
\begin{equation*}
\begin{aligned} 
&\left.\frac{\partial S^*_{1}}{\partial x_{1,k}}\right|_{T^k_O(M^+)}=0 \quad  \text{at}\quad \varepsilon=0, \quad \text{and} \quad \left.\frac{\partial S^*_{2}}{\partial x_{1,k}}\right|_{T^k_O(M^+)}=0 \quad \text{for all small $\varepsilon$},
\end{aligned}
\end{equation*}
and
\begin{equation}\label{lemmachoiceofcoordinatesandhomoclinicpointseq6}
\begin{aligned} 
&\left.\frac{\partial S^*_{1}}{\partial y_{1,k}}\right|_{T^k_O(M^+)}=\gamma_{\varepsilon}^{-k}\cdot\left(\left.\frac{\partial S_1}{\partial y_1}\right|_{M^+}\cdot \cos(k\psi_{\varepsilon})-\left.\frac{\partial S_1}{\partial y_2}\right|_{M^+}\cdot \sin(k\psi_{\varepsilon})\right).
\end{aligned}
\end{equation}
We make the coordinate transformation \eqref{lemmachoiceofcoordinatesandhomoclinicpointseq3} in the whole neighborhood $U(O)$, so the equation of $\widetilde{\mathcal W}^s(O)$ near $\widetilde M^-$ changes. It takes the form 
\begin{equation*}
\widetilde{\mathcal W}^s(O): \;\;
\begin{aligned} 
0&= \widetilde S^*_1\left(\tilde x,\tilde u,\tilde y_1-\tilde y_1^-,\varepsilon\right),\\
\tilde y_2-\tilde y_2^-&=\widetilde S^*_2\left(\tilde x,\tilde u,\tilde y_1-\tilde y_1^-,\varepsilon\right),\\
\tilde v-\tilde v^-&=\widetilde S^*_3\left(\tilde x,\tilde u,\tilde y_1-\tilde y_1^-,\varepsilon\right),\\
\end{aligned}
\end{equation*}
where $\widetilde S^*$ are $C^r$-smooth functions satisfying
\begin{equation*}
\begin{aligned} 
&\left.\frac{\partial \widetilde S^*_1}{\partial \tilde y_1}\right|_{\widetilde M^-}=0 \quad  \text{at}\quad \varepsilon=0, \quad \text{and} \quad \left.\frac{\partial \widetilde S^*_2}{\partial \tilde y_1}\right|_{\widetilde M^-}=0 \quad \text{for all small $\varepsilon$},
\end{aligned}
\end{equation*}
and
\begin{equation}\label{lemmachoiceofcoordinatesandhomoclinicpointseq9}
\begin{aligned} 
&\left.\frac{\partial \widetilde S^*_{1}}{\partial \tilde x_1}\right|_{\widetilde M^-}=\lambda_{\varepsilon}^{k}\cdot\left(\left.\frac{\partial \widetilde S_1}{\partial \tilde x_1}\right|_{\widetilde M^-}\cdot \cos(k\varphi_{\varepsilon})+\left.\frac{\partial \widetilde S_1}{\partial \tilde x_2}\right|_{\widetilde M^-}\cdot \sin(k\varphi_{\varepsilon})\right).
\end{aligned}
\end{equation}
Thus, to complete the proof of the lemma, it remains to determine the number $k\in\mathbb N$ such that
\begin{equation*}
\begin{aligned} 
&\left.\frac{\partial S^*_{1}}{\partial y_{1,k}}\right|_{T^k_O(M^+)}\not=0 \quad \text{and} \quad \left.\frac{\partial \widetilde S^*_{1}}{\partial \tilde x_1}\right|_{\widetilde M^-}\not=0 \quad  \text{at} \;\; \varepsilon=0.
\end{aligned}
\end{equation*}
For the rest of the proof, we fix $\varepsilon=0$.

By condition \textbf{\hyperlink{C3*.}{C3*}}, the case $(\varphi_0;\psi_0)=\left(\frac{\pi}{2};\frac{\pi}{2}\right)$ is excluded. Without loss of generality, assume $\psi_0\not=\frac{\pi}{2}$. If $\left.\frac{\partial S_{1}}{\partial y_{1}}\right|_{M^+}\not=0$ and $\left.\frac{\partial \widetilde S_{1}}{\partial \tilde x_1}\right|_{\widetilde M^-}\not=0$, we need to do nothing and the lemma is proved. In view of \eqref{lemmachoiceofcoordinateseq18} and \eqref{lemmachoiceofcoordinateseq19}, the following conditions hold:
\begin{itemize}
\item If $\left.\frac{\partial S_1}{\partial y_1}\right|_{M^+}=0$, then $\left.\frac{\partial S_1}{\partial y_2}\right|_{M^+}\not=0$;

\item If $\left.\frac{\partial \widetilde S_1}{\partial \tilde x_1}\right|_{\widetilde M^-}=0$, then $\left.\frac{\partial \widetilde S_1}{\partial \tilde x_2}\right|_{\widetilde M^-}\not=0$.
\end{itemize}
The proof reduces to the consideration of the following three cases. 

\vspace{0.2cm}

\textbf{Case 1:} $\left.\frac{\partial S_{1}}{\partial y_{1}}\right|_{M^+}=0 \quad \text{and} \quad \left.\frac{\partial \widetilde S_{1}}{\partial \tilde x_1}\right|_{\widetilde M^-}=0$. Choose $k=1$. By \eqref{lemmachoiceofcoordinatesandhomoclinicpointseq6},
\begin{equation}\label{lemmachoiceofcoordinatesandhomoclinicpointseq10}
\begin{aligned} 
&\left.\frac{\partial S^*_{1}}{\partial y_{1,1}}\right|_{T_O(M^+)}=(-\gamma_{0}^{-1})\cdot\left.\frac{\partial S_1}{\partial y_2}\right|_{M^+}\cdot \sin \psi_{0}\not=0,\\
\end{aligned}
\end{equation}
and, by \eqref{lemmachoiceofcoordinatesandhomoclinicpointseq9},
\begin{equation}\label{lemmachoiceofcoordinatesandhomoclinicpointseq11}
\begin{aligned} 
&\left.\frac{\partial \widetilde S^*_{1}}{\partial \tilde x_1}\right|_{\widetilde M^-}=\lambda_{0}\cdot\left.\frac{\partial \widetilde S_1}{\partial \tilde x_2}\right|_{\widetilde M^-}\cdot \sin \varphi_{0}\not=0,
\end{aligned}
\end{equation}
which proves the lemma. 

\textbf{Case 2:} $\left.\frac{\partial S_{1}}{\partial y_{1}}\right|_{M^+}=0 \quad \text{and} \quad \left.\frac{\partial \widetilde S_{1}}{\partial \tilde x_1}\right|_{\widetilde M^-}\not=0$. Choose $k=1$. We obtain that \eqref{lemmachoiceofcoordinatesandhomoclinicpointseq10} holds. If $\left.\frac{\partial \widetilde S^*_{1}}{\partial \tilde x_1}\right|_{\widetilde M^-}\not=0$, the lemma is proved. Otherwise, suppose $\left.\frac{\partial \widetilde S^*_{1}}{\partial \tilde x_1}\right|_{\widetilde M^-}=0$, which implies 
\begin{equation}\label{lemmachoiceofcoordinatesandhomoclinicpointseq12}
\begin{aligned} 
&\left.\frac{\partial \widetilde S_1}{\partial \tilde x_1}\right|_{\widetilde M^-}\cdot \cos\varphi_{0}+\left.\frac{\partial \widetilde S_1}{\partial \tilde x_2}\right|_{\widetilde M^-}\cdot \sin\varphi_{0}=0.
\end{aligned}
\end{equation}
Choose $k=2$. Then, according to \eqref{lemmachoiceofcoordinatesandhomoclinicpointseq6},
\begin{equation*}
\begin{aligned} 
&\left.\frac{\partial S^*_{1}}{\partial y_{1,2}}\right|_{T^2_O(M^+)}=(-\gamma_{0}^{-2})\cdot\left.\frac{\partial S_1}{\partial y_2}\right|_{M^+}\cdot \sin \left(2\psi_{0}\right)\not=0.
\end{aligned}
\end{equation*}
By \eqref{lemmachoiceofcoordinatesandhomoclinicpointseq9} and \eqref{lemmachoiceofcoordinatesandhomoclinicpointseq12}, we obtain
\begin{equation*}
\begin{aligned} 
&\left.\frac{\partial \widetilde S^*_{1}}{\partial \tilde x_1}\right|_{\widetilde M^-}=\lambda_{0}^{2}\cdot\left(\left.\frac{\partial \widetilde S_1}{\partial \tilde x_1}\right|_{\widetilde M^-}\cdot \cos(2\varphi_{0})+\left.\frac{\partial \widetilde S_1}{\partial \tilde x_2}\right|_{\widetilde M^-}\cdot \sin(2\varphi_{0})\right)=(-\lambda_{0}^{2})\cdot\left.\frac{\partial \widetilde S_1}{\partial \tilde x_1}\right|_{\widetilde M^-}\not=0,
\end{aligned}
\end{equation*}
which proves the lemma.

\textbf{Case 3:} $\left.\frac{\partial S_{1}}{\partial y_{1}}\right|_{M^+}\not=0 \quad \text{and} \quad \left.\frac{\partial \widetilde S_{1}}{\partial \tilde x_1}\right|_{\widetilde M^-}=0$. Choose $k=1$. We get that \eqref{lemmachoiceofcoordinatesandhomoclinicpointseq11} is fulfilled. If $\left.\frac{\partial S^*_{1}}{\partial y_{1,1}}\right|_{T_O(M^+)}\not=0$, the lemma is proved. Suppose $\left.\frac{\partial S^*_{1}}{\partial y_{1,1}}\right|_{T_O(M^+)}=0$, which entails
\begin{equation}\label{lemmachoiceofcoordinatesandhomoclinicpointseq15}
\begin{aligned} 
&\left.\frac{\partial S_1}{\partial y_1}\right|_{M^+}\cdot \cos\psi_{0}-\left.\frac{\partial S_1}{\partial y_2}\right|_{M^+}\cdot \sin\psi_{0}=0.
\end{aligned}
\end{equation}
Choose $k=2$. Then, it follows from \eqref{lemmachoiceofcoordinatesandhomoclinicpointseq6} and \eqref{lemmachoiceofcoordinatesandhomoclinicpointseq15} that 
\begin{equation*}
\begin{aligned} 
&\left.\frac{\partial S^*_{1}}{\partial y_{1,2}}\right|_{T^2_O(M^+)}=\gamma_{0}^{-2}\cdot\left(\left.\frac{\partial S_1}{\partial y_1}\right|_{M^+}\cdot \cos(2\psi_{0})-\left.\frac{\partial S_1}{\partial y_2}\right|_{M^+}\cdot \sin(2\psi_{0})\right)=(-\gamma_{0}^{-2})\cdot\left.\frac{\partial S_1}{\partial y_1}\right|_{M^+}\not=0.
\end{aligned}
\end{equation*}
According to \eqref{lemmachoiceofcoordinatesandhomoclinicpointseq9},
\begin{equation*}
\begin{aligned} 
&\left.\frac{\partial \widetilde S^*_{1}}{\partial \tilde x_1}\right|_{\widetilde M^-}=\lambda^2_{0}\cdot\left.\frac{\partial \widetilde S_1}{\partial \tilde x_2}\right|_{\widetilde M^-}\cdot \sin(2\varphi_{0}).
\end{aligned}
\end{equation*}
Thus, if $\varphi_0\not=\frac{\pi}{2}$, then $\left.\frac{\partial \widetilde S^*_{1}}{\partial \tilde x_1}\right|_{\widetilde M^-}\not=0$, and the lemma is proved. Suppose $\varphi_0=\frac{\pi}{2}$ and choose $k=3$. Then, formula \eqref{lemmachoiceofcoordinatesandhomoclinicpointseq9} implies
\begin{equation*}
\begin{aligned} 
&\left.\frac{\partial \widetilde S^*_{1}}{\partial \tilde x_1}\right|_{\widetilde M^-}=\lambda^3_{0}\cdot\left.\frac{\partial \widetilde S_1}{\partial \tilde x_2}\right|_{\widetilde M^-}\cdot \sin(3\varphi_{0})=(-\lambda^3_{0})\cdot\left.\frac{\partial \widetilde S_1}{\partial \tilde x_2}\right|_{\widetilde M^-}\not=0.
\end{aligned}
\end{equation*}
By \eqref{lemmachoiceofcoordinatesandhomoclinicpointseq6} and \eqref{lemmachoiceofcoordinatesandhomoclinicpointseq15}, we obtain
\begin{equation*}
\begin{aligned} 
&\left.\frac{\partial S^*_{1}}{\partial y_{1,3}}\right|_{T^3_O(M^+)}=\gamma_{0}^{-3}\cdot\left(\left.\frac{\partial S_1}{\partial y_1}\right|_{M^+}\cdot \cos(3\psi_{0})-\left.\frac{\partial S_1}{\partial y_2}\right|_{M^+}\cdot \sin(3\psi_{0})\right)=\\
&=(-\gamma_{0}^{-3})\cdot 2\cdot\left.\frac{\partial S_1}{\partial y_1}\right|_{M^+}\cdot \cos\psi_{0}\not=0,
\end{aligned}
\end{equation*}
which proves the lemma.
\end{demo}

\begin{remark}
By applying the same strategy, it can be shown that Lemma \ref{lemmachoiceofcoordinatesandhomoclinicpoints} holds for the main cases $(1,1), (2,1), (1,2)$ $($without requiring condition \textbf{\hyperlink{C3*.}{C3*}}$)$.
\end{remark}

\section{Perturbations of tangencies}\label{sectionperturbationsoftangencies}

\subsection{Free splitting}\label{subsectionfreesplitting}

Let $f_0:\mathcal M^d\rightarrow \mathcal M^d$ be a $C^r$ diffeomorphism, where $d\ge 2$ and $r=2,.\ldots,\infty,\omega$. Let $f_0$ possess a non-trivial basic set $\Lambda$ of type $(d_s,d_u)$. Let $p$ and $q$ be points (possibly coinciding) in $\Lambda$ such that the stable manifold $W^s(p)$ and the unstable manifold $W^u(q)$ have a corank-$c$ tangency at the points of an orbit $\Gamma$. Note that $c$ satisfies $c\le\min\{d_s,d_u\}\le[\frac{d}{2}]$. Take a point $M\in\Gamma$ and a small neighborhood $\Pi$ of $M$. Let $\mathcal W^{s}_{M}(p)$ and $\mathcal W^{u}_{M}(q)$ be path-connected components of $M$ in $W^{s}(p)\cap\Pi$ and $W^{u}(q)\cap\Pi$, respectively.

Assume $f_0$ is embedded into a finite-parameter $C^r$ family of diffeomorphisms $(f_{\varepsilon})$. In the neighborhood $\Pi$, we choose coordinates $(z_1,z_2,z_3,z_4)\in\mathbb R^{c}\times\mathbb R^{d_s-c}\times\mathbb R^{c}\times\mathbb R^{d_u-c}$ such that, for all small $\varepsilon$, the manifolds $\mathcal W^{s}_{M}(p)$ and $\mathcal W^{u}_{M}(q)$ are represented as
\begin{equation}\label{splittingoftangenciesbasicseteq1}
\begin{aligned} 
\mathcal W^{s}_{M}(p):\;\;&(z_3,z_4)=(0,0),\\
\mathcal W^{u}_{M}(q):\;\;&z_3=F(z_1,\varepsilon)+\theta_{1}(z_1,z_4,\varepsilon), \;\; z_2=\theta_{2}(z_1,z_4,\varepsilon),
\end{aligned}
\end{equation}
where the functions $F$ and $\theta$ are $C^r$-smooth with respect to the coordinates and parameters. They satisfy (cf. \eqref{orderoftangencyeq2})
\begin{equation}\label{splittingoftangenciesbasicseteq2}
\begin{aligned} 
F\left(0,0\right)&=0, &\left.\frac{\partial F}{\partial z_1}\right|_{(0,0)}&=0, &\theta_1(z_1,0,\varepsilon)&=0.
\end{aligned}
\end{equation}
We call $F$ the \textit{splitting function} for the orbit $\Gamma$. Fix an integer $n'$ such that $0\le n'\le r-1$. The function $F$ can be expanded as
\begin{equation*}
\begin{aligned} 
&F(z_1,\varepsilon)=F(z_1,0)+\sum\limits_{|i|=0}^{n'} \eta_{i}\cdot z_1^i+O(||z_1||^{n'+1}),
\end{aligned}
\end{equation*}
where $\eta_{i}$ are smooth functions of $\varepsilon$ vanishing at $\varepsilon=0$; here $i=(i_1,\ldots,i_c)$ is a multi-index, so $z_1^i=z_1^{i_1}\cdot \ldots\cdot z_1^{i_c}$. We call $\eta_i$ the \textit{splitting coefficients}. We denote by $\mathbf H_{n'}$ the vector of all such coefficients from order $0$ to order $n'$.
\begin{definition}\label{definitionfreeperturbationorbit}
The family $(f_{\varepsilon})$ perturbs the orbit $\Gamma$ freely up to order $n'$ $($or we say that the tangency splits freely up to order $n')$ if 
\begin{equation*} 
\left.\text{rank}\;\frac{\partial\mathbf H_{n'}}{\partial\varepsilon}\right |_{\varepsilon=0}=\dim\mathbf H_{n'}. 
\end{equation*}
\end{definition}
While there is no formal link between the order $n$ of tangency $\Gamma$ and the order $n'$ of the perturbation freedom, for the most of the paper we consider the cases $n'=n$ and $n'=n-1$. 

Let the tangency $\Gamma$ be of order $n$, where $n\le r-1$. It is directly verified that the following lemma holds.  
\begin{lemma}\label{lemmasplittingdoesnotdependoncoordinates}
For any two choices of $C^r$ coordinates keeping $\mathcal W^{s}_{M}(p)$ and $\mathcal W^{u}_{M}(q)$ in the form \eqref{splittingoftangenciesbasicseteq1}, the vectors $\mathbf H_{n-1}$ of all splitting coefficients from order $0$ to order $n-1$ are related by a diffeomorphism.
\end{lemma}

\begin{remark}
Thus, the notion of the free perturbation up to order $n-1$ does not depend on the choice of coordinates near $M$, and also on the choice of the point $M$ in the orbit $\Gamma$ $($since the change of the point is equivalent to the change of coordinates$)$. On the other hand, the free perturbation of a higher order is coordinate dependent and also depends on the choice of $M$. Consequently, whenever a free perturbation of order $n$ or higher is discussed, it is assumed that both a point of tangency and a coordinate system have been fixed.
\end{remark}

We can define free perturbation with a different choice of coordinates in $\Pi$, when the manifold $\mathcal W^{u}_{M}(q)$ is straightened instead of $\mathcal W^{s}_{M}(p)$. Specifically, we choose coordinates $(\tilde z_1,\tilde z_2,\tilde z_3,\tilde z_4)\in\mathbb R^{c}\times\mathbb R^{d_s-c}\times\mathbb R^{c}\times\mathbb R^{d_u-c}$ such that, for all small $\varepsilon$, the manifolds $\mathcal W^{s}_{M}(p)$ and $\mathcal W^{u}_{M}(q)$ are given by 
\begin{equation}\label{splittingoftangenciesinverseeq1}
\begin{aligned} 
\mathcal W^{s}_{M}(p):\;\;&\tilde z_1=\widetilde F(\tilde z_3,\varepsilon)+\tilde\theta_{1}(\tilde z_3,\tilde z_2,\varepsilon), \;\; \tilde z_4=\tilde\theta_{2}(\tilde z_3,\tilde z_2,\varepsilon),\\
\mathcal W^{u}_{M}(q):\;\;&(\tilde z_1,\tilde z_2)=(0,0),
\end{aligned}
\end{equation}
where the functions $\widetilde F$ and $\tilde\theta$ are $C^r$-smooth with respect to the coordinates and parameters. They satisfy
\begin{equation*}
\begin{aligned} 
\widetilde F(0,0)&=0, &\left.\frac{\partial \widetilde F}{\partial \tilde z_3}\right|_{(0,0)}&=0, &\tilde\theta_1(\tilde z_3,0,\varepsilon)&=0.
\end{aligned}
\end{equation*}
The function $\widetilde F$ can be expanded as
\begin{equation*}
\begin{aligned} 
&\widetilde F(\tilde z_3,\varepsilon)=\widetilde F(\tilde z_3,0)+\sum\limits_{|i|=0}^{n'} \tilde\eta_{i}\cdot \tilde z_3^i+O(||\tilde z_3||^{n'+1}),
\end{aligned}
\end{equation*}
where $\tilde\eta_{i}$ are smooth functions of $\varepsilon$ vanishing at $\varepsilon=0$. We also call $\tilde\eta_i$ the splitting coefficients and say that the family $(f_{\varepsilon})$ perturbs the orbit $\Gamma$ freely up to order $n'$ if 
\begin{equation*} 
\left.\text{rank}\;\frac{\partial\widetilde{\mathbf H}_{n'}}{\partial\varepsilon}\right |_{\varepsilon=0}=\dim\widetilde{\mathbf H}_{n'},
\end{equation*}
where $\widetilde{\mathbf H}_{n'}$ is a vector of all coefficients $\tilde\eta_i$ from order $0$ to order $n'$.

Note that free perturbation in this sense is equivalent to the one defined before (see Definition \ref{definitionfreeperturbationorbit}). Indeed, we make a change of variables:
\begin{equation*}
\begin{aligned} 
z_1&=\tilde z_3 , &z_2&=\tilde z_2, &z_3&=\tilde z_1-\widetilde F(\tilde z_3,\varepsilon)-\tilde\theta_{1}(\tilde z_3,\tilde z_2,\varepsilon), &z_4&=\tilde z_4-\tilde\theta_{2}(\tilde z_3,\tilde z_2,\varepsilon).
\end{aligned}
\end{equation*}
As a result, we rewrite \eqref{splittingoftangenciesinverseeq1} as
\begin{equation*}
\begin{aligned} 
\mathcal W^{s}_{M}(p):\;\;&(z_3,z_4)=(0,0),\\
\mathcal W^{u}_{M}(q):\;\;&z_3=-\widetilde F(z_1,\varepsilon)-\tilde\theta_{1}(z_1,0,\varepsilon), \;\; z_2=0,
\end{aligned}
\end{equation*}
i.e., $\mathcal W^{s}_{M}(p)$ and $\mathcal W^{u}_{M}(q)$ are given by the formula of type \eqref{splittingoftangencieseq1} with the splitting coefficients
\begin{equation*}
\begin{aligned} 
&\eta_i=-\tilde\eta_i \quad \text{for all} \;\;\; |i|\in\{0,\ldots,n'\}.
\end{aligned}
\end{equation*}

In the following lemma (the case $c=2$ was done in \cite{Mints}), one can think of free splitting in any of the two above meanings. For simplicity, we take the number of parameters $\varepsilon$ to be the same as the number of the splitting coefficients up to order $n$. 
\begin{lemma}\label{lemmasplittangencyinbasicset}
Let $f_0:\mathcal M^d\rightarrow \mathcal M^d$ be a $C^r$ diffeomorphism, where $d\ge 2$ and $r=2,.\ldots,\infty,\omega$. Let $f_0$ have a non-trivial basic set $\Lambda$, let $p,q\in\Lambda$ be points such that $W^s(p)$ and $W^u(q)$ contain an orbit $\Gamma$ of corank-$c$ tangency of order $n$ $(n\le r-1)$, and let $\mathcal O_1,\mathcal O_2\in\Lambda$ be two periodic orbits. Suppose $f_0$ is embedded into a finite-parameter $C^r$ family $(f_{\varepsilon})$ in which the tangency $\Gamma$ splits freely up to order $n$. Then there exists a sequence $\langle\varepsilon_\alpha\rangle_{\alpha\in\mathbb N}$ with $\varepsilon_{\alpha}\xrightarrow{}0$ such that $f_{\varepsilon_{\alpha}}$ has an orbit of corank-$c$ homoclinic tangency of order $n$ between $W^s(\mathcal O_1)$ and $W^u(\mathcal O_2)$.
\end{lemma}

\begin{demo}
For all small $\varepsilon$, the manifolds $\mathcal W^s_{M}(p)$ and $\mathcal W^u_{M}(q)$ are specified by
\begin{equation*}
\begin{aligned} 
\mathcal W^s_{M}(p):\;\;&(z_3,z_4)=(0,0),\\
\mathcal W^u_{M}(q):\;\;&z_3=F(z_1,\varepsilon)+\theta_1(z_1,z_4,\varepsilon), \;\; z_2=\theta_{2}(z_1,z_4,\varepsilon),
\end{aligned}
\end{equation*}
with the estimates \eqref{splittingoftangenciesbasicseteq2} and the splitting function $F$ given by 
\begin{equation*}
\begin{aligned} 
F(z_1,\varepsilon)&=\sum\limits_{|i|=0}^{n} \eta_{i}\cdot z_1^i+O(||z_1||^{n+1}),
\end{aligned}
\end{equation*}
where the splitting coefficients $\eta_i$ vanish at $\varepsilon=0$, and $\left.\frac{\partial^{j} F}{\partial z_1^{j}}\right|_{z_1=0}\not=0$ for some $|j|=n+1$. The tangency $\Gamma$ splits freely, so one has
\begin{equation}\label{lemmasplittangencyinbasicseteq2} 
\left.\text{rank}\;\frac{\partial\mathbf H_{n}}{\partial\varepsilon}\right |_{\varepsilon=0}=\dim\mathbf H_{n}. 
\end{equation}
Since the matrix $\left.\frac{\partial\mathbf H_{n}}{\partial\varepsilon}\right |_{\varepsilon=0}$ is quadratic, condition \eqref{lemmasplittangencyinbasicseteq2} means that 
\begin{equation}\label{lemmasplittangencyinbasicseteq3} 
\left.\det\frac{\partial\mathbf H_{n}}{\partial\varepsilon}\right |_{\varepsilon=0}\not=0. 
\end{equation}

The periodic orbits $\mathcal O_1$, $\mathcal O_2$ belong to the basic set $\Lambda$, therefore the manifolds $W^s(\mathcal O_1)$ and $W^u(\mathcal O_2)$ accumulate on $W^s(p)$ and $W^u(q)$, respectively. We fix path-connected components of the sets $W^s(\mathcal O_1)\cap\Pi$ and $W^u(\mathcal O_2)\cap\Pi$ that are $C^r$-close to $\mathcal W^{s}_{M}(p)$ and $\mathcal W^{u}_{M}(q)$, and denote them as $\mathcal W^s(\mathcal O_1)$ and $\mathcal W^u(\mathcal O_2)$, respectively. Their equations can be written as
\begin{equation*}
\begin{aligned} 
\mathcal W^s(\mathcal O_1):\;\;&z_3=\beta_3(z_1,z_2,\varepsilon) \;\; z_4=\beta_4(z_1,z_2,\varepsilon),\\
\mathcal W^u(\mathcal O_2):\;\;&z_3=G(z_1,\varepsilon)+\beta_1(z_1,z_4,\varepsilon), \;\; z_2=\beta_{2}(z_1,z_4,\varepsilon),
\end{aligned}
\end{equation*}
where $\beta_3, \beta_4$ are $C^r$-close to the zero function; $G, \beta_1, \beta_2$ are $C^r$-close to $F, \theta_1, \theta_2$, respectively. 

We now perform a coordinate transformation:
\begin{equation*}
\begin{aligned} 
&z_3^{new}=z_3-\beta_3(z_1,z_2,\varepsilon), \qquad z_4^{new}=z_4-\beta_4(z_1,z_2,\varepsilon).
\end{aligned}
\end{equation*}
Under these new variables, $\mathcal W^s(\mathcal O_1)$ becomes the coordinate subspace $(z_3,z_4)=(0,0)$. Applying the Implicit Function Theorem to the transformed equations of $\mathcal W^u(\mathcal O_2)$, we express it as  
\begin{equation*}
\begin{aligned} 
\mathcal W^u(\mathcal O_2):\;\;&z_3=\widehat G(z_1,\varepsilon)+\hat \beta_1(z_1,z_4,\varepsilon), \;\; z_2=\hat \beta_{2}(z_1,z_4,\varepsilon),
\end{aligned}
\end{equation*}
where the functions $\widehat G,\hat \beta_1,\hat \beta_2$ are $C^r$-close to $G,\beta_1,\beta_2$. We expand $\widehat G$ as 
\begin{equation}\label{lemmasplittangencyinbasicseteq5}
\begin{aligned} 
&\widehat G(z_1,\varepsilon)=\sum\limits_{|i|=0}^{n} \rho_{i}\cdot z_1^i+O(||z_1||^{n+1}),
\end{aligned}
\end{equation}
where $\rho_i$ are smooth functions of $\varepsilon$. Let $\mathbf P_n$ be the vector of all the coefficients $\rho_i$. The functions $\widehat G$ and $F$ are $C^r$-close, so the vectors $\mathbf P_n$ and $\mathbf H_n$ are close. Considering \eqref{lemmasplittangencyinbasicseteq3}, we conclude that 
\begin{equation}\label{lemmasplittangencyinbasicseteq6} 
\left.\det\frac{\partial\mathbf P_{n}}{\partial\varepsilon}\right |_{\varepsilon=0}\not=0. 
\end{equation}

By \eqref{lemmasplittangencyinbasicseteq5}, the diffeomorphism $f_{\varepsilon}$ has an orbit of corank-$c$ tangency of order not lower than $n$ between $\mathcal W^s(\mathcal O_1)$ and $\mathcal W^u(\mathcal O_2)$ if 
\begin{equation}\label{lemmasplittangencyinbasicseteq7} 
\mathbf P_n=0.
\end{equation}
Employing \eqref{lemmasplittangencyinbasicseteq6} and the Implicit Function Theorem, we obtain the solution to \eqref{lemmasplittangencyinbasicseteq7} that yields the desired tangency. Finally, since $\widehat G$ is $C^r$-close to $F$ and the $(n+1)$-th derivative of $F$ at $z_1=0$ is non-zero, the tangency is exactly of order $n$.
\end{demo}

\subsection{Free and independent splitting}\label{subsectionfreeandindependentsplitting}
Let $p_1,q_1,\ldots,p_h,q_h$ be points (possibly coinciding) in $\Lambda$ such that for each $j\in\{1,\ldots,h\}$ the stable manifold $W^s(p_j)$ and the unstable manifold $W^u(q_j)$ intersect over an orbit $\Gamma_j$ of corank-$c_j$ tangency. Take a point $M_j\in\Gamma_j$ and a small neighborhood $\Pi_j$ of $M_j$. We define $\mathcal W^{s}_{M_j}(p_j)$ and $\mathcal W^{u}_{M_j}(q_j)$ as the path-connected components of $M_j$ in $W^s(p_j)\cap\Pi_j$ and $W^u(q_j)\cap\Pi_j$, respectively.

Assume $f_0$ is embedded into a finite-parameter $C^r$ family of diffeomorphisms $(f_{\varepsilon})$. As in Section \ref{subsectionfreesplitting}, we choose coordinates $(z_1,z_2,z_3,z_4)\in\mathbb R^{c_j}\times\mathbb R^{d_s-c_j}\times\mathbb R^{c_j}\times\mathbb R^{d_u-c_j}$ in the neighborhood $\Pi_j$ such that, for all small $\varepsilon$, the manifolds $\mathcal W^{s}_{M_j}(p_j)$ and $\mathcal W^{u}_{M_j}(q_j)$ are given by 
\begin{equation}\label{splittingoftangencieseq1}
\begin{aligned} 
\mathcal W^{s}_{M_j}(p_j):\;\;&(z_3,z_4)=(0,0),\\
\mathcal W^{u}_{M_j}(q_j):\;\;&z_3=F^j(z_1,\varepsilon)+\theta^j_{1}(z_1,z_4,\varepsilon), \;\; z_2=\theta^j_{2}(z_1,z_4,\varepsilon),
\end{aligned}
\end{equation}
where the functions $F^j,\theta^j$ are $C^r$-smooth with respect to the coordinates and parameters, and 
\begin{equation}\label{splittingoftangencieseq2}
\begin{aligned} 
F^j(0,0)&=0, &\left.\frac{\partial F^j}{\partial z_1}\right|_{(0,0)}&=0, &\theta^j_1(z_1,0,\varepsilon)&=0.
\end{aligned}
\end{equation}
Fix an integer $n'_j$ such that $0\le n'_j\le r-1$. The splitting function $F^j$ can be expanded as
\begin{equation*}
\begin{aligned} 
&F^j(z_1,\varepsilon)=F^j(z_1,0)+\sum\limits_{|i|=0}^{n'_j} \eta^j_{i}\cdot z_1^i+O(||z_1||^{n'_j+1}),
\end{aligned}
\end{equation*}
where $\eta^j_{i}$ are splitting coefficients. We denote by $\mathbf H^j_{n'_j}$ the vector of all such coefficients from order $0$ to order $n'_j$.
\begin{definition}\label{definitionfreeindependentperturbationorbit}
The family $(f_{\varepsilon})$ perturbs the orbits $\Gamma_1,\ldots,\Gamma_h$ freely and independently up to orders $n'_1,\ldots, n'_h$ $($or we say that the tangencies split freely and independently up to orders $n'_1,\ldots, n'_h)$ if 
\begin{equation}\label{corank2splitindependently}
\left.\text{rank}\;\frac{\partial\Big(\mathbf H^{1}_{n'_1},\ldots,\mathbf H^{h}_{n'_h}\Big)}{\partial\varepsilon}\right|_{\varepsilon=0}=\sum\limits_{j=1}^h\dim\mathbf H^{j}_{n'_j}.
\end{equation} 
\end{definition}
\begin{remark}
For finite $r$, the above constructions can be extended to the case $n'_j=r$ $(j\in\{1,\ldots,h\})$ if the $C^r$ family $(f_{\varepsilon})$ is such that the splitting coefficients $\eta^j_{i}$ are $C^1$ functions of $\varepsilon$ for all $|i|=r$. 
\end{remark}

The existence of parametric families satisfying \eqref{corank2splitindependently} follows from the lemma below. A similar result was established for two-dimensional diffeomorphisms in \cite{arbitraryorderconservative}, and for the case $c=2$ in \cite{Mints}. 
\begin{lemma}\label{lemmasplittingoftangency}
Let $f_0:\mathcal M^d\rightarrow \mathcal M^d$ be a $C^r$ diffeomorphism, where $d\ge 2$ and $r=2,.\ldots,\infty,\omega$. Let $f_0$ have a basic set $\Lambda$, and let $p_1,q_1,\ldots,p_h,q_h$ be points in $\Lambda$ such that for each $j\in\{1,\ldots,h\}$, the manifolds $W^s(p_j)$ and $W^u(q_j)$ intersect over an orbit $\Gamma_j$ of corank-$c_j$ tangency. Then, for any integers $n'_1,\ldots,n'_h$ $(0\le n'_j\le r)$, there exists a finite-parameter $C^{\infty}$ family $(\mathcal F_{\varepsilon})$ such that $\mathcal F_{0}=id$ and the family $(\mathcal F_{\varepsilon}\circ f_0)$ perturbs the orbits $\Gamma_1,\ldots,\Gamma_h$ freely and independently up to orders $n'_1,\ldots,n'_h$. If all $n'_{j}\le r-1$, then the family $(\mathcal F_{\varepsilon})$ can be chosen analytic.
\end{lemma}

\begin{demo} \textit{Construction of the $C^{\infty}$ family $(\mathcal F_{\varepsilon})$.} Let $B^{\rho}_{M}$ denote a ball of radius $\rho$ centered at the point $M$. Choose a sufficiently small $\rho>0$ such that the $2\rho$-balls $B^{2\rho}_{M_1},\ldots,B^{2\rho}_{M_h}$ are mutually disjoint. For each $j\in\{1,\ldots h\}$, we write the manifolds $\mathcal W^{s}_{M_j}(p_j), \mathcal W^{u}_{M_j}(q_j)$ inside the $\rho$-ball $B^{\rho}_{M_j}$ in the form \eqref{splittingoftangencieseq1},\eqref{splittingoftangencieseq2}.

Let $\xi^{\rho}_j$ be a $C^{\infty}$-smooth cut-off function such that for all small $\varepsilon$ it is equal to zero outside $B^{2\rho}_{M_j}$ and is equal to one inside $B^{\rho}_{M_j}$. Let $(\mathcal F_{\varepsilon})$ be a finite-parameter $C^{\infty}$ family such that for all $\varepsilon$ the diffeomorphism $\mathcal F_{\varepsilon}$ is an identity map outside the $2\rho$-balls $B^{2\rho}_{M_1},\ldots,B^{2\rho}_{M_h}$, and inside each ball $B^{2\rho}_{M_j}$ it acts as
\begin{equation}\label{lemmasplittingoftangencyeq0}
\begin{aligned} 
&(z_1,z_2,z_3,z_4)\mapsto(z_1,z_2, z_3+\xi^j_{\rho}\cdot\sum\limits_{|i|=0}^{n'_t} \eta^j_{i}\cdot z_1^i,\,z_4),
\end{aligned} 
\end{equation} 
where $\varepsilon=\big(\mathbf H^1_{n'_1},\ldots,\mathbf H^h_{n'_h}\big)$. By construction, $\mathcal F_{\varepsilon}=id$ at $\varepsilon=0$.

By definition, $\mathcal W^{s}_{M_j}(p_j)$ is an iteration of a small piece of $W^{s}_{loc}(p_j)$ by $f_0^{-1}$, while $\mathcal W^{u}_{M_j}(q_j)$ is an iteration of a small piece of $W^{u}_{loc}(q_j)$ by $f_0$. Since the diffeomorphism $\mathcal F_{\varepsilon}$ is the identity map everywhere outside $B^{2\rho}_{M_j}$, the manifolds $W^{s}_{loc}(p_j)$ and $W^{u}_{loc}(q_j)$ do not move as $\varepsilon$ varies. The same holds true for $\mathcal W^{s}_{M_j}(p_j)$ and $f_0\big(\mathcal W^{u}_{f_0^{-1}(M_j)}(q_j)\big)$. Therefore, inside the ball $B^{\rho}_{M_j}$, the equation for the manifold $\mathcal W^{s}_{M_j}(p_j)$ remains unchanged, whereas the equation for the manifold $\mathcal W^{u}_{M_j}(q_j)=\mathcal F_{\varepsilon}\circ f_0\big(\mathcal W^{u}_{f_0^{-1}(M_j)}(q_j)\big)$ becomes
\begin{equation}\label{lemmasplittingoftangencyeq1}
\begin{aligned} 
\mathcal W^{u}_{M_j}(q_j):\;\;& z_3=\sum\limits_{|i|=0}^{n'_j} \eta^j_{i}\cdot z_1^i+F^j(z_1,0)+\theta^j_{1}(z_1,z_4,0), \;\; z_2=\theta^j_{2}(z_1,z_4,0).
\end{aligned} 
\end{equation} 
By assumption, the splitting coefficients $\eta^j_i$ (where $0\le |i|\le n'_j\le r$) depend smoothly on the parameters $\varepsilon$, so the derivative of each coefficient with respect to $\varepsilon$ is well-defined. Considering that $\varepsilon=\big(\mathbf H^1_{n'_1},\ldots,\mathbf H^h_{n'_h}\big)$, the matrix $\left.\frac{\partial\big(\mathbf H^{1}_{n'_1},\ldots,\mathbf H^{h}_{n'_h}\big)}{\partial\varepsilon}\right|_{\varepsilon=0}$ is the identity matrix. Consequently, condition \eqref{corank2splitindependently} is satisfied, and the parametric family $(\mathcal F_{\varepsilon}\circ f_0)$ perturbs the orbits $\Gamma_1,\ldots,\Gamma_h$ freely and independently up to orders $n'_1,\ldots,n'_h$, as required.

\textit{Construction of the $C^{\omega}$ family $(\mathcal F_{\varepsilon})$.} Let $r'=\max\limits_{j\in\{1,\ldots,h\}} n'_j\le r-1$. Let $(\mathcal F_{\varepsilon})^{new}$ be a finite-parameter analytic family that is $C^{\infty}$-close to the family $(\mathcal F_{\varepsilon})$. By condition \eqref{corank2splitindependently}, the square matrix $\left.\frac{\partial\big(\mathbf H^{1}_{n'_1},\ldots,\mathbf H^{h}_{n'_h}\big)}{\partial\varepsilon}\right|_{\varepsilon=0}$ has full rank. Therefore, condition \eqref{corank2splitindependently} can be recast as 
\begin{equation*}
\begin{aligned} 
&\det\left.\frac{\partial\Big(\mathbf H^{1}_{n'_1},\ldots,\mathbf H^{h}_{n'_h}\Big)}{\partial\varepsilon}\right|_{\varepsilon=0}\not=0.
\end{aligned}
\end{equation*}
Since this determinant is non-zero and the parametric family $(\mathcal F^{new}_{\varepsilon}\circ f_0)$ is $C^{r'+1}$-close to $(\mathcal F_{\varepsilon}\circ f_0)$, the property is preserved for the family $(\mathcal F^{new}_{\varepsilon}\circ f_0)$. Thus, this family perturbs the orbits $\Gamma_1,\ldots,\Gamma_h$ freely and independently up to orders $n'_1,\ldots, n'_h$.
\end{demo}
 
 \subsection{Bi-focus case}
In Section \ref{subsectioncorank2tangencies}, we consider the case where the orbit $L$ is a bi-focus. We examine parametric families that not only perturb freely and independently an arbitrary (finite) number of orbits of tangencies but also allow for the independent variation of the arguments $\varphi_{\varepsilon},\psi_{\varepsilon}$ of the leading multipliers of $L$. That is:
\begin{equation}\label{corank2splitindependentlyandbifocus}
\left.\text{rank}\;\frac{\partial\Big(\mathbf H^{1}_{n'_1},\ldots,\mathbf H^{h}_{n'_h},\varphi_{\varepsilon},\psi_{\varepsilon}\Big)}{\partial\varepsilon}\right|_{\varepsilon=0}=\sum\limits_{j=1}^h\dim\Big(\mathbf H^{j}_{n'_j}\Big)+2.
\end{equation} 

\begin{lemma}\label{lemmasplittingoftangencybifocus}
If the orbit $L$ is bi-focus, then the family $(\mathcal F_{\varepsilon}\circ f_0)$ of Lemma \ref{lemmasplittingoftangency} can be chosen to satisfy condition \eqref{corank2splitindependentlyandbifocus}.
\end{lemma}

\begin{demo}
Throughout the proof, we maintain the same notations as in Lemma \ref{lemmasplittingoftangency}. We construct the $C^{\infty}$ family $(\mathcal F_{\varepsilon})$ essentially in the same way as in Lemma \ref{lemmasplittingoftangency}, with a modification within a small $2\rho$-ball $B^{2\rho}_{O}$ centered at the point $O\in L$. To obtain the $C^{\omega}$ family $(\mathcal F_{\varepsilon})$, we apply the same approximation argument used previously; we omit the proof here for brevity.

Choose a sufficiently small $\rho>0$ such that the $2\rho$-balls $B^{2\rho}_{M^1},\ldots,B^{2\rho}_{M^h},B^{2\rho}_{O}$ are mutually disjoint. In the ball $B^{2\rho}_{O}$, we introduce coordinates $(x,u,y,v)$, where $(x,u,y,v)\in\mathbb R^2\times\mathbb R^{d_s-2}\times\mathbb R^2\times\mathbb R^{d_u-2}$, such that the local map $T_O$ is in the special main normal form. Let $\xi^{\rho}$ be a $C^{\infty}$-smooth cut-off function such that, for all small $\varepsilon$, it is equal to zero outside $B^{2\rho}_{O}$ and is equal to one inside $B^{\rho}_{O}$. 

Let $(\mathcal F_{\varepsilon})$ be a finite-parameter $C^{\infty}$ family such that for all $\varepsilon$ the diffeomorphism $\mathcal F_{\varepsilon}$ is an identity map outside the $2\rho$-balls $B^{2\rho}_{M_1},\ldots,B^{2\rho}_{M_h}, B^{2\rho}_{O}$, inside each ball $B^{2\rho}_{M_j}$ acts by \eqref{lemmasplittingoftangencyeq0}, and inside the ball $B^{2\rho}_{O}$ acts as
\begin{equation*}
\begin{aligned} 
&(x,u,y,v)\mapsto(\mathcal R_1\cdot x,u,\mathcal R_2\cdot y,v),
\end{aligned} 
\end{equation*} 
where
\begin{equation*}
\begin{aligned} 
&\mathcal R_1=\begin{pmatrix}\cos(\xi^{\rho}\cdot\varepsilon_1)& -\sin(\xi^{\rho}\cdot\varepsilon_1)\\ \sin(\xi^{\rho}\cdot\varepsilon_1)& \cos(\xi^{\rho}\cdot\varepsilon_1) \end{pmatrix} \quad \text{and} \quad \mathcal R_2=\begin{pmatrix}\cos(\xi^{\rho}\cdot\varepsilon_2)& -\sin(\xi^{\rho}\cdot\varepsilon_2)\\ \sin(\xi^{\rho}\cdot\varepsilon_2)& \cos(\xi^{\rho}\cdot\varepsilon_2) \end{pmatrix}.
\end{aligned} 
\end{equation*} 
Here, $\varepsilon=\big(\mathbf H^1_{n'_1},\ldots,\mathbf H^h_{n'_h},\varepsilon_1,\varepsilon_2\big)$. Inside the ball $B^{\rho}_{O}$, the new local map $\mathcal F_{\varepsilon}\circ T_{O}$ remains in the special main normal form. Direct computation immediately gives 
\begin{equation*}
\begin{aligned} 
&\varphi_{\varepsilon}=\varphi_0+\varepsilon_1, \quad \psi_{\varepsilon}=\psi_0+\varepsilon_2.
\end{aligned} 
\end{equation*} 
Note that the position of the local stable and local unstable manifolds of the point $O$ (in a small neighborhood independent of $\rho$) does not change as $\varepsilon$ varies. Therefore, the manifold $\mathcal W^{s}_{M_j}(p_j)$ remains straightened and the manifold $\mathcal W^{u}_{M_j}(q_j)$ preserves the form \eqref{lemmasplittingoftangencyeq1}.
\end{demo}

\begin{remark}\label{remarkperturbationofcoefficients}
In the same way, one can construct finite-parameter families such that any finite set of the Taylor coefficients of the local return maps near any set of periodic orbits and of the Taylor coefficients of the equations of any small pieces of the stable and unstable invariant manifolds change independently as the parameters vary. This imlpies that the genericity conditions that we require throughout the paper can indeed be achieved in every regularity class.
\end{remark}

\section{Homoclinic tangencies of all orders}\label{section3}

In this section, we prove Theorem \ref{theorem4}. As it is explained in the introduction, Theorems \ref{theorem1}--\ref{theorem3} follow from Theorem \ref{theorem4} immediately.

\subsection{Corank-1 homoclinic tangencies}\label{subsectionhighordertangencies}

\begin{lemma}\label{lemma1highordertangencies}
Let $f_0:\mathcal M^d\rightarrow \mathcal M^d$ be a $C^r$ diffeomorphism, where $d\ge 2$ and $r=4,\ldots,\infty,\omega$. Let $f_0$ have a saddle periodic orbit $L$ and two orbits $\Gamma,\widetilde\Gamma$ of corank-1 homoclinic tangency between $W^s(L)$ and $W^u(L)$. Let $L$ satisfy condition \textbf{\hyperlink{C1.}{C1}}, and $\Gamma,\widetilde\Gamma$ satisfy conditions \textbf{\hyperlink{C2.}{C2}},\textbf{\hyperlink{C3.}{C3}}. Let $\Gamma$ be of order $n$ $(n\le r-3)$, and $\widetilde\Gamma$ be quadratic. Suppose $f_0$ is embedded into an $\ell$-parameter $(\ell\ge n+1)$ $C^r$ family $(f_{\zeta,\varepsilon})$ with $\zeta\in\mathbb R^{\ell-n-1},\varepsilon\in\mathbb R^{n+1}$ such that 
\begin{equation}\label{lemma1highordertangencieseq0}
\begin{aligned} 
&\left.\det\frac{\partial(\eta_0,\ldots,\eta_{n-1},\tilde\eta_0)}{\partial\varepsilon}\right|_{(\zeta,\varepsilon)=0}\not=0,
\end{aligned}
\end{equation}
where $\eta_0,\ldots,\eta_{n-1}$ are the splitting coefficients of orders $0,\ldots,n-1$ for $\Gamma$, and $\tilde\eta_0$ is the splitting coefficient of order $0$ for $\widetilde\Gamma$. Then there exists a sequence $\langle\Delta_\alpha\rangle_{\alpha\in\mathbb N}$ of $(\ell-n-1)$-dimensional, centered at $\zeta=0$ balls and a sequence $\langle\varepsilon_\alpha(\zeta)\rangle_{\alpha\in\mathbb N}\to0$ of $C^{r-n-1}$ functions defined for $\zeta\in\Delta_{\alpha}$ such that the diffeomorphism $f_{\zeta,\varepsilon_{\alpha}(\zeta)}$ has, for all $\zeta\in\Delta_{\alpha}$, an orbit of corank-1 homoclinic tangency of order $(n+1)$ between $W^s(L)$ and $W^u(L)$. For each $\alpha$, this orbit depends on $\zeta$ continuously and the tangency splits freely up to order $n$ as $\varepsilon$ varies.
\end{lemma}

\begin{demo}
Let $O$ be a point of the periodic orbit $L$, and let $T_{O}$ be the local map defined in a small neighborhood $U(O)$. By Lemma \ref{lemmachoiceofcoordinates}, we choose $(T_{O},\Gamma,\widetilde\Gamma)$-adapted coordinates in $U(O)$. In these coordinates, for all small $(\zeta,\varepsilon)$, the manifolds $\mathcal W^s(O)$ and $\widetilde{\mathcal W}^u(O)$ are straightened, while the manifolds $\mathcal W^u(O)$ and $\widetilde{\mathcal W}^s(O)$ are written as 
\begin{equation}\label{lemma1highordertangencieseq0.25}
\mathcal W^u(O): \;\;
\begin{aligned} 
y_1&=\sum\limits_{i=0}^{n} \eta_{i}\cdot (x_1-x_1^+)^i+a\cdot (x_1-x_1^+)^{n+1}+b\cdot y_2+\\
&+o(||x_1-x_1^+||^{n+1})+\theta_1(x_1-x_1^+,y_2,v,\zeta,\varepsilon),\\
x_2- x_2^+&=\theta_2(x_1-x_1^+,y_2,v,\zeta,\varepsilon),\\
u-u^+&=\theta_3(x_1-x_1^+,y_2,v,\zeta,\varepsilon),
\end{aligned}
\end{equation}
and 
\begin{equation}\label{lemma1highordertangencieseq0.5}
\widetilde{\mathcal W}^s(O): \;\;
\begin{aligned} 
\tilde x_1&=\tilde\eta_0+\tilde\eta_{1}\cdot (\tilde y_1-\tilde y_1^-)+\tilde a\cdot (\tilde y_1-\tilde y_1^-)^{2}+\tilde b\cdot x_2+\\
&+o((\tilde y_1-\tilde y_1^-)^{2})+\tilde\theta_1(\tilde x_2,\tilde u,\tilde y_1-\tilde y_1^-,\zeta,\varepsilon),\\
\tilde y_2-\tilde y_2^-&=\tilde\theta_2(\tilde x_2,\tilde u,\tilde y_1-\tilde y_1^-,\zeta,\varepsilon),\\
\tilde v-\tilde v^-&=\tilde\theta_3(\tilde x_2,\tilde u,\tilde y_1-\tilde y_1^-,\zeta,\varepsilon),
\end{aligned}
\end{equation}
where the functions $\theta$ and $\tilde\theta$ satisfy \eqref{lemmachoiceofcoordinatesestimateseq2}; $a, \tilde a$ are non-vanishing functions of $\zeta$; $\eta_0,\ldots,\eta_{n}$ and $\tilde\eta_0,\tilde\eta_1$ are the splitting coefficients for the orbits $\Gamma$ and $\widetilde\Gamma$, respectively. By \eqref{lemma1highordertangencieseq0}, we may take $\eta_0,\ldots,\eta_{n-1}$ and $\tilde\eta_0$ as the new independent parameters, i.e., we assume $\varepsilon=\eta:=(\eta_0,\ldots,\eta_{n-1},\tilde\eta_0)$ without loss of generality. 

By Corollary \ref{corollarylocalmap}, there exists a sequence $\langle k_\alpha\rangle_{\alpha\in\mathbb N}$ such that $k_\alpha\xrightarrow{}+\infty$ and for any $\langle k_{\alpha}\rangle$-tame sequence $\langle\mathfrak{D}_{\alpha}\rangle_{\alpha\in\mathbb N}$ and all $(\zeta,\eta)\in\mathfrak{D}_{\alpha}$, the equations \eqref{eqlocal3tamesequence} are satisfied. Further, we will show that there exists a $\langle k_{\alpha}\rangle$-tame sequence $\langle\mathfrak{D}_{\alpha}\rangle_{\alpha\in\mathbb N}$, a sequence $\langle\Delta_\alpha\rangle_{\alpha\in\mathbb N}$ of centered at $\zeta=0$ balls, and a sequence $\langle\varepsilon_\alpha(\zeta)\rangle_{\alpha\in\mathbb N}\to0$ of $C^{r-n-1}$ functions defined on $\zeta\in\Delta_{\alpha}$ such that $(\zeta,\varepsilon_\alpha(\zeta))\in\mathfrak D_{\alpha}$ and, for all sufficiently large $\alpha$, the manifolds $\mathcal W^u(O)$ and $T_{O}^{-k_\alpha}(\widetilde{\mathcal W}^s(O))$ have a corank-1 tangency of order $n+1$ for $\eta=\varepsilon_{\alpha}(\zeta)$. This will prove the lemma.

Substituting formulas \eqref{eqlocal3tamesequence} for the local map into \eqref{lemma1highordertangencieseq0.25} and \eqref{lemma1highordertangencieseq0.5}, we conclude that a corank-1 tangency of order $(n+1)$ between the manifolds $\mathcal W^u(O)$ and $T_{O}^{-k_{\alpha}}(\widetilde{\mathcal W}^s(O))$ exists if and only if there is a corank-1 tangency of order $(n+1)$ between the following manifolds:
\begin{equation}\label{lemma1highordertangencieseq1}
\begin{aligned} 
&\gamma^{-k_{\alpha}}\cdot (Y_1+\tilde y_1^-)=\sum\limits_{i=0}^{n} \eta_{i}\cdot X_1^i+a\cdot X_1^{n+1}+\gamma^{-k_{\alpha}}\cdot b\cdot(Y_2+\tilde y_2^-)+\Theta_1(X_1,Y,V,\zeta,\eta),\\
&X_2=\Theta_2(X_1,Y,V,\zeta,\eta), \quad U=\Theta_3(X_1,Y,V,\zeta,\eta),
\end{aligned}
\end{equation}
and
\begin{equation}\label{lemma1highordertangencieseq2}
\begin{aligned} 
&\lambda^{k_{\alpha}}\cdot (X_1+x_1^+)=\tilde\eta_0+ \tilde\eta_{1}\cdot Y_1+\tilde a\cdot Y_1^2+\lambda^{k_{\alpha}}\cdot \tilde b\cdot(X_2+x_2^+)+\widetilde\Theta_1(X,U,Y_1,\zeta,\eta),\\
&Y_2=\widetilde\Theta_2(X,U,Y_1,\zeta,\eta), \quad V=\widetilde\Theta_3(X,U,Y_1,\zeta,\eta),
\end{aligned}
\end{equation}
where 
\begin{equation}\label{lemma1highordertangencieseq2.5}
\begin{aligned} 
&X=x-x^+, \qquad U=u-u^+, \qquad Y=\tilde y-\tilde y^-, \qquad V=\tilde v-\tilde v^-,\\
\end{aligned}
\end{equation}
and $\Theta, \widetilde\Theta$ are some functions that satisfy
\begin{equation*}
\begin{aligned} 
\Theta_1&=o(X_1^{n+1})+o(\gamma^{-k_{\alpha}})\cdot(Y,V)+\gamma^{-k_{\alpha}} \cdot O(X_1\cdot Y),\\
\Theta_2&=O(X_1^2)+O(\gamma^{-k_{\alpha}}), \quad\quad \Theta_3=O(X_1)+O(\gamma^{-k_{\alpha}}),\\
\widetilde\Theta_1&=o(Y_1^{2})+o(\lambda^{k_{\alpha}})\cdot(X,U)+\lambda^{k_{\alpha}}\cdot O(Y_1\cdot X),\\
\widetilde\Theta_2&=O(Y_1^2)+O(\lambda^{k_{\alpha}}), \quad\quad\;\;\; \widetilde\Theta_3=O(Y_1)+O(\lambda^{k_{\alpha}}).
\end{aligned}
\end{equation*}
We need to nullify the constant terms in the functions $\Theta_{2,3}, \widetilde\Theta_{2,3}$, as well as the coefficient of $X_1^n$ in the first equation of \eqref{lemma1highordertangencieseq1}, and the coefficient of $Y_1$ in the first equation of \eqref{lemma1highordertangencieseq2}. This can be done via the shift
\begin{equation*}
\begin{aligned} 
X_1&=X_1^{new}+O(\eta_n)+o(\gamma^{-k_{\alpha}}), \quad &Y_1&=Y_1^{new}+O(\tilde\eta_1)+o(\lambda^{k_{\alpha}}),\\
X_2&=X_2^{new}+O(\eta^2_n)+O(\gamma^{-k_{\alpha}}), \quad &Y_2&=Y_2^{new}+O(\tilde\eta^2_1)+O(\lambda^{k_{\alpha}}),\\
U&=U^{new}+O(\eta_n)+O(\gamma^{-k_{\alpha}}), \quad &V&=V^{new}+O(\tilde\eta_1)+O(\lambda^{k_{\alpha}}).
\end{aligned}
\end{equation*}
As a result, we bring the equations \eqref{lemma1highordertangencieseq1} and \eqref{lemma1highordertangencieseq2} to the form
\begin{equation}\label{lemma1highordertangencieseq4}
\begin{aligned} 
&\gamma^{-k_{\alpha}}\cdot (Y_1-b\cdot Y_2)=\sum\limits_{i=0}^{n-1} \eta^*_{i}\cdot X_1^i+a\cdot X_1^{n+1}+\Theta^*_1(X_1,Y,V,\zeta,\eta^*),\\
&X_2=\Theta^*_2(X_1,Y,V,\zeta,\eta^*), \quad U=\Theta^*_3(X_1,Y,V,\zeta,\eta^*),
\end{aligned}
\end{equation}
and 
\begin{equation}\label{lemma1highordertangencieseq5}
\begin{aligned} 
&\lambda^{k_{\alpha}}\cdot(X_1-\tilde b\cdot X_2)=\tilde\eta^*_0+\tilde a\cdot Y_1^2+\widetilde\Theta^*_1(X,U,Y_1,\zeta,\eta^*),\\
&Y_2=\widetilde\Theta^*_2(X,U,Y_1,\zeta,\eta^*), \quad V=\widetilde\Theta^*_3(X,U,Y_1,\zeta,\eta^*),
\end{aligned}
\end{equation}
where 
\begin{equation*}
\begin{aligned} 
\eta^*_0&=\eta_0+\gamma^{-k_{\alpha}}\cdot(b\cdot\tilde y_2^--\tilde y_1^-+O(\eta^*))+o(\gamma^{-k_{\alpha}}),\\
\eta^*_i&=\eta_i+O(\gamma^{-k_{\alpha}}) \quad \text{for} \quad i\in\{1,\ldots,n-1\},\\
\tilde\eta^*_0&=\tilde\eta_0+\lambda^{k_{\alpha}}\cdot(\tilde b\cdot x_2^+-x_1^++O(\eta^*))+o(\lambda^{k_{\alpha}}),
\end{aligned}
\end{equation*}
and $\Theta^*, \widetilde\Theta^*$ are functions satisfying
\begin{equation*}
\begin{aligned} 
\Theta^*_1&=o(X_1^{n+1})+o(\gamma^{-k_{\alpha}})\cdot(Y,V)+\gamma^{-k_{\alpha}} \cdot O(X_1\cdot Y),\\
\Theta^*_2&=O(X_1^2)+\gamma^{-k_{\alpha}}\cdot O(X_1,Y,V), \quad\quad \Theta^*_3=O(X_1)+\gamma^{-k_{\alpha}}\cdot O(X_1,Y,V),\\
\widetilde\Theta^*_1&=o(Y_1^{2})+o(\lambda^{k_{\alpha}})\cdot(X,U)+\lambda^{k_{\alpha}}\cdot O(Y_1\cdot X),\\
\widetilde\Theta^*_2&=O(Y_1^2)+\lambda^{k_{\alpha}}\cdot O(X,U,Y_1), \quad\quad\;\;\; \widetilde\Theta^*_3=O(Y_1)+\lambda^{k_{\alpha}}\cdot O(X,U,Y_1).
\end{aligned}
\end{equation*}
Let us rescale the variables $(X,U,Y,V)$ as
\begin{equation*}
\begin{aligned} 
X_1&=-\frac{1}{(\tilde a \cdot a^2)^{\frac{1}{2n+1}}}\cdot\lambda^{\frac{k_{\alpha}}{2n+1}}\cdot\gamma^{-\frac{2k_{\alpha}}{2n+1}}\cdot X_1^{new}, &Y_1&=\frac{(-1)^{n+1}}{(\tilde a^{n+1}\cdot a)^{\frac{1}{2n+1}}}\cdot\lambda^{\frac{k_{\alpha}(n+1)}{2n+1}}\cdot\gamma^{-\frac{k_{\alpha}}{2n+1}}\cdot Y_1^{new},\\
X_2&=\varrho_{k_\alpha}\cdot\lambda^{\frac{k_\alpha}{2n+1}}\cdot\gamma^{-\frac{2k_\alpha}{2n+1}}\cdot X_2^{new} &Y_2&=\varrho_{k_\alpha}\cdot\lambda^{\frac{k_\alpha(n+1)}{2n+1}}\cdot\gamma^{-\frac{k_\alpha}{2n+1}}\cdot Y_2^{new},\\
U&=\frac{1}{\varrho_{k_{\alpha}}}\cdot\lambda^{\frac{k_{\alpha}}{2n+1}}\cdot\gamma^{-\frac{2k_{\alpha}}{2n+1}}\cdot U^{new}, &V&=\frac{1}{\varrho_{k_{\alpha}}}\cdot\lambda^{\frac{k_{\alpha}(n+1)}{2n+1}}\cdot\gamma^{-\frac{k_{\alpha}}{2n+1}}\cdot V^{new},
\end{aligned}
\end{equation*}
where $\varrho_{k_{\alpha}}$ tends sufficiently slowly to zero as $\alpha\rightarrow+\infty$. We also rescale the parameters $\eta^*_0,\ldots,\eta^*_{n-1},\tilde\eta^*_0$ as
\begin{equation}\label{lemma1highordertangencieseq7.5}
\begin{aligned} 
\eta^*_i&=\frac{(-1)^{n+1-i}}{(\tilde a\cdot a^2)^{\frac{n+1-i}{2n+1}}}\cdot a\cdot\lambda^{\frac{k_{\alpha}(n+1-i)}{2n+1}}\cdot\gamma^{-\frac{2k_{\alpha}(n+1-i)}{2n+1}}\cdot H_i \quad \text{for} \quad i\in\{0,\ldots,n-1\},\\
\tilde\eta^*_0&=\frac{1}{(\tilde a \cdot a^2)^{\frac{1}{2n+1}}}\cdot\lambda^{\frac{2k_{\alpha}(n+1)}{2n+1}}\cdot\gamma^{-\frac{2k_{\alpha}}{2n+1}}\cdot \widetilde H_0.
\end{aligned}
\end{equation}
As a result, the equations \eqref{lemma1highordertangencieseq4} and \eqref{lemma1highordertangencieseq5} change to
\begin{equation}\label{lemma1highordertangencieseq8}
\begin{aligned} 
Y_1&=\sum\limits_{i=0}^{n-1} H_{i}\cdot X_1^i+X_1^{n+1}+\vartheta_{1,k_{\alpha}}(X_1,Y,V,\zeta,H),\\
(X_2,U)&= \vartheta_{2,k_{\alpha}}(X_1,Y,V,\zeta,H),\\
\end{aligned}
\end{equation}
and
\begin{equation}\label{lemma1highordertangencieseq9}
\begin{aligned} 
X_1&=\widetilde H_0-Y_1^2+\tilde\vartheta_{1,k_{\alpha}}(X,U,Y_1,\zeta,H),\\
(Y_2,V)&=\tilde\vartheta_{2,k_{\alpha}}(X,U,Y_1,\zeta,H),\\
\end{aligned}
\end{equation}
where the functions $\vartheta_{k_{\alpha}},\tilde\vartheta_{k_{\alpha}}$ and their derivatives (including at least one derivative with respect to $\zeta$ and $H$) tend to zero uniformly on any compact set of the rescaled variables $(X,U,Y,V)$ and of the rescaled parameters $H=\big(H_0,\ldots,H_{n-1},\widetilde H_0\big)$ as $\alpha\rightarrow+\infty$. Making use of the Implicit Function Theorem, we rewrite the equations \eqref{lemma1highordertangencieseq8} and \eqref{lemma1highordertangencieseq9} as
\begin{equation}\label{lemma1highordertangencieseq10}
\begin{aligned} 
Y_1&=\sum\limits_{i=0}^{n-1} H_{i}\cdot X_1^i+X_1^{n+1}+\vartheta^*_{1,k_{\alpha}}(X_1,Y_2,V,\zeta,H),\\
(X_2,U)&= \vartheta^*_{2,k_{\alpha}}(X_1,Y_2,V,\zeta,H),\\
\end{aligned}
\end{equation}
and
\begin{equation}\label{lemma1highordertangencieseq11}
\begin{aligned} 
X_1&=\widetilde H_0-Y_1^2+\tilde\vartheta^*_{1,k_{\alpha}}(X_2,U,Y_1,\zeta,H),\\
(Y_2,V)&=\tilde\vartheta^*_{2,k_{\alpha}}(X_2,U,Y_1,\zeta,H),\\
\end{aligned}
\end{equation}
where the functions $\vartheta^*_{k_{\alpha}},\tilde\vartheta^*_{k_{\alpha}}$, along with their derivatives, tend to zero as $\alpha\rightarrow+\infty$. 

We recast \eqref{lemma1highordertangencieseq11} as
\begin{equation}\label{lemma1highordertangencieseq13}
\begin{aligned} 
Y_1&=\pm\sqrt{\widetilde H_0-X_1}+\tilde\vartheta^{**}_{1,k_{\alpha}}(X_2,U,\zeta,H),\\
\left(Y_2,V\right)&=\tilde\vartheta^{**}_{2,k_{\alpha}}(X,U,\zeta,H),\\
\end{aligned}
\end{equation}
where the functions $\tilde\vartheta^{**}_{k_{\alpha}}$, together with their derivatives, decay to zero as $\alpha\rightarrow+\infty$, uniformly in the region where $\widetilde H_0-X_1$ is bounded away from zero. Further, we will consider only the "minus" branch of this equation. Let us make the following change of variables
\begin{equation*}
\begin{aligned} 
Y_1^{new}&=Y_1+\sqrt{\widetilde H_0-X_1}-\tilde\vartheta^{**}_{1,k_{\alpha}}(X_2,U,\zeta,H),\\
(Y_2,V)^{new}&=(Y_2,V)-\tilde\vartheta^{**}_{2,k_{\alpha}}(X,U,\zeta,H).\\
\end{aligned}
\end{equation*}
As a result, the equations \eqref{lemma1highordertangencieseq10} and \eqref{lemma1highordertangencieseq13} take the form 
\begin{equation}\label{lemma1highordertangencieseq14}
\begin{aligned} 
Y_1&=\sqrt{\widetilde H_0-X_1}+\sum\limits_{i=0}^{n-1} H_{i}\cdot X_1^i+X_1^{n+1}+\vartheta^{**}_{1,k_{\alpha}}(X_1,Y_2,V,\zeta,H),\\
\left(X_2,U\right)&= \vartheta^{**}_{2,k_{\alpha}}(X_1,Y_2,V,\zeta,H),\\
\end{aligned}
\end{equation}
and
\begin{equation}\label{lemma1highordertangencieseq15}
\begin{aligned} 
(Y,V)&=0,\\
\end{aligned}
\end{equation}
where the functions $\vartheta^{**}_{k_{\alpha}}$, along with their derivatives, decay to zero as $\alpha\rightarrow+\infty$. 

By Definition \ref{coordinatedefinitionorder}, the manifolds given by \eqref{lemma1highordertangencieseq14} and \eqref{lemma1highordertangencieseq15} have a corank-1 tangency of order $(n+1)$ if one can choose the parameters $H_0,...,H_{n-1},\widetilde H_0$ and the $X_1$-coordinate $\widehat X_1$ of the point of the sought tangency in such a way that the first $(n+2)$ terms $\nu_0,\ldots,\nu_{n+1}$ of the Taylor expansion at $\widehat X_1$ for the first equation in \eqref{lemma1highordertangencieseq14} vanish. They are given by
\begin{equation}\label{lemma1highordertangencieseq15.5}
\begin{aligned} 
\nu_0&=\sqrt{\widetilde H_0-\widehat X_1}+\sum\limits_{i=0}^{n-1} H_{i}\cdot \widehat X_1^i+\widehat X_1^{n+1}+o(1)_{\alpha\rightarrow+\infty},\\
\nu_j&=-\frac{\left(2j-3\right)!!}{2^j\cdot j!}\cdot(\widetilde H_0-\widehat X_1)^{\frac{1}{2}-j}+\sum\limits_{i=j}^{n-1}\cdot\binom{i}{j}\cdot H_i\cdot\widehat X_1^{i-j}+\\
&+\binom{n+1}{j}\cdot\widehat X_1^{n+1-j}+o(1)_{\alpha\rightarrow+\infty} \;\; \text{for} \;\; j\in\{1,\ldots,n-1\},\\
\nu_n&=-\frac{\left(2n-3\right)!!}{2^n\cdot n!}\cdot (\widetilde H_0-\widehat X_1)^{\frac{1}{2}-n}+(n+1)\cdot \widehat X_1+o(1)_{\alpha\rightarrow+\infty},\\
\nu_{n+1}&=-\frac{\left(2n-1\right)!!}{2^{n+1}\cdot\left(n+1\right)!}\cdot (\widetilde H_0-\widehat X_1)^{-\frac{1}{2}-n}+1+o(1)_{\alpha\rightarrow+\infty}.
\end{aligned}
\end{equation}
For $\widetilde H_0>\widehat X_1>0$, the first $n$ terms can obviously be nullified by the choice of $H_0,...,H_{n-1}$ since the first $n$ equations constitute a linear triangular system for $H_0,...,H_{n-1}$. The vanishing of the remaining two terms is achieved by
\begin{equation}\label{lemma1highordertangencieseq15.75}
\widetilde H_0=\left(\frac{2n+1}{2}\right)\cdot\widehat X_1+o(1)_{\alpha\rightarrow+\infty}, \quad \widehat X_1=\frac{2}{2n-1}\cdot\left(\frac{\left(2n-1\right)!!}{2^{n+1}\cdot\left(n+1\right)!}\right)^{\frac{2}{2n+1}}+o(1)_{\alpha\rightarrow+\infty}.
\end{equation}
Note that $\widetilde H_0-\widehat X_1>0$ here, uniformly for all large $\alpha$, as required.

Thus, we proved the existence of a corank-1 tangency of order $(n+1)$ between the manifolds \eqref{lemma1highordertangencieseq1} and \eqref{lemma1highordertangencieseq2}. This tangency splits freely up to order $n$ as $H$ varies: by \eqref{lemma1highordertangencieseq15.5}),
\begin{equation}\label{lemma1highordertangencieseq16}
\det\frac{\partial (\nu_0,\ldots,\nu_n)}{\partial H}\not=0.
\end{equation}

The corresponding values of $H_0,...,H_{n-1},\widetilde H_0$ are uniformly bounded for all large $\alpha$. Since the original parameters $\eta=(\eta_0,...,\eta_{n-1},\tilde\eta_0)$ are related to $H=\big(H_0,...,H_{n-1},\widetilde H_0\big)$ by a scaling diffeomorphism \eqref{lemma1highordertangencieseq7.5}, we get that $\varepsilon_{\alpha}(\zeta)$ corresponding to the found tangency of order $(n+1)$ tend to zero as $\alpha\rightarrow+\infty$. They decay exponentially fast, so there exists a $\langle k_{\alpha}\rangle$-tame sequence $\langle\mathfrak{D}_{\alpha}\rangle_{\alpha\in\mathbb N}$ and a sequence of balls $\langle\Delta_\alpha\rangle_{\alpha\in\mathbb N}$ such that, for all sufficiently large $\alpha$, one has $(\zeta,\varepsilon_{\alpha}(\zeta))\in\mathfrak{D}_{\alpha}$. The balls $\Delta_\alpha$ are centered at $\zeta=0$ because $\mathfrak{D}_{\alpha}$ are centered at zero, by Definition \ref{definitiontamesequence}.
\end{demo}

\begin{remark}\label{remarkhighordertangencies}
By \eqref{lemma1highordertangencieseq15.5},\eqref{lemma1highordertangencieseq15.75}, we can always choose $\widehat X_1$ such that $\nu_{n+1}=0$ for all parameter values in s small neighborhood $(\zeta,\varepsilon_\alpha(\zeta))$. In the later use (see Lemma \ref{lemmakeepingparatangencyorder}), this will be our choice. It also follows from \eqref{lemma1highordertangencieseq15.5},\eqref{lemma1highordertangencieseq15.75} that the rescaled splitting coefficients $H$ are smooth functions of the remaining splitting coefficients $\nu=(\nu_0,\ldots,\nu_{n})$ for the newly found homoclinic tangency of order $(n+1)$ for this choice of $\widehat X_1$. Therefore, by \eqref{lemma1highordertangencieseq7.5}, one can find the following relation between the splitting coefficients $\eta=(\eta_0,\ldots,\eta_{n-1},\tilde\eta_0)$ of the original homoclinic tangencies $\Gamma,\widetilde\Gamma$ and $\nu$:
\begin{equation}\label{remarkhighordertangencieseq1}
\begin{aligned} 
\eta(\zeta,\varepsilon)&=Q_{\alpha}(\zeta,\nu),
\end{aligned}
\end{equation}
where the functions $Q_{\alpha}$ are such that  $Q_{\alpha}\to 0$ in $C^{r-n-1}$ and, by \eqref{lemma1highordertangencieseq16}, one has $\left.\det\frac{\partial Q_{\alpha}}{\partial\nu}\right|_{(\zeta,\nu)=0}\not=0$.
\end{remark}

\begin{lemma}\label{lemma2highordertangencies}
Let $f:\mathcal M^d\rightarrow \mathcal M^d$ be a $C^r$ diffeomorphism, where $d\ge 2$ and $r=2,.\ldots,\infty,\omega$. Let $f$ have a saddle periodic orbit $L$ and $n$ orbits $\Gamma_1,\ldots,\Gamma_n$ of corank-1 quadratic homoclinic tangency between $W^s(L)$ and $W^u(L)$, where $n\le r-1$. Then, $C^r$-close to $f$, there exists a diffeomorphism with an orbit of corank-1 homoclinic tangency of order $n$ between $W^s(L)$ and $W^u(L)$.
\end{lemma}

\begin{demo} 
For $n=1$, the lemma holds trivially. Assume $n>1$. Using Lemma \ref{lemmasplittingoftangency}, we embed $f$ into a finite-parameter $C^r$ family of diffeomorphisms such that the tangencies $\Gamma_1,\Gamma_2$ split freely and independently up to order $0$, whereas the remaining $n-2$ tangencies do not split. According to Lemma \ref{lemma1highordertangencies}, in this family, there exists a diffeomorphism $f_1$ which is $C^r$ close to $f$ and has an orbit $\Gamma_{1,2}$ of corank-1 homoclinic tangency of order $2$ (cubic) between $W^s(L)$ and $W^u(L)$. Thus, if $n=2$, the lemma is proved. 

When $n>2$, we continue by embedding $f_1$ into a finite-parameter $C^r$ family of diffeomorphisms such that the tangencies $\Gamma_{1,2}$ and $\Gamma_3$ split freely and independently up to orders $1$ and $0$, respectively, whereas the remaining $n-3$ tangencies do not split. By Lemma \ref{lemma1highordertangencies}, in this family, there exists a diffeomorphism $f_2$ which is $C^r$ close to $f_1$ and possesses an orbit $\Gamma_{1,2,3}$ of corank-1 homoclinic tangency of order $3$ between $W^s(L)$ and $W^u(L)$. Thus, if $n=3$, the lemma is proved. 

By repeating this procedure, we construct a diffeomorphism $f_{n-1}$ which is $C^r$-close to $f$ and possesses an orbit $\Gamma_{1,\ldots,n}$ of corank-1 homoclinic tangency of order $n$ between $W^s(L)$ and $W^u(L)$.
\end{demo}

\subsection{Corank-2 homoclinic tangencies}\label{subsectioncorank2tangencies}

\begin{lemma}\label{lemmacorank2tangency}
Let $f_0:\mathcal M^d\rightarrow \mathcal M^d$ be a $C^r$ diffeomorphism, where $d\ge 4$ and $r=2,\ldots,\infty,\omega$. Let $f_0$ have a bi-focus periodic orbit $L$ and two orbits $\Gamma,\widetilde\Gamma$ of corank-1 quadratic homoclinic tangency between $W^s(L)$ and $W^u(L)$. Let $L$ satisfy condition \textbf{\hyperlink{C3*.}{C3*}}, and $\Gamma,\widetilde\Gamma$ satisfy condition \textbf{\hyperlink{C2.}{C2}}. Suppose $f_0$ is embedded into a four-parameter $C^r$ family $(f_{\varepsilon})$ such that 
\begin{equation}\label{lemmacorank2tangencyeq0}
\begin{aligned} 
&\left.\det\frac{\partial(\eta_0,\tilde\eta_0,\varphi_{\varepsilon},\psi_{\varepsilon})}{\partial\varepsilon}\right|_{\varepsilon=0}\not=0,
\end{aligned}
\end{equation}
where $\eta_0,\tilde\eta_0$ are the splitting coefficients of order $0$ for $\Gamma,\widetilde\Gamma$, respectively, and $\varphi_{\varepsilon},\psi_{\varepsilon}$ are the arguments of the leading stable and leading unstable multipliers of $L$, respectively. Then there exists a sequence $\langle\varepsilon_k\rangle_{k\in\mathbb N}$ with $\varepsilon_{k}\xrightarrow{}0$ such that $f_{\varepsilon_{k}}$ has an orbit of corank-2 homoclinic tangency between $W^s(L)$ and $W^u(L)$.
\end{lemma}

\begin{demo}
First, we choose a point $O$ in the periodic orbit $L$, and define the local map $T_{O}$ in a small neighborhood $U(O)$. By Lemma \ref{lemmachoiceofcoordinatesandhomoclinicpoints}, we choose $(T_{O},\Gamma,\widetilde\Gamma)$-adapted coordinates in $U(O)$. In these coordinates, for all small $\varepsilon$, the manifolds $\mathcal W^s(O)$ and $\widetilde{\mathcal W}^u(O)$ are given by \eqref{lemmachoiceofcoordinateseq0}, while the manifolds $\mathcal W^u(O)$ and $\widetilde{\mathcal W}^s(O)$ can be written as 
\begin{equation}\label{lemmacorank2tangencyeq0.25}
\mathcal W^u(O): \;\;
\begin{aligned} 
y_1&=\eta_0+\eta_{1}\cdot (x_1-x_1^+)+a\cdot(x_1-x_1^+)^{2}+b\cdot y_2+\\
&+o((x_1-x_1^+)^{2})+\theta_1(x_1-x_1^+,y_2,v,\varepsilon),\\
x_2- x_2^+&=\theta_2(x_1-x_1^+,y_2,v,\varepsilon),\\
u-u^+&=\theta_3(x_1-x_1^+,y_2,v,\varepsilon),
\end{aligned}
\end{equation}
and 
\begin{equation}\label{lemmacorank2tangencyeq0.5}
\widetilde{\mathcal W}^s(O): \;\;
\begin{aligned} 
\tilde x_1&=\tilde\eta_0+\tilde\eta_{1}\cdot (\tilde y_1-\tilde y_1^-)+\tilde a\cdot(\tilde y_1-\tilde y_1^-)^{2}+\tilde b\cdot x_2+\\
&+o((\tilde y_1-\tilde y_1^-)^{2})+\tilde\theta_1(\tilde x_2,\tilde u,\tilde y_1-\tilde y_1^-,\varepsilon),\\
\tilde y_2-\tilde y_2^-&=\tilde\theta_2(\tilde x_2,\tilde u,\tilde y_1-\tilde y_1^-,\varepsilon),\\
\tilde v-\tilde v^-&=\tilde\theta_3(\tilde x_2,\tilde u,\tilde y_1-\tilde y_1^-,\varepsilon),
\end{aligned}
\end{equation}
where the functions $\theta$ and $\tilde\theta$ satisfy \eqref{lemmachoiceofcoordinatesestimateseq2}; $a,\tilde a$ are non-zero numbers, $\eta_0,\eta_1$ and $\tilde\eta_0,\tilde\eta_1$ are the splitting coefficients for the orbits $\Gamma$ and $\widetilde\Gamma$, respectively. By condition \eqref{lemmacorank2tangencyeq0}, we may choose $\eta_0,\tilde\eta_0$, and $\widehat\varphi=\varphi_{\varepsilon}-\varphi_0,\widehat\psi=\psi_{\varepsilon}-\psi_0$ as the new independent parameters, so we assume $\varepsilon=(\eta_0,\tilde\eta_0,\widehat\varphi,\widehat\psi)$ without loss of generality. Below, we will prove that there exists a sequence $\langle\varepsilon_k\rangle_{k\in\mathbb N}$ with $\varepsilon_{k}\xrightarrow{}0$ such that, for all sufficiently large $k$, the manifolds $\mathcal W^u(O)$ and $T_{O}^{-k}(\widetilde{\mathcal W}^s(O))$ have a corank-2 tangency. This will prove the lemma.

By combining formulas \eqref{eqlocal3} for the local map and equations \eqref{lemmacorank2tangencyeq0.25},\eqref{lemmacorank2tangencyeq0.5}, we deduce that a corank-2 tangency between the manifolds $\mathcal W^u(O)$ and $T_{O}^{-k}(\widetilde{\mathcal W}^s(O))$ exists if and only if there is a corank-2 tangency between the manifolds defined by formulas \eqref{lemmacorank2tangencyeq1} and \eqref{lemmacorank2tangencyeq2} below:
\begin{equation}\label{lemmacorank2tangencyeq1}
\begin{aligned} 
&\gamma_{\varepsilon}^{-k}\cdot (\cos(k\psi_{\varepsilon})\cdot(Y_1+\tilde y_1^-)+\sin(k\psi_{\varepsilon})\cdot (Y_2+\tilde y_2^-))=\eta_0+\eta_1\cdot X_1+a\cdot X_1^2+\\
&+\gamma_{\varepsilon}^{-k}\cdot b\cdot (-\sin(k\psi_{\varepsilon})\cdot(Y_1+\tilde y_1^-)+\cos(k\psi_{\varepsilon})\cdot (Y_2+\tilde y_2^-))+\Theta_1(X,U,Y,V,\varepsilon),\\
&X_2=\Theta_2(X,U,Y,V,\varepsilon), \quad U=\Theta_3(X,U,Y,V,\varepsilon),
\end{aligned}
\end{equation}
and
\begin{equation}\label{lemmacorank2tangencyeq2}
\begin{aligned} 
&\lambda_{\varepsilon}^{k}\cdot(\cos(k\varphi_{\varepsilon})\cdot(X_1+x_1^+)-\sin(k\varphi_{\varepsilon})\cdot (X_2+x_2^+))=\tilde\eta_0+ \tilde\eta_{1}\cdot Y_1+\tilde a\cdot Y_1^2+\\
&+\lambda_{\varepsilon}^{k}\cdot \tilde b\cdot(\sin(k\varphi_{\varepsilon})\cdot(X_1+x_1^+)+\cos(k\varphi_{\varepsilon})\cdot (X_2+x_2^+))+\widetilde\Theta_1(X,U,Y,V,\varepsilon),\\
&Y_2=\widetilde\Theta_2(X,U,Y,V,\varepsilon), \quad V=\widetilde\Theta_3(X,U,Y,V,\varepsilon),
\end{aligned}
\end{equation}
where $X,U,Y,V$ are given by \eqref{lemma1highordertangencieseq2.5}, and $\Theta, \widetilde\Theta$ are some functions satisfying 
\begin{equation*}
\begin{aligned} 
\Theta_1&=o(X_1^2)+o(\gamma_{\varepsilon}^{-k})+\gamma_{\varepsilon}^{-k}\cdot O(X_1),\\
\Theta_2&=O(X_1^2)+O(\gamma_{\varepsilon}^{-k}), \quad\quad \Theta_3=O(X_1)+O(\gamma_{\varepsilon}^{-k}),\\
\widetilde\Theta_1&=o(Y_1^{2})+o(\lambda_{\varepsilon}^{k})+\lambda_{\varepsilon}^{k}\cdot O(Y_1),\\
\widetilde\Theta_2&=O(Y_1^2)+O(\lambda_{\varepsilon}^{k}), \quad\quad\;\;\; \widetilde\Theta_3=O(Y_1)+O(\lambda_{\varepsilon}^{k}).
\end{aligned}
\end{equation*}
By Definition \ref{coordinatedefinitioncorank}, the problem of finding a corank-2 tangency between the manifolds above is equivalent to finding a corank-2 solution to the system consisting of \eqref{lemmacorank2tangencyeq1} and \eqref{lemmacorank2tangencyeq2}, therefore, in what follows, we will consider such a system. Applying the Implicit Function Theorem to this system, we express $X_2,U,Y_2,V$ as functions of $X_1,Y_1,\varepsilon$. As a result, we obtain
\begin{equation}\label{lemmacorank2tangencyeq3}
\begin{aligned} 
0&=\eta_0+\alpha+\eta^*_1\cdot X_1+a\cdot X_1^2-\gamma_{\varepsilon}^{-k}\cdot(\cos(k\psi_{\varepsilon})+b\cdot\sin(k\psi_{\varepsilon}))\cdot Y_1+\Theta^*_1(X_1,Y_1,\varepsilon),\\
0&=\tilde\eta_0+\beta-\lambda_{\varepsilon}^{k}\cdot(\cos(k\varphi_{\varepsilon})-\tilde b\cdot \sin(k\varphi_{\varepsilon}))\cdot X_1+ \tilde\eta^*_{1}\cdot Y_1+\tilde a\cdot Y_1^2+\widetilde\Theta^*_1(X_1,Y_1,\varepsilon),\\
X_2&=\Theta^*_2(X_1,Y_1,\varepsilon), \quad Y_2=\widetilde\Theta^*_2(X_1,Y_1,\varepsilon), \quad U=\Theta^*_3(X_1,Y_1,\varepsilon), \quad V=\widetilde\Theta^*_3(X_1,Y_1,\varepsilon),
\end{aligned}
\end{equation}
where
\begin{equation*}
\begin{aligned} 
\alpha&=(-\gamma_{\varepsilon}^{-k})\cdot(\tilde y_1^-\cdot(\cos(k\psi_{\varepsilon})+b\cdot\sin(k\psi_{\varepsilon}))+\tilde y_2^-\cdot(\sin(k\psi_{\varepsilon})-b\cdot\cos(k\psi_{\varepsilon})))+o(\gamma_{\varepsilon}^{-k}),\\
\beta&=\lambda_{\varepsilon}^{k}\cdot(x_1^+\cdot(-\cos(k\varphi_{\varepsilon})+\tilde b\cdot\sin(k\varphi_{\varepsilon}))+x_2^+\cdot(\sin(k\varphi_{\varepsilon})+\tilde b\cdot\cos(k\varphi_{\varepsilon})))+o(\lambda_{\varepsilon}^{k}),\\
\eta^*_1&=\eta_1+O(\gamma_{\varepsilon}^{-k}), \quad \tilde\eta^*_1=\tilde\eta_1+O(\lambda_{\varepsilon}^{k}),
\end{aligned}
\end{equation*}
and the functions $\Theta^*, \widetilde\Theta^*$ satisfy
\begin{equation*}
\begin{aligned} 
\Theta^*_1&=o(X_1^2)+o(\gamma_{\varepsilon}^{-k})+\gamma_{\varepsilon}^{-k}\cdot (O(X_1)+O(Y_1^2)),\\
\Theta^*_2&=O(X_1^2)+O(\gamma_{\varepsilon}^{-k}), \quad\quad \Theta^*_3=O(X_1)+O(\gamma_{\varepsilon}^{-k}),\\
\widetilde\Theta^*_1&=o(Y_1^{2})+o(\lambda_{\varepsilon}^{k})+\lambda_{\varepsilon}^{k}\cdot (O\left(Y_1\right)+O(X_1^2)),\\
\widetilde\Theta^*_2&=O(Y_1^2)+O(\lambda_{\varepsilon}^{k}), \quad\quad\;\;\; \widetilde\Theta^*_3=O(Y_1)+O(\lambda_{\varepsilon}^{k}).
\end{aligned}
\end{equation*}
Since $X_2,U,Y_2,V$ are now functions of $X_1,Y_1$ and $\varepsilon$, it is enough to prove that the system consisting of the first two equations of \eqref{lemmacorank2tangencyeq3} has a corank-2 solution. Namely, we need to find $X_1,Y_1$ and $\varepsilon$ which solve
\begin{equation}\label{lemmacorank2tangencyeq3.5}
\begin{aligned} 
\eta_0&=-\alpha-\eta^*_1\cdot X_1-a\cdot X_1^2+\gamma_{\varepsilon}^{-k}\cdot(\cos(k\psi_{\varepsilon})+b\cdot\sin(k\psi_{\varepsilon}))\cdot Y_1-\Theta^*_1(X_1,Y_1,\varepsilon),\\
\tilde\eta_0&=-\beta+\lambda_{\varepsilon}^{k}\cdot(\cos(k\varphi_{\varepsilon})-\tilde b\cdot \sin(k\varphi_{\varepsilon}))\cdot X_1- \tilde\eta^*_{1}\cdot Y_1-\tilde a\cdot Y_1^2-\widetilde\Theta^*_1(X_1,Y_1,\varepsilon),
\end{aligned}
\end{equation}
and nullify the Jacobian matrix of the right-hand side of this system. This means that we need to solve \eqref{lemmacorank2tangencyeq3.5} augmented by
\begin{equation}\label{lemmacorank2tangencyeq4}
\begin{aligned} 
0&=\gamma_{\varepsilon}^{-k}\cdot(\cos(k\psi_{\varepsilon})+b\cdot\sin(k\psi_{\varepsilon}))+o(\gamma_{\varepsilon}^{-k})+\gamma_{\varepsilon}^{-k}\cdot O(Y_1),\\
0&=\lambda_{\varepsilon}^k\cdot(\cos(k\varphi_{\varepsilon})-\tilde b\cdot\sin(k\varphi_{\varepsilon}))+o(\lambda_{\varepsilon}^k)+\lambda_{\varepsilon}^{k}\cdot O(X_1),\\
0&=2a\cdot X_1+\eta_1^*+o(X_1)+O(\gamma_{\varepsilon}^{-k}),\\ 
0&=2\tilde a\cdot Y_1+\tilde\eta_1^*+o(Y_1)+O(\lambda_{\varepsilon}^k).
\end{aligned}
\end{equation}
Taking into account that $a\not=0$ and $\tilde a\not=0$ (the tangencies $\Gamma,\widetilde\Gamma$ are quadratic), we recast the system \eqref{lemmacorank2tangencyeq4} as
\begin{equation}\label{lemmacorank2tangencyeq5}
\begin{aligned} 
\widehat\psi&=-\psi_0-\frac{\text{arccot}(b)+\pi l_1}{k}+o\left(\frac{1}{k}\right), &\widehat\varphi&=-\varphi_0+\frac{\text{arccot}(\tilde b)+\pi l_2}{k}+o\left(\frac{1}{k}\right),\\
X_1&=-\frac{\eta_1^*}{2a}+O(\gamma_{\varepsilon}^{-k}), &Y_1&=-\frac{\tilde\eta_1^*}{2\tilde a}+O(\lambda_{\varepsilon}^k),
\end{aligned}
\end{equation}
where $l_1,l_2\in\mathbb Z$. Here, the small terms are functions of $X_1,Y_1$ and $\varepsilon$, and all the coefficients are functions of $\varepsilon$. We choose the integers $l_1(k)$ and $l_2(k)$ in \eqref{lemmacorank2tangencyeq5} such that 
\begin{equation*}
\begin{aligned} 
&-\frac{\pi l_1(k)}{k}-\psi_0\to 0 \quad \text{and} \quad \frac{\pi l_2(k)}{k}-\varphi_0\to 0 \quad \text{as} \quad k\to+\infty.
\end{aligned}
\end{equation*}
Then, by the Implicit Function Theorem, the system of \eqref{lemmacorank2tangencyeq3.5} and \eqref{lemmacorank2tangencyeq5} has, for every sufficiently large $k$, a unique solution $\eta_0,\tilde\eta_0,\widehat\varphi,\widehat\psi,X_1, Y_1$ which converges to zero as $k\to+\infty$.
\end{demo}

\begin{remark}
Lemma \ref{lemmacorank2tangency} can be reformulated for the number of parameters larger than four (like Lemma \ref{lemma1highordertangencies}). In this case, the existence of corank-2 tangency would correspond to smooth codimension-four manifolds in the space of parameters.
\end{remark}

\begin{lemma}\label{lemmahighordercorank2tangencies}
Let $f:\mathcal M^d\rightarrow \mathcal M^d$ be a $C^r$ diffeomorphism, where $d\ge 4$ and $r=2,.\ldots,\infty,\omega$. Let $f$ have a bi-focus periodic orbit $L$ and $h=2\cdot2^{\frac{(n-1)(n+4)}{2}}$ orbits $\Gamma_1,\ldots,\Gamma_h$ of corank-1 quadratic homoclinic tangency between $W^s(L)$ and $W^u(L)$, where $n\in\mathbb N, n\le r-1$. Then, $C^r$-close to $f$, there exists a diffeomorphism with an orbit of corank-2 homoclinic tangency of order $n$ between $W^s(L)$ and $W^u(L)$.  
\end{lemma}

\begin{demo}
Using Lemma \ref{lemmasplittingoftangency}, we embed $f$ into a four-parameter $C^r$ family of diffeomorphisms, which splits the orbits $\Gamma_1,\Gamma_2$ freely and independently up to order $0$, changes the arguments $\varphi,\psi$ of $L$ independently as well, and does not split the remaining orbits $\Gamma_3,\ldots,\Gamma_h$. By Lemma \ref{lemmacorank2tangency}, in this family, there exists a diffeomorphism $f_1$ which is $C^r$-close to $f$ and has an orbit of corank-2 homoclinic tangency between $W^s(L)$ and $W^u(L)$. Repeating this procedure $\frac{h}{2}$ times with pairs $\{\Gamma_3,\Gamma_4\},\ldots,\{\Gamma_{h-1},\Gamma_h\}$ (without splitting the corank-2 tangencies obtained at the previous steps), we arrive at the diffeomorphism $\hat f$ which is $C^r$-close to $f$ and has $\frac{h}{2}$ orbits of corank-2 homoclinic tangency between $W^s(L)$ and $W^u(L)$. Employing Theorem 2 of \cite{Mints}, we find a diffeomorphism that is $C^r$-close to $\hat f$ and possesses an orbit of corank-2 homoclinic tangency of order $n$ between $W^s(L)$ and $W^u(L)$.
\end{demo}

\subsection{Wild basic sets}

\begin{lemma}\label{lemmamanytangencies}
Let $f:\mathcal M^d\rightarrow \mathcal M^d$ be a $C^r$ diffeomorphism, where $d\ge 2$ and $r=1,.\ldots,\infty,\omega$. Let $f$ belong to a Newhouse region $\mathcal N(L,\Lambda)$, and let $h\in\mathbb N$. Then, $C^r$-close to $f$, there exists a diffeomorphism with $h$ different orbits of corank-$1$ quadratic homoclinic tangency between $W^s(L)$ and $W^u(L)$.
\end{lemma}

\begin{demo}
The periodic orbit $L$ and the wild basic set $\Lambda$ are homoclinically related, so $W^s(L)$ and $W^u(L)$ accumulate on the stable and unstable manifolds, respectively, of every point of $\Lambda$. This property is fulfilled for all diffeomorphisms in a sufficiently small $C^r$ neighborhood of $f$.

Since $f\in\mathcal N(L,\Lambda)$, there exist points $p_0,q_0\in\Lambda$ such that $W^s(p_0)$ and $W^u(q_0)$ contain an orbit $\Gamma_0$ of tangency. Using Lemma \ref{lemmasplittingoftangency}, we embed $f$ into a finite-parameter $C^{r}$ family of diffeomorphisms, which splits the tangency $\Gamma_0$ freely up to order $2$. By the definition of free splitting, a slight variation of parameters in this family ensures that this tangency is of corank $1$ and quadratic. By Lemma \ref{lemmasplittangencyinbasicset}, in the same family, we find a diffeomorphism $f_1$ that is $C^{r}$-close to $f$ and has an orbit $\Gamma_1$ of corank-$1$ quadratic homoclinic tangency between the stable and unstable manifolds of $L$. Furthermore, because $f_1$ remains within the region $\mathcal N(L,\Lambda)$, there exist points $p_1,q_1\in\Lambda$ such that $W^s(p_1)$ and $W^u(q_1)$ contain an orbit $\Gamma_{00}$ of tangency.

Next, we apply the same reasoning for $f_1$. Specifically, employing Lemma \ref{lemmasplittingoftangency}, we embed $f_1$ into a finite-parameter $C^{r}$ family of diffeomorphisms that splits the tangency $\Gamma_{00}$ freely up to order $2$, and does not split the tangency $\Gamma_1$. By slightly changing parameters in this family, we make the tangency $\Gamma_{00}$ be of corank $1$ and quadratic. By Lemma \ref{lemmasplittangencyinbasicset}, in this family, we find a diffeomorphism $f_2$ that is $C^{r}$-close to $f_1$ and has orbits $\Gamma_1,\Gamma_2$ of corank-$1$ quadratic homoclinic tangency between the stable and unstable manifolds of $L$. Additionally, there exist points $p_2,q_2\in\Lambda$ such that $W^s(p_2)$ and $W^u(q_2)$ contain an orbit $\Gamma_{000}$ of tangency.

By applying this algorithm $h$ times, we obtain a $C^{r}$ diffeomorphism $f_h$ with $h$ different orbits $\Gamma_1,\ldots\Gamma_h$ of corank-$1$ quadratic homoclinic tangency between $W^s(L)$ and $W^u(L)$. By construction, $f_h$ is $C^r$-close to $f$, which completes the proof.
\end{demo}

\begin{remark}
This result can be proven applying the original approach \cite{Gonchenkoonmodelswith},\cite{arbitraryorderconservative} which does not refer to the Newhouse Theorem and uses the existence of the so-called $\Omega$-moduli - continuous invariants of $\Omega$-equivalence \cite{GonchenkoShilnikov1987},\cite{Gonchenkoonmoduli1993},\cite{dynamicalphenomena}. The proof we employ here follows Lemma 2.10 in \cite{Mints}, and is similar in spirit to Propositions 6.1,6.9 in \cite{Berger1}. 
\end{remark}

\begin{demo4}
Let $f$ be a diffeomorphism from a Newhouse region $\mathcal N(L,\Lambda)$. We assume $f$ is real-analytic (if $f$ is only $C^r$-smooth, then one can always find a real-analytic diffeomorphism $f^{new}\in\mathcal N(L,\Lambda)$ that is $C^r$-close to $f$). The non-trivial basic set $\Lambda$ is wild, so there exists a compact set $\mathcal L\subset W^s(\Lambda)\cup W^u(\Lambda)$ with the following property: given any compact neighborhood $U$ of $\mathcal L$ there exists a neighborhood $\mathcal U$ of $f$ such that for any diffeomorphism $g\in\mathcal U$ the manifolds $W^s(\Lambda)$ and $W^u(\Lambda)$ contain an orbit of tangency lying entirely in $U$. We choose $U$ such that it does not intersect the orbit $L$, thereby fixing the neighborhood $\mathcal U$ of $f$ for the remainder of the proof. 

Let $\Omega$ be a basic set containing both $L$ and $\Lambda$. We define the following sequences:
\begin{itemize}
    \item $\langle\mathcal O_{i,1};\mathcal O_{i,2}\rangle_{i\in\mathbb N}$: a sequence of pairs of periodic orbits from $\Omega$; 
    
    \item $\langle c_i\rangle_{i\in\mathbb N}$: a sequence of integers, where $c_i=1$ if $L$ is not bi-focus, and $c_i\in\{1,2\}$ if $L$ is bi-focus, representing coranks;
    
    \item $\langle n_i\rangle_{i\in\mathbb N}$: a sequence of natural numbers representing orders.  
\end{itemize}
These sequences are constructed such that every possible combination comprising any pair of periodic orbits from $\Omega$, any admissible number $c_i$, and any natural number $n_i$ appears for infinitely many indices $i$. Our goal is to prove the existence of a $C^{\omega}$-small perturbation of $f$ within its neighborhood $\mathcal U$ that produces a diffeomorphism possessing, for each $i\in\mathbb N$, an orbit of corank-$c_i$ tangency of order $n_i$ between $W^s(\mathcal O_{i,1})$ and $W^u(\mathcal O_{i,2})$.

Let $\delta>0$ be arbitrarily small. Take a sequence of positive values $\langle\delta_i\rangle_{i\in\mathbb N}$ such that $\sum\limits_{i=1}^{\infty}\delta_i=\delta$. We construct a sequence of diffeomorphisms $\langle f_i\rangle_{i\in\mathbb N}$, starting with $f_0\equiv f$, such that: 
\begin{itemize}

\item $f_i\in\mathcal U$, ensuring the existence of an orbit $\Gamma'_i$ of tangency between $W^s(\Lambda)$ and $W^u(\Lambda)$ that lies entirely in $U$;

\item $f_i$ keeps the $i$ orbits $\{\Gamma_1,\ldots,\Gamma_i\}$ of tangency generated in the previous steps, with their prescribed coranks and orders; 

\item the distance between $f_i$ and $f_{i-1}$ is less than $\delta_i$. 
\end{itemize}
The limit $f^*$ of this sequence is at a distance at most $\delta$ from $f$, and possesses the required infinite sequence of orbits of tangency.

Thus, we have to show that given $f_i$, one can produce $f_{i+1}$ that preserves corank and order of the existing $i$ orbits $\{\Gamma_1,\ldots,\Gamma_i\}$ of tangency and also has a new orbit $\Gamma_{i+1}$ of corank-$c_{i+1}$ tangency of order $n_{i+1}$ between $W^s(\mathcal O_{i+1,1})$ and $W^u(\mathcal O_{i+1,2})$. This is achieved via a finite sequence of $C^{\omega}$-small perturbations with a total size less than $\delta_{i+1}$, in the following way. 

By Lemma \ref{lemmamanytangencies}, $C^{\omega}$-close to $f_i$, there exists a diffeomorphism $\hat f_i$ with an orbit of tangency between $W^s(\Lambda)$ and $W^u(\Lambda)$, and $h$ orbits of corank-1 quadratic homoclinic tangency between $W^s(L)$ and $W^u(L)$, where $h$ can be chosen to be any given integer. 

\underline{\textbf{Case 1: $c_{i+1}=1$.}} We set $h=n_{i+1}$. 
Employing Lemma \ref{lemma2highordertangencies}, we obtain a diffeomorphism $\overline f_i$ with an orbit $\overline\Gamma$ of corank-1 homoclinic tangency of order $n_{i+1}$ between $W^s(L)$ and $W^u(L)$. 

\underline{\textbf{Case 2: $c_{i+1}=2$.}} We set $h=2\cdot 2^{\frac{\left(n_{i+1}-1\right)\left(n_{i+1}+4\right)}{2}}$. Applying Lemma \ref{lemmahighordercorank2tangencies}, we obtain a diffeomorphism $\overline f_i$ with an orbit $\overline\Gamma$ of corank-2 homoclinic tangency of order $n_{i+1}$ between $W^s(L)$ and $W^u(L)$. 

Since the orbits $L, \mathcal O_{i+1,1}$ and $\mathcal O_{i+1,2}$ belong to the same non-trivial basic set $\Omega$, the manifolds $W^s(\mathcal O_{i+1,1})$ and $W^u(\mathcal O_{i+1,2})$ accumulate on $W^s(L)$ and $W^u(L)$, respectively. Making use of Lemma \ref{lemmasplittingoftangency}, we embed $\overline f_i$ into a finite-parameter analytic family of diffeomorphisms in which the tangency $\overline\Gamma$ splits freely up to order $n_{i+1}$. By Lemma \ref{lemmasplittangencyinbasicset}, in this family, there exists a diffeomorphism $f_{i+1}$ that is $C^{\omega}$-close to $\overline f_i$ and has the sought orbit $\Gamma_{i+1}$ of corank-$c_{i+1}$ tangency of order $n_{i+1}$ between $W^s(\mathcal O_{i+1,1})$ and $W^u(\mathcal O_{i+1,2})$. By construction, at least one point of the orbit $\Gamma_{i+1}$ is close to $L$ (proximity to $L$ is controlled by the size of the perturbation that splits the tangency $\overline\Gamma$). Consequently, $\Gamma_{i+1}$ does not lie entirely in $U$. Since $f_{i+1}\in\mathcal U$, it also possesses an orbit $\Gamma'_{i+1}$ of tangency between $W^s(\Lambda)$ and $W^u(\Lambda)$ lying entirely in $U$. Thus, the orbits $\Gamma_{i+1}$ and $\Gamma'_{i+1}$ are distinct.

The construction of $f_{i+1}$ given above preserves the $i$ orbits $\{\Gamma_1,\ldots,\Gamma_i\}$ of tangency, along with their coranks and orders, inherited from the diffeomorphism $f_i$, thereby completing the inductive step and the proof of the theorem.
\end{demo4}

\section{Paratangencies of all orders}

In this section, we prove Theorem \ref{theorem5}.

\subsection{Paratangency functionals}\label{subsectionparatangencyfunctionals}
Let $(f_{\zeta})$ be an $\ell$-parameter $C^r$ family of diffeomorphisms of the manifold $\mathcal M^d$, where $\ell\ge 1$, $d\ge 2$, and $r=\infty,\omega$. Let $(f_{\zeta})$ have a paratangency of order $m$ $(m\in\mathbb N)$ with the bifurcation set $\mathcal X$. We also consider the case $m=0$ and, for convenience, say that a parametric family $(f_{\zeta})$ has a paratangency of order $0$ with $\mathcal X$ if there exists a smooth manifold $\mathcal Y\subseteq \mathcal X$ such that the graph $(\zeta, f_{\zeta})\subset D^\ell\times\operatorname{Diff}^r(\mathcal M^d)$ is transverse to $D^\ell\times\mathcal Y$ at the point of intersection (cf. Definition \ref{definitionparatangency}). 

Let $\eta: \mathrm{Diff}^r(\mathcal M^d)\to\mathbb R$ be a smooth functional such that the set $\{\eta=0\}$ locally contains the manifold $\mathcal Y\subseteq \mathcal X$ (see Definition \ref{definitionparatangency}). We refer to the manifold $\mathcal Y$ as the \textit{bifurcation stratum} and to the functional $\eta$ as a \textit{paratangency functional}. We interpret  Definition \ref{definitionparatangency} as follows: for every paratangency functional $\eta$, one has
\begin{equation}\label{paratangencyordereq0}
\begin{aligned} 
&\left.\frac{\partial^{j}}{\partial \zeta^{j}}\,\eta(f_{\zeta})\right|_{\zeta=0}=0 \;\;\; \text{for} \;\;\; 0\le |j|\le m, \;\; \text{and} \;\;\; \left.\frac{\partial^{j}}{\partial \zeta^{j}}\,\eta(f_{\zeta})\right|_{\zeta=0}\not=0 \;\; \text{for some} \;\; j \;\; \text{with} \;\; |j|=m+1.
\end{aligned}
\end{equation}
Moreover, we will consider only bifurcation strata $\mathcal Y$ of finite codimension. Namely, we always assume the existence of $C^{\infty}$ paratangency functionals $\eta_1,\ldots,\eta_s$ that define $\mathcal Y$, i.e., $\mathcal Y=\{\eta_1=0,\ldots,\eta_s=0\}$.

Let us fix a system of coordinates in the ball $D^\ell$ of the parameters $\zeta$. A paratangency functional $\eta$ can be expanded as 
\begin{equation*}
\begin{aligned} 
&\eta(f_{\zeta})=\sum\limits_{|j|=m+1}c_{j}\cdot \zeta^j+O(||\zeta||^{m+2}),
\end{aligned}
\end{equation*}
where each term $c_j\cdot \zeta^j=c_j\cdot\zeta_1^{j_1}\cdot\ldots\cdot\zeta_\ell^{j_\ell}$ is a monomial of degree $m+1$. 

We define an order $\prec$ on the set of multi-indices $j=(j_1,\ldots,j_\ell)$ as follows: 
\begin{definition}\label{definitionmultiindexorder}
We say that $j\prec j'$ \, if \, $|j|<|j'|$, or $|j|=|j'|$ and there exists $k\in\{1,\ldots,\ell\}$ such that $j_i=j'_i$ for all $0<i<k$, and $j_k<j'_k$.
\end{definition}
Note that the number of multi-indices of degree $m+1$ with $\ell$ entries equals to $\binom{\ell+m}{m+1}$. For $t\in\{1,\ldots,\binom{\ell+m}{m+1}\}$, we denote by $j^{(m+1,t)}$ the $t$-th multi-index in the set of multi-indices of degree $m+1$ that are lexicographically ordered according to Definition \ref{definitionmultiindexorder}.

\begin{definition}\label{definitionindexofparatangencyfunctional}
We say, for $t\in\{1,\ldots,\binom{\ell+m}{m+1}\}$, that the paratangency functional $\eta$ is of index $(m,t)$ (relative to $f_{\zeta}$) if $c_{j}=0$ for all $j\prec j^{(m,t)}$ and $c_{j^{(m,t)}}\not=0$.
\end{definition}

We introduce the natural ordering for pairs:
$$(m,t)<(m',t') \Longleftrightarrow\left[\begin{aligned} m&<m', \\ m&=m' \;\; \text{and} \;\; t<t'.\end{aligned}\right.$$
For convenience of notation, we adopt the following convention: 
$$(m,\binom{\ell+m}{m+1}+1)\equiv (m+1,1).$$

\begin{definition}\label{definitionindexofparatangency}
We say that the paratangency is of index $(m,t)$ if there exist paratangency functionals of index $(m,t)$ that locally define the bifurcation stratum $\mathcal Y$. 
\end{definition}

\subsection{Perturbations of paratangency functionals}
Assume that the family $(f_{\zeta})$ is embedded into a finite-parameter $C^r$ family $(f_{\zeta,\varepsilon})$. Fix positive integers $m'$ and $t'$. The paratangency functional $\eta$ can be expanded as
\begin{equation*}
\begin{aligned} 
&\eta(f_{\zeta,\varepsilon})=\eta(f_{\zeta,0})+\sum\limits_{j\prec j^{(m',t')}}\xi_{j}\cdot \zeta^j+O_{(m',t')}(\zeta),
\end{aligned}
\end{equation*}
where $\xi_{j}$ are $C^r$ functions of $\varepsilon$ vanishing at $\varepsilon=0$, and $O_{(m',t')}$ denotes a function whose Taylor expansion starts with the term $\zeta^{j^{(m',t')}}$. We call $\xi_j$ the \textit{parasplitting coefficients} and denote by $\mathbf{\Xi}_{(m',t')}$ the vector of all such coefficients with $j\prec j^{(m',t')}$.
\begin{definition}
The parametric family $(f_{\zeta,\varepsilon})$ perturbs the paratangency functional $\eta$ freely up to index $(m',t')$ if
\begin{equation*} 
\left.\text{rank}\;\frac{\partial\mathbf \Xi_{(m',t')}}{\partial\varepsilon}\right |_{\varepsilon=0}=\dim\mathbf \Xi_{(m',t')}. 
\end{equation*}
\end{definition}

Now let $\eta_1,\ldots,\eta_{s}: \mathrm{Diff}^r(\mathcal M^d)\to\mathbb R$ be paratangency functionals corresponding to the bifurcation stratum $\mathcal Y\subseteq \mathcal X$. Fix positive integers $m'_1,\ldots,m'_s$ and $t'_1,\ldots,t'_s$. Expanding the functional $\eta_k$, we obtain
\begin{equation*}
\begin{aligned} 
&\eta_k(f_{\zeta,\varepsilon})=\eta_k(f_{\zeta,0})+\sum\limits_{j\prec j^{(m'_k,t'_k)}}\xi_{k,j}\cdot \zeta^j+O_{(m'_k,t'_k)}(\zeta),
\end{aligned}
\end{equation*}
where $\xi_{k,j}$ are parasplitting coefficients. We denote by $\mathbf{\Xi}^k_{(m'_k,t'_k)}$ the vector of all such coefficients with $j\prec j^{(m'_k,t'_k)}$.

\begin{definition}
The parametric family $(f_{\zeta,\varepsilon})$ perturbs the paratangency functionals $\eta_1,\ldots,\eta_s$ freely up to indices $(m'_1,t'_1),\ldots,(m'_s,t'_s)$ if
\begin{equation*}
\left.\text{rank}\;\frac{\partial\Big(\mathbf \Xi^{1}_{(m'_1,t'_1)},\ldots,\mathbf \Xi^{s}_{(m'_s,t'_s)}\Big)}{\partial\varepsilon}\right|_{\varepsilon=0}=\sum\limits_{k=1}^s\dim\mathbf \Xi^{k}_{(m'_k,t'_k)}.
\end{equation*} 
\end{definition}

We focus on a particular type of bifurcation strata $\mathcal Y=\mathcal Y(\Gamma_1,\ldots,\Gamma_h)$ corresponding to the existence of a finite number of homoclinic tangencies $\Gamma_1,\ldots,\Gamma_h$ of finite orders. In this case, each functional $\eta$ is a splitting coefficient of a homoclinic tangency. Note that in the definition of the splitting coefficients, we considered $\eta$ as a function of a finite-dimensional parameter $\varepsilon$ but the definition also holds when $\varepsilon$ runs a ball in the infinite-dimensional space of $C^r$ diffeomorphisms $f$. 

Similar to Lemma \ref{lemmasplittingoftangency}, we obtain
\begin{lemma}\label{lemmaperturbationofparatangencyfunctionals}
Let $(f_{\zeta})$ be an $\ell$-parameter $C^r$ family of diffeomorphisms of $\mathcal M^d$, where $\ell\ge 1$, $d\ge 2$ and $r=\infty,\omega$. Let $\eta_1,\ldots,\eta_{s}$ be paratangency functionals, each representing some splitting coefficient of one of a finite collection of homoclinic tangencies. Then, for any non-negative integers $m'_1,\ldots,m'_{s}$ and positive integers $t'_1,\ldots,t'_{s}$, there exists a finite-parameter analytic family $(\mathcal F_{\zeta,\varepsilon})$ such that $\mathcal F_{0,0}=id$, and the family $(\mathcal F_{\zeta,\varepsilon}\circ f_{\zeta})$ perturbs $\eta_1,\ldots,\eta_{s}$ freely up to indices $(m'_1,t'_1),\ldots,(m'_{s},t'_{s})$.
\end{lemma}

\subsection{Bifurcation sets of homoclinic tangencies}

Let $f:\mathcal M^d\rightarrow \mathcal M^d$ be a $C^r$ diffeomorphism, where $d\ge 2$ and $r=\infty,\omega$. Let $f$ have a saddle periodic orbit $L$. We set $\lambda\equiv|\lambda_1|$ and $\gamma\equiv|\gamma_1|$, where $\lambda_1$ and $\gamma_1$ are leading stable and leading unstable multipliers of $L$, respectively. Below, we assume the following condition:

\vspace{2mm}

\textbf{\hypertarget{C4.}{C4.}}  $\lambda\gamma<1$.

\vspace{2mm}

This is not a restriction of generality since the expanding case $\lambda\gamma>1$ is treated by considering the inverse map $f^{-1}$. 

Let us fix a coordinate system in the ball $D^\ell$ (where the parameters $\zeta$ are defined), which we will keep unchanged throughout the whole section. In what follows, we assume that the saddle periodic orbit $L$ satisfies conditions \textbf{\hyperlink{C1.}{C1}},\textbf{\hyperlink{C4.}{C4}}, and the orbits of homoclinic tangency satisfy conditions \textbf{\hyperlink{C2.}{C2}},\textbf{\hyperlink{C3.}{C3}}. 

\begin{lemma}\label{lemmaincreasingparatangencyindex}
Let $f_0:\mathcal M^d\rightarrow \mathcal M^d$ be a $C^r$ diffeomorphism, where $d\ge 2$ and $r=\infty,\omega$. Let $f_0$ have a saddle periodic orbit $L$ and two orbits $\Gamma,\widetilde\Gamma$ of corank-1 quadratic homoclinic tangency between $W^s(L)$ and $W^u(L)$. Let $f_0$ be embedded into an $\ell$-parameter $C^r$ family $(f_{\zeta})$ which has a paratangency of index $(m,t)$ with a bifurcation stratum $\mathcal Y(\Gamma,\widetilde\Gamma)$, where $\ell\ge 1,m\ge0,t\ge1$. Then, $C^{r}$-close to $(f_{\zeta})$, there exists a parametric family $(g_{\zeta})$ which, at $\zeta=0$, has a paratangency of index at least $(m,t+1)$ with the bifurcation stratum corresponding to a corank-1 quadratic homoclinic tangency between $W^s(L)$ and $W^u(L)$. 
\end{lemma}

\begin{demo}
Let us choose a point $O\in L$ and define the local map $T_O$ in a small neighborhood $U(O)$. Applying Lemma \ref{lemmachoiceofcoordinates}, we introduce $(T_{O},\Gamma,\widetilde\Gamma)$-adapted coordinates in $U(O)$. In these coordinates, for all small $\zeta$, the manifolds $\mathcal W^s(O)$ and $\widetilde{\mathcal W}^u(O)$ are straightened, whereas the manifolds $\mathcal W^u(O)$ and $\widetilde{\mathcal W}^s(O)$ have the form 
\begin{equation*}
\mathcal W^u(O): \;\;
\begin{aligned}
y_1&=\eta'_0+\eta'_1\cdot(x_1-x_1^+)+a\cdot(x_1-x_1^+)^{2}+b\cdot y_2+\\
&+o((x_1-x_1^+)^{2})+\theta'_1(x_1-x_1^+,y_2,v,\zeta),\\
x_2- x_2^+&=\theta'_2(x_1-x_1^+,y_2,v,\zeta),\\ 
u-u^+&=\theta'_3(x_1-x_1^+,y_2,v,\zeta),
\end{aligned}
\end{equation*}
and 
\begin{equation*}
\widetilde{\mathcal W}^s(O): \;\;
\begin{aligned} 
\tilde x_1&=\tilde\eta'_0+\tilde\eta'_1\cdot(\tilde y_1-\tilde y_1^-)+\tilde a\cdot(\tilde y_1-\tilde y_1^-)^{2}+\tilde b\cdot x_2+\\
&+o((\tilde y_1-\tilde y_1^-)^{2})+\tilde\theta'_1(\tilde x_2,\tilde u,\tilde y_1-\tilde y_1^-,\zeta),\\
\tilde y_2-\tilde y_2^-&=\tilde\theta'_2(\tilde x_2,\tilde u,\tilde y_1-\tilde y_1^-,\zeta),\\
\tilde v-\tilde v^-&=\tilde\theta'_3(\tilde x_2,\tilde u,\tilde y_1-\tilde y_1^-,\zeta),
\end{aligned}
\end{equation*}
where the functions $\theta',\tilde\theta'$ satisfy \eqref{lemmachoiceofcoordinatesestimateseq2}; $a\not=0,\tilde a\not =0,b,\tilde b$ are $C^r$ in $\zeta$. The paratangency is of index $(m,t)$, therefore the adapted coordinates can be chosen such that the splitting coefficients $\eta',\tilde\eta'$ are expressed by
\begin{equation*}
\begin{aligned}
\eta'_0&=c_{0}\cdot\zeta^{j^{(m,t)}}+O_{(m,t+1)}(\zeta), &\tilde\eta'_0&=\tilde c_{0}\cdot\zeta^{j^{(m,t)}}+O_{(m,t+1)}(\zeta),\\
\eta'_1&=c_{1}\cdot\zeta^{j^{(m,t)}}+O_{(m,t+1)}(\zeta), &\tilde\eta'_1&=\tilde c_{1}\cdot\zeta^{j^{(m,t)}}+O_{(m,t+1)}(\zeta),
\end{aligned}
\end{equation*}
with $c_{0}\not=0, c_{1}\not=0$ and $\tilde c_{0}\not=0,\tilde c_{1}\not=0$. Making use of Lemma \ref{lemmaperturbationofparatangencyfunctionals}, we find a $C^r$ family $(f_{\zeta,\varepsilon,\tilde\varepsilon})$ such that $(f_{\zeta,0,0})=(f_{\zeta})$, the manifolds $\mathcal W^s(O)$ and $\widetilde{\mathcal W}^u(O)$ remain straightened, whereas the manifolds $\mathcal W^u(O)$ and $\widetilde{\mathcal W}^s(O)$ take the form
\begin{equation}\label{lemmaincreasingparatangencyindexeq4}
\mathcal W^u(O): \;\;
\begin{aligned}
y_1&=\eta_0+\eta_1\cdot(x_1-x_1^+)+a\cdot(x_1-x_1^+)^{2}+b\cdot y_2+\\
&+o((x_1-x_1^+)^{2})+\theta_1(x_1-x_1^+,y_2,v,\zeta,\varepsilon,\tilde\varepsilon),\\
x_2- x_2^+&=\theta_2(x_1-x_1^+,y_2,v,\zeta,\varepsilon,\tilde\varepsilon),\\
u-u^+&=\theta_3(x_1-x_1^+,y_2,v,\zeta,\varepsilon,\tilde\varepsilon),
\end{aligned}
\end{equation}
and 
\begin{equation}\label{lemmaincreasingparatangencyindexeq5}
\widetilde{\mathcal W}^s(O): \;\;
\begin{aligned} 
\tilde x_1&=\tilde\eta_0+\tilde\eta_1\cdot(\tilde y_1-\tilde y_1^-)+\tilde a\cdot(\tilde y_1-\tilde y_1^-)^{2}+\tilde b\cdot x_2+\\
&+o((\tilde y_1-\tilde y_1^-)^{2})+\tilde\theta_1(\tilde x_2,\tilde u,\tilde y_1-\tilde y_1^-,\zeta,\varepsilon,\tilde\varepsilon),\\
\tilde y_2-\tilde y_2^-&=\tilde\theta_2(\tilde x_2,\tilde u,\tilde y_1-\tilde y_1^-,\zeta,\varepsilon,\tilde\varepsilon),\\
\tilde v-\tilde v^-&=\tilde\theta_3(\tilde x_2,\tilde u,\tilde y_1-\tilde y_1^-,\zeta,\varepsilon,\tilde\varepsilon),
\end{aligned}
\end{equation}
where the functions $\theta,\tilde\theta$ satisfy \eqref{lemmachoiceofcoordinatesestimateseq2}; the splitting coefficients $\eta,\tilde\eta$ are specified by
\begin{equation}\label{lemmaincreasingparatangencyindexeq6}
\begin{aligned}
\eta_0&=\sum\limits_{j\prec j^{(m,t)}}\varepsilon_{0,j}\cdot\zeta^j+c_{0}\cdot\zeta^{j^{(m,t)}}+O_{(m,t+1)}(\zeta),\\
\eta_1&=\sum\limits_{j\prec j^{(m,t)}}\varepsilon_{1,j}\cdot\zeta^j+c_{1}\cdot\zeta^{j^{(m,t)}}+O_{(m,t+1)}(\zeta),
\end{aligned}
\end{equation} 
and
\begin{equation}\label{lemmaincreasingparatangencyindexeq7}
\begin{aligned}
\tilde\eta_0&=\sum\limits_{j\prec j^{(m,t)}}\tilde\varepsilon_{0,j}\cdot\zeta^j+\tilde c_{0}\cdot\zeta^{j^{(m,t)}}+O_{(m,t+1)}(\zeta),\\
\tilde\eta_1&=\tilde\varepsilon_1+\tilde c_{1}\cdot\zeta^{j^{(m,t)}}+O_{(m,t+1)}(\zeta).
\end{aligned}
\end{equation}

By Corollary \ref{corollarylocalmap}, there exists a sequence $\langle k_\alpha\rangle_{\alpha\in\mathbb N}$ such that $k_\alpha\xrightarrow{}+\infty$ and the equations \eqref{eqlocal3tamesequence} hold for any $\langle k_{\alpha}\rangle$-tame sequence of open balls $\langle\mathfrak D_{\alpha}\rangle_{\alpha\in\mathbb N}$ centered at zero and for all $(\zeta,\varepsilon,\tilde\varepsilon)\in\mathfrak D_{\alpha}$. We will show below that there exists a $\langle k_{\alpha}\rangle$-tame sequence $\langle\mathfrak D_{\alpha}\rangle_{\alpha\in\mathbb N}$ and $(\zeta,\varepsilon_{\alpha},\tilde\varepsilon_{\alpha})\in\mathfrak D_{\alpha}$ (so, $\varepsilon_{\alpha},\tilde\varepsilon_{\alpha}\xrightarrow{}0$) such that
\begin{enumerate}
\item for all sufficiently large $\alpha$, the manifolds $\mathcal W^u(O)$ and $T_{O}^{-k_\alpha}(\widetilde{\mathcal W}^s(O))$ have a corank-1 quadratic tangency at $(\zeta,\varepsilon,\tilde\varepsilon)=(0,\varepsilon_{\alpha},\tilde\varepsilon_{\alpha})$.

\item the family $(g_{\zeta})\equiv(f_{\zeta,\varepsilon_{\alpha},\tilde\varepsilon_{\alpha}})$, at $\zeta=0$, has a paratangency of index at least $(m,t+1)$ with the bifurcation stratum corresponding to this homoclinic tangency. 
\end{enumerate}
This will prove the lemma.
 
Substituting formulas \eqref{eqlocal3tamesequence} into  \eqref{lemmaincreasingparatangencyindexeq4} and \eqref{lemmaincreasingparatangencyindexeq5}, we rewrite these equations as
\begin{equation}\label{lemmaincreasingparatangencyindexeq8}
\begin{aligned} 
&\gamma^{-k_{\alpha}}\cdot (Y_1+\tilde y_1^-)=\eta_0+\eta_1\cdot X_1+a\cdot X_1^{2}+\gamma^{-k_{\alpha}}\cdot b\cdot (Y_2+\tilde y_2^-)+\Theta_1(X_1,Y,V,\zeta,\varepsilon,\tilde\varepsilon),\\
&X_2=\Theta_2(X_1,Y,V,\zeta,\varepsilon,\tilde\varepsilon), \quad U=\Theta_3(X_1,Y,V,\zeta,\varepsilon,\tilde\varepsilon),
\end{aligned}
\end{equation}
and
\begin{equation}\label{lemmaincreasingparatangencyindexeq9}
\begin{aligned} 
&\lambda^{k_{\alpha}}\cdot (X_1+x_1^+)=\tilde\eta_0+\tilde\eta_1\cdot Y_1+\tilde a\cdot Y_1^{2}+\lambda^{k_{\alpha}}\cdot \tilde b\cdot(X_2+x_2^+)+\widetilde\Theta_1(X,U,Y_1,\zeta,\varepsilon,\tilde\varepsilon),\\
&Y_2=\widetilde\Theta_2(X,U,Y_1,\zeta,\varepsilon,\tilde\varepsilon), \quad V=\widetilde\Theta_3(X,U,Y_1,\zeta,\varepsilon,\tilde\varepsilon),
\end{aligned}
\end{equation}
where $X,U,Y,V$ are given by \eqref{lemma1highordertangencieseq2.5} and $\Theta, \widetilde\Theta$ are some $C^r$ functions satisfying
\begin{equation*}
\begin{aligned} 
\Theta_1&=o(X_1^{2})+o(\gamma^{-k_{\alpha}})\cdot(Y,V)+\gamma^{-k_{\alpha}} \cdot O(X_1\cdot Y),\\
\Theta_2&=O(X_1^2)+O(\gamma^{-k_{\alpha}}), \quad\quad \Theta_3=O(X_1)+O(\gamma^{-k_{\alpha}}),\\
\widetilde\Theta_1&=o(Y_1^{2})+o(\lambda^{k_{\alpha}})\cdot(X,U)+\lambda^{k_{\alpha}}\cdot O(Y_1\cdot X),\\
\widetilde\Theta_2&=O(Y_1^2)+O(\lambda^{k_{\alpha}}), \quad\quad\;\;\; \widetilde\Theta_3=O(Y_1)+O(\lambda^{k_{\alpha}}).
\end{aligned}
\end{equation*}
We make a shift: 
\begin{equation*}
\begin{aligned} 
(X_2,U)&=(X_2,U)^{new}+O(\gamma^{-k_{\alpha}}), &&(Y_2,V)=(Y_2,V)^{new}+O(\lambda^{k_{\alpha}}),
\end{aligned}
\end{equation*}
to nullify the constant terms in the second and the third equations of both \eqref{lemmaincreasingparatangencyindexeq8} and \eqref{lemmaincreasingparatangencyindexeq9}, so the new functions $\Theta_{2,3},\widetilde\Theta_{2,3}$ satisfy \eqref{lemmaincreasingparatangencyindexeq18}. Also, we offset the contributions to the splitting coefficients coming from functions $\Theta_1$,$\widetilde\Theta_1$ by making a change of parameters $\varepsilon$ and $\tilde\varepsilon$:
\begin{equation}\label{lemmaincreasingparatangencyindexeq11}
\begin{aligned} 
\varepsilon_{0,j}&=\varepsilon_{0,j}^*+O(\gamma^{-k_{\alpha}}), &\varepsilon_{1,j}&=\varepsilon_{1,j}^*+O(\gamma^{-k_{\alpha}}), &&\tilde\varepsilon_{0,j}=\tilde\varepsilon_{0,j}^*+O(\lambda^{k_{\alpha}}) &&\text{for} && j\prec j^{(m,t)},\\
\tilde\varepsilon_{1}&=\tilde\varepsilon_{1}^*+O(\lambda^{k_{\alpha}}).
\end{aligned}
\end{equation}
As a result, we recast \eqref{lemmaincreasingparatangencyindexeq8} and \eqref{lemmaincreasingparatangencyindexeq9} as
\begin{equation}\label{lemmaincreasingparatangencyindexeq14}
\begin{aligned} 
&\gamma^{-k_{\alpha}}\cdot (Y_1-b\cdot Y_2)=\eta^*_0+\eta^*_1\cdot X_1+a\cdot X_1^{2}+\Theta^*_1(X_1,Y,V,\zeta,\varepsilon^*,\tilde\varepsilon^*),\\
&X_2=\Theta^*_2(X_1,Y,V,\zeta,\varepsilon^*,\tilde\varepsilon^*), \quad U=\Theta^*_3(X_1,Y,V,\zeta,\varepsilon,^*\tilde\varepsilon^*),
\end{aligned}
\end{equation}
and
\begin{equation}\label{lemmaincreasingparatangencyindexeq15}
\begin{aligned} 
&\lambda^{k_{\alpha}}\cdot (X_1-\tilde b\cdot X_2)=\tilde\eta^*_0+\tilde\eta^*_1\cdot Y_1+\tilde a\cdot Y_1^{2}+\widetilde\Theta^*_1(X,U,Y_1,\zeta,\varepsilon^*,\tilde\varepsilon^*),\\
&Y_2=\widetilde\Theta^*_2(X,U,Y_1,\zeta,\varepsilon^*,\tilde\varepsilon^*), \quad V=\widetilde\Theta^*_3(X,U,Y_1,\zeta,\varepsilon^*,\tilde\varepsilon^*),
\end{aligned}
\end{equation}
where the splitting coefficients $\eta^*,\tilde\eta^*$ are expressed by 
\begin{equation*}
\begin{aligned}
\eta_0^*&=\sum\limits_{j\prec j^{(m,t)}}\varepsilon^*_{0,j}\cdot\zeta^j+(c_{0}+O(\gamma^{-k_{\alpha}}))\cdot\zeta^{j^{(m,t)}}+O_{(m,t+1)}(\zeta),\\
\eta_1^*&=\sum\limits_{j\prec j^{(m,t)}}\varepsilon^*_{1,j}\cdot\zeta^j+(c_{1}+O(\gamma^{-k_{\alpha}}))\cdot\zeta^{j^{(m,t)}}+O_{(m,t+1)}(\zeta),
\end{aligned}
\end{equation*}
and
\begin{equation*}
\begin{aligned}
\tilde\eta^*_0&=\sum\limits_{j\prec j^{(m,t)}}\tilde\varepsilon^*_{0,j}\cdot\zeta^j+(\tilde c_{0}+O(\lambda^{k_{\alpha}}))\cdot\zeta^{j^{(m,t)}}+O_{(m,t+1)}(\zeta),\\
\tilde\eta^*_1&=\tilde\varepsilon^*_1+\lambda^{k_{\alpha}}\cdot O(\zeta)+(\tilde c_{1}+O(\lambda^{k_{\alpha}}))\cdot\zeta^{j^{(m,t)}}+O_{(m,t+1)}(\zeta),
\end{aligned}
\end{equation*}
and $\Theta^*,\widetilde\Theta^*$ satisfy
\begin{equation}\label{lemmaincreasingparatangencyindexeq18}
\begin{aligned} 
\Theta^*_1&=o(X_1^{2})+o(\gamma^{-k_{\alpha}})\cdot(Y,V)+\gamma^{-k_{\alpha}} \cdot O(X_1\cdot Y),\\
\Theta^*_2&=O(X_1^2)+\gamma^{-k_{\alpha}}\cdot O(X_1,Y,V), \quad\quad \Theta^*_3=O(X_1)+\gamma^{-k_{\alpha}}\cdot O(X_1,Y,V),\\
\widetilde\Theta^*_1&=o(Y_1^{2})+o(\lambda^{k_{\alpha}})\cdot(X,U)+\lambda^{k_{\alpha}}\cdot O(Y_1\cdot X),\\
\widetilde\Theta^*_2&=O(Y_1^2)+\lambda^{k_{\alpha}}\cdot O(X,U,Y_1), \quad\quad\;\;\; \widetilde\Theta^*_3=O(Y_1)+\lambda^{k_{\alpha}}\cdot O(X,U,Y_1).
\end{aligned}
\end{equation}
We put 
\begin{equation}\label{lemmaincreasingparatangencyindexeq18.5}
\begin{aligned} 
\tilde\varepsilon^*_{0,j}&=0 &&\text{for} && j\prec j^{(m,t)}.
\end{aligned}
\end{equation}
We rescale the coordinates $(X,U,Y,V)$ as
\begin{equation*}
\begin{aligned} 
X_1&=\lambda^{k_{\alpha}}\cdot X_1^{new}, &Y_1&=\lambda^{2k_{\alpha}}\cdot\gamma^{k_{\alpha}}\cdot Y_1^{new},\\
X_2&=\varrho_{k_{\alpha}}\cdot \lambda^{k_{\alpha}} \cdot X_1^{new}, &Y_2&=\varrho_{k_{\alpha}}\cdot\lambda^{2k_{\alpha}}\cdot\gamma^{k_{\alpha}}\cdot Y_1^{new},\\
U&=\frac{1}{\varrho_{k_{\alpha}}}\cdot\lambda^{k_{\alpha}}\cdot X_1^{new}, &V&=\frac{1}{\varrho_{k_{\alpha}}}\cdot\lambda^{2k_{\alpha}}\cdot\gamma^{k_{\alpha}}\cdot Y_1^{new},\\
\end{aligned}
\end{equation*}
and the parameters $\varepsilon,\tilde\varepsilon$ as
\begin{equation}\label{lemmaincreasingparatangencyindexeq21}
\begin{aligned} 
\varepsilon_{0,j}^*&=\lambda^{k_{\alpha}}\cdot E_{0,j} && \text{and} &&\varepsilon^*_{1,j}=\lambda^{k_{\alpha}}\cdot E_{1,j} && \text{for}  \;\; j\prec j^{(m,t)},\\
\tilde\varepsilon_{1}^*&=\gamma^{-k_{\alpha}}\cdot\widetilde E_{1},
\end{aligned}
\end{equation}
where $\varrho_{k_{\alpha}}$ tends sufficiently slowly to zero as $\alpha\to+\infty$. In this way, using condition \textbf{\hyperlink{C4.}{C4}}, we reduce equations \eqref{lemmaincreasingparatangencyindexeq14} and \eqref{lemmaincreasingparatangencyindexeq15} to
\begin{equation}\label{lemmaincreasingparatangencyindexeq23}
\begin{aligned} 
Y_1&=H_0+H_1\cdot X_1+a\cdot X_1^{2}+\vartheta_{1,k_{\alpha}}(X_1,Y,V,\zeta,E,\widetilde E),\\
(X_2,U)&=\vartheta_{2,k_{\alpha}}(X_1,Y,V,\zeta,E,\widetilde E),
\end{aligned}
\end{equation}
and
\begin{equation}\label{lemmaincreasingparatangencyindexeq24}
\begin{aligned} 
X_1&=\widetilde{H}_0+\widetilde{H}_1\cdot Y_1+\tilde\vartheta_{1,k_{\alpha}}(X,U,Y_1,\zeta,E,\widetilde E),\\
(Y_2,V)&=\tilde\vartheta_{2,k_{\alpha}}(X,U,Y_1,\zeta,E,\widetilde E),
\end{aligned}
\end{equation}
where the functions $\vartheta_{k_{\alpha}},\tilde\vartheta_{k_{\alpha}}$, along with their derivatives up to any given order, uniformly tend to zero for all sufficiently small $\zeta$, on any compact set of the rescaled variables $(X,U,Y,V)$ and of the rescaled parameters $E,\widetilde E$ as $\alpha\to+\infty$. Moreover, $\vartheta_{1,k_{\alpha}}$ does not contain constant and linear terms in $X_1$, and $\tilde\vartheta_{1,k_{\alpha}}$ does not contain constant and linear terms in $Y_1$. The splitting coefficients $H,\widetilde{H}$ are specified by
\begin{equation}\label{lemmaincreasingparatangencyindexeq25}
\begin{aligned}
H_0&=\lambda^{-k_{\alpha}}\cdot\sum\limits_{j\prec j^{(m,t)}}E_{0,j}\cdot\zeta^j+\lambda^{-2k_{\alpha}}\cdot (c_{0}+O(\gamma^{-k_{\alpha}}))\cdot\zeta^{j^{(m,t)}}+O_{(m,t+1)}(\zeta),\\
H_1&=\sum\limits_{j\prec j^{(m,t)}}E_{1,j}\cdot\zeta^j+\lambda^{-k_{\alpha}}\cdot (c_{1}+O(\gamma^{-k_{\alpha}}))\cdot\zeta^{j^{(m,t)}}+O_{(m,t+1)}(\zeta),
\end{aligned}
\end{equation}
and
\begin{equation}\label{lemmaincreasingparatangencyindexeq26}
\begin{aligned}
\widetilde{H}_0&=\lambda^{-2k_{\alpha}}\cdot (\tilde c_{0}+O(\lambda^{k_{\alpha}}))\cdot\zeta^{j^{(m,t)}}+O_{(m,t+1)}(\zeta),\\
\widetilde{H}_1&=\widetilde E_1+(\lambda\gamma)^{k_{\alpha}}\cdot O(\zeta)+\gamma^{k_{\alpha}}\cdot(\tilde c_{1}+O(\lambda^{k_{\alpha}}))\cdot\zeta^{j^{(m,t)}}+O_{(m,t+1)}(\zeta).
\end{aligned}
\end{equation}
For each $\alpha$, we consider these equations for $\zeta$ so small that the splitting coefficients $H,\widetilde{H}$ stay uniformly bounded. Note that, like the other coefficients in \eqref{lemmaincreasingparatangencyindexeq25},\eqref{lemmaincreasingparatangencyindexeq26}, the $O_{(m,t+1)}(\zeta)$-terms may tend to infinity as $\alpha\to+\infty$, but for each fixed $\alpha$ they decay to zero as $\zeta\to0$. 

Expressing $Y_1$ in \eqref{lemmaincreasingparatangencyindexeq24}, we obtain
\begin{equation}\label{lemmaincreasingparatangencyindexeq27}
\begin{aligned} 
Y_1&=-\widetilde{H}_1^{-1}\cdot\widetilde{H}_0+\widetilde H_1^{-1}\cdot X_1+\tilde\vartheta^*_{1,k_{\alpha}}(X,U,\zeta,E,\widetilde E,\widetilde{H}_0,\widetilde{H}_1),\\
(Y_2,V)&=\tilde\vartheta^*_{2,k_{\alpha}}(X,U,\zeta,E,\widetilde E,\widetilde{H}_0,\widetilde{H}_1),
\end{aligned}
\end{equation}
where the functions $\tilde\vartheta^*_{k_{\alpha}}$ with their derivatives decay to zero as $\alpha\to+\infty$ provided $\widetilde{H}_0$ stays bounded and $\widetilde{H}_1$ stays bounded from zero; the function $\tilde\vartheta^*_{1,k_{\alpha}}$ does not contain constant and linear terms in $X_1$. Applying the change of variables
\begin{equation*}
\begin{aligned} 
Y_1^{new}&=Y_1+\widetilde H_1^{-1}\cdot\widetilde{H}_0-\widetilde H_1^{-1}\cdot X_1-\tilde\vartheta^*_{1,k_{\alpha}}(X,U,\zeta,E,\widetilde E,\widetilde{H}_0,\widetilde{H}_1),\\
(Y_2,V)^{new}&=(Y_2,V)-\tilde\vartheta^*_{2,k_{\alpha}}(X,U,\zeta,E,\widetilde E,\widetilde{H}_0,\widetilde{H}_1),
\end{aligned}
\end{equation*}
we transform \eqref{lemmaincreasingparatangencyindexeq23} and \eqref{lemmaincreasingparatangencyindexeq27} into the form
\begin{equation}\label{lemmaincreasingparatangencyindexeq30}
\begin{aligned} 
Y_1&=H^*_0+H^*_1\cdot X_1+a\cdot X_1^{2}+\vartheta^*_{1,k_{\alpha}}(X_1,Y_2,V,\zeta,E,\widetilde E,H_0,H_1,\widetilde{H}_0,\widetilde{H}_1),\\
(X_2,U)&=\vartheta^*_{2,k_{\alpha}}(X_1,Y_2,V,\zeta,E,\widetilde E,H_0,H_1,\widetilde{H}_0,\widetilde{H}_1),
\end{aligned}
\end{equation}
and
\begin{equation*}
\begin{aligned} 
(Y,V)&=0,
\end{aligned}
\end{equation*}
where the functions $\vartheta^*_{k_{\alpha}}$ with their derivatives tend to zero as $\alpha\to+\infty$; the function $\vartheta^*_{1,k_{\alpha}}$ does not contain constant and linear terms in $X_1$. The coefficients $H^*$ satisfy
\begin{equation*}
\begin{aligned} 
H^*_0&=H_0+\widetilde H_1^{-1}\cdot\widetilde{H}_0+\sum\limits_{j\prec j^{(m,t+1)}}\Delta_{0,j}(E,\widetilde E,H_0,H_1,\widetilde H_0,\widetilde H_1)\cdot\zeta^j+O_{(m,t+1)},\\
H^*_1&=H_1-\widetilde H_1^{-1}+\sum\limits_{j\prec j^{(m,t+1)}}\Delta_{1,j}(E,\widetilde E,H_0,H_1,\widetilde H_0,\widetilde H_1)\cdot\zeta^j+O_{(m,t+1)},
\end{aligned}
\end{equation*}
where the small terms $\Delta_{i,j}$ tend to zero with their derivatives for all small $\zeta$, on any bounded set of $E,\widetilde E,H_0,H_1,\widetilde H_0,\widetilde H_1$ (with $\widetilde H_1$ staying bounded away from zero) as $\alpha\to+\infty$. We make a shift 
\begin{equation}\label{lemmaincreasingparatangencyindexeq30.75}
\begin{aligned} 
X_1&=X_1^{new}+\frac{1}{2a}\cdot\lambda^{-k_{\alpha}}\cdot E_{1,j^{(m,t)}}\cdot\zeta^{j^{(m,t)}},
\end{aligned}
\end{equation}
where $E_{1,j^{(m,t)}}$ is a bounded coefficient to be determined. Then, equations \eqref{lemmaincreasingparatangencyindexeq30} keep its form: 
\begin{equation}\label{lemmaincreasingparatangencyindexeq31}
\begin{aligned} 
Y_1&=H^{**}_0+H^{**}_1\cdot X_1+a\cdot X_1^{2}+\vartheta^{**}_{1,k_{\alpha}}(X_1,Y_2,V,\zeta,E,\widetilde E,H_0,H_1,\widetilde{H}_0,\widetilde{H}_1),\\
(X_2,U)&=\vartheta^{**}_{2,k_{\alpha}}(X_1,Y_2,V,\zeta,E,\widetilde E,H_0,H_1,\widetilde{H}_0,\widetilde{H}_1),
\end{aligned}
\end{equation}
with the new splitting coefficients acquiring $O_{(m,t)}$ corrections, namely,  
\begin{equation*}
\begin{aligned} 
H^{**}_0&=H_0+\widetilde H_1^{-1}\cdot\widetilde{H}_0+\lambda^{-k_{\alpha}}\cdot H_1^{**}\cdot E_{1,j^{(m,t)}}\cdot\zeta^{j^{(m,t)}}+\\
&+\sum\limits_{j\prec j^{(m,t+1)}}\Delta_{0,j}(E,\widetilde E,H_0,H_1,\widetilde H_0,\widetilde H_1)\cdot\zeta^j+O_{(m,t+1)},\\
H^{**}_1&=H_1-\widetilde H_1^{-1}+\lambda^{-k_{\alpha}}\cdot E_{1,j^{(m,t)}}\cdot\zeta^{j^{(m,t)}}+\\
&+\sum\limits_{j\prec j^{(m,t+1)}}\Delta_{1,j}(E,\widetilde E,H_0,H_1,\widetilde H_0,\widetilde H_1)\cdot\zeta^j+O_{(m,t+1)}.
\end{aligned}
\end{equation*}
By Definitions \ref{definitionindexofparatangencyfunctional} and \ref{definitionindexofparatangency}, we ensure the required paratangency at $\zeta=0$ to the bifurcation stratum that corresponds to a quadratic homoclinc tangency at $X_1=0$ (in the shifted coordinates \eqref{lemmaincreasingparatangencyindexeq30.75}) by eliminating the coefficients of $\zeta^j$ and $\zeta^j\cdot X_1$ for $j\prec j^{(m,t+1)}$ in the first equation of \eqref{lemmaincreasingparatangencyindexeq31}, i.e., by making $H_0^{**},H_1^{**}$ vanish for $j\prec j^{(m,t+1)}$. This gives us the equations
\begin{equation}\label{lemmaincreasingparatangencyindexeq33}
\begin{aligned} 
H_0+\widetilde H_1^{-1}\cdot\widetilde{H}_0+\sum\limits_{j\prec j^{(m,t+1)}}\Delta_{0,j}(E,\widetilde E,H_0,H_1,\widetilde H_0,\widetilde H_1)\cdot\zeta^j&=O_{(m,t+1)},\\
H_1-\widetilde H_1^{-1}+\lambda^{-k_{\alpha}}\cdot E_{1,j^{(m,t)}}\cdot\zeta^{j^{(m,t)}}+\sum\limits_{j\prec j^{(m,t+1)}}\Delta_{1,j}(E,\widetilde E,H_0,H_1,\widetilde H_0,\widetilde H_1)\cdot\zeta^j&=O_{(m,t+1)}.
\end{aligned}
\end{equation}
By the Implicit Function Theorem, we rewrite \eqref{lemmaincreasingparatangencyindexeq33} as
\begin{equation}\label{lemmaincreasingparatangencyindexeq34}
\begin{aligned} 
H_0+\widetilde H_1^{-1}\cdot\widetilde{H}_0+\lambda^{-k_{\alpha}}\cdot \delta_{0,j}(E,\widetilde E,\widetilde H_0,\widetilde H_1)\cdot E_{1,j^{(m,t)}}\cdot\zeta^{j^{(m,t)}}+\\
+\sum\limits_{j\prec j^{(m,t+1)}}\Delta^*_{0,j}(E,\widetilde E,\widetilde H_0,\widetilde H_1)\cdot\zeta^j=O_{(m,t+1)},\\
H_1-\widetilde H_1^{-1}+\lambda^{-k_{\alpha}}\cdot(1+\delta_{1,j}(E,\widetilde E,\widetilde H_0,\widetilde H_1))\cdot E_{1,j^{(m,t)}}\cdot\zeta^{j^{(m,t)}}+\\
+\sum\limits_{j\prec j^{(m,t+1)}}\Delta^*_{1,j}(E,\widetilde E,\widetilde H_0,\widetilde H_1)\cdot\zeta^j=O_{(m,t+1)},
\end{aligned}
\end{equation}
where the functions $\Delta^*_{i,j}$ and $\delta_{i,j}$ tend to zero as $\alpha\to+\infty$ along with their derivatives for all small $\zeta$, on any bounded set of $E,\widetilde E,\widetilde H_0,\widetilde H_1$ with $\widetilde H_1$ staying bounded away from zero. We multiply the first equation of \eqref{lemmaincreasingparatangencyindexeq34} to $\lambda^{k_{\alpha}}$ and expand in $\zeta$ using \eqref{lemmaincreasingparatangencyindexeq25},\eqref{lemmaincreasingparatangencyindexeq26}. This gives the following system for $E,\widetilde E$: 
\begin{equation}\label{lemmaincreasingparatangencyindexeq35}
\begin{aligned} 
E_{0,j}&=o(1)_{\alpha\to+\infty} \;\; \text{for}  \;\; j\prec j^{(m,t)},\\
\lambda^{-k_{\alpha}}\cdot (c_{0}+O(\gamma^{-k_{\alpha}}))+\lambda^{-k_{\alpha}}\cdot \widetilde E_1^{-1}\cdot(\tilde c_{0}+O(\lambda^{k_{\alpha}}))&=o(1)_{\alpha\to+\infty}.
\end{aligned}
\end{equation}
Note that $\lambda$ and $\gamma$ are functions of $\zeta$, but in the last equation one should replace them by their values at $\zeta=0$. Expanding the second equation of \eqref{lemmaincreasingparatangencyindexeq34} in $\zeta$, we obtain
\begin{equation}\label{lemmaincreasingparatangencyindexeq36}
\begin{aligned} 
\widetilde E_1\cdot E_{1,0}&=1+o(1)_{\alpha\to+\infty},\\
\widetilde E_1\cdot E_{1,j}&=o(1)_{\alpha\to+\infty} \;\; \text{for}  \;\; 0\prec j\prec j^{(m,t)},\\
\lambda^{-k_{\alpha}}\cdot \widetilde E_1\cdot (c_1+O(\gamma^{-k_{\alpha}}))+\gamma^{k_{\alpha}}\cdot(\tilde c_{1}+O(\lambda^{k_{\alpha}}))\cdot E_{1,0}+\\
+\lambda^{-k_{\alpha}}\cdot\widetilde E_1\cdot(1+o(1)_{\alpha\to+\infty})\cdot E_{1,j^{(m,t)}}&=o(1)_{\alpha\to+\infty}.
\end{aligned}
\end{equation}
In \eqref{lemmaincreasingparatangencyindexeq35} and \eqref{lemmaincreasingparatangencyindexeq36}, the terms $o(1)_{\alpha\to+\infty}$ stand for functions of $E,\widetilde E$ that tend to zero as $\alpha\to+\infty$, along with their derivatives, uniformly on any bounded set of $E,\widetilde E$ with $\widetilde E_1$ bounded away from zero. 

The solution of \eqref{lemmaincreasingparatangencyindexeq35},\eqref{lemmaincreasingparatangencyindexeq36} is 
\begin{equation*}
\begin{aligned} 
\widetilde E_1&=-\frac{\tilde c_{0}}{c_{0}}+o(1)_{\alpha\to+\infty}, &E_{1,0}&=-\frac{c_{0}}{\tilde c_{0}}+o(1)_{\alpha\to+\infty}, \;\;\; E_{1,j^{(m,t)}}=-c_1+o(1)_{\alpha\to+\infty},\\
E_{0,j}&=o(1)_{\alpha\to+\infty} \;\; \text{for} \;\; j\prec j^{(m,t)}, &E_{1,j}&=o(1)_{\alpha\to+\infty} \;\; \text{for}  \;\; 0\prec j\prec j^{(m,t)}.
\end{aligned}
\end{equation*}
Thus, the sought paratangency of index at least $(m,t+1)$ exists for these values of $E,\widetilde E$.

By \eqref{lemmaincreasingparatangencyindexeq11},\eqref{lemmaincreasingparatangencyindexeq18.5}, and \eqref{lemmaincreasingparatangencyindexeq21}, the values $\varepsilon_{\alpha},\tilde\varepsilon_{\alpha}$ corresponding to the found paratangency tend to zero exponentially fast as $\alpha\to+\infty$. Therefore, there exists a $\langle k_{\alpha}\rangle$-tame sequence $\langle\mathfrak{D}_{\alpha}\rangle_{\alpha\in\mathbb N}$ such that, for all sufficiently large $\alpha$, one has $(\zeta,\varepsilon_{\alpha},\tilde\varepsilon_{\alpha})\in\mathfrak{D}_{\alpha}$. 
\end{demo}

\begin{remark}
One can see from the proof that it remains valid for every parametric family $(f_{\zeta,\varepsilon,\tilde\varepsilon,\varepsilon'})$ of perturbations (that also depends on certain additional parameters $\varepsilon'$) if formulas \eqref{lemmaincreasingparatangencyindexeq4}--\eqref{lemmaincreasingparatangencyindexeq7} for $\mathcal W^u(O)$ and $\widetilde{\mathcal W}^s(O)$ are satisfied, i.e., for every family which perturbs the paratangency functionals $\eta_0,\eta_1,\tilde\eta_0,\tilde\eta_1$ freely.
\end{remark}

By this remark and Lemma \ref{lemmaperturbationofparatangencyfunctionals}, we immediately obtain 
\begin{lemma}\label{lemmaincreasingparatangencyindexmanyorbits}
Let $f_0:\mathcal M^d\rightarrow \mathcal M^d$ be a $C^r$ diffeomorphism, where $d\ge 2$ and $r=\infty,\omega$. Let $f_0$ have a saddle periodic orbit $L$ and an even number $h$ of orbits $\Gamma_1,\ldots,\Gamma_h$ of corank-$1$ quadratic homoclinic tangency between $W^s(L)$ and $W^u(L)$. Let $f_0$ be embedded into an $\ell$-parameter $C^r$ family $(f_{\zeta})$ which has a paratangency of index $(m,t)$ with a bifurcation stratum $\mathcal Y(\Gamma_1,\ldots, \Gamma_h)$, where $\ell\ge 1,m\ge0,t\ge1$. Then, $C^{r}$-close to $(f_{\zeta})$, there exists a parametric family $(g_{\zeta})$ which, at $\zeta=0$, has a paratangency of index at least $(m,t+1)$ with the bifurcation stratum corresponding to $\frac{h}{2}$ corank-$1$ quadratic homoclinic tangencies between $W^s(L)$ and $W^u(L)$. 
\end{lemma}

Adding, if necessary, a $C^r$-small perturbation to the parametric family $(g_{\zeta})$, we make these paratangencies be exactly of index $(m,t+1)$, i.e., the corresponding coefficients $c_{j^{(m,t+1)}}\not=0$ (see Definition \ref{definitionindexofparatangencyfunctional}).

\begin{lemma}\label{lemmaincreasingparatangencyorder}
Let $f_0:\mathcal M^d\rightarrow \mathcal M^d$ be a $C^r$ diffeomorphism, where $d\ge 2$ and $r=\infty,\omega$. Let $f_0$ have a saddle periodic orbit $L$ and $h = 2^{\raisebox{0.5ex}{$\binom{\ell+m}{m} - 1$}}$ orbits $\Gamma_1,\ldots,\Gamma_h$ of corank-$1$ quadratic homoclinic tangency between $W^s(L)$ and $W^u(L)$. Let $f_0$ be embedded into an $\ell$-parameter $C^r$ family $(f_{\zeta})$ which has a paratangency of index $(0,1)$ with a bifurcation stratum $\mathcal Y(\Gamma_1,\ldots,\Gamma_h)$, where $\ell\ge 1,m\ge0,t\ge1$. Then, $C^{r}$-close to $(f_{\zeta})$, there exists a parametric family $(g_{\zeta})$ which, at $\zeta=0$, has a paratangency of index $(m,1)$ with the bifurcation stratum corresponding to corank-$1$ quadratic homoclinic tangency between $W^s(L)$ and $W^u(L)$. 
\end{lemma}

\begin{demo}
First, observe that the number of multi-indices of degree from $1$ to $m$ with $\ell$ entries is 
$$p=\binom{\ell+m}{m} - 1.$$ 

If $m=0$, then the diffeomorphism $f_0$ has one orbit of corank-$1$ quadratic homoclinic tangency $\Gamma_1$, and the family $(f_{\zeta})$, at $\zeta=0$, has a paratangency of index $(0,1)$ with a bifurcation stratum $\mathcal Y(\Gamma_1)$. The lemma is trivially true.

Assume $m>0$. By Lemma \ref{lemmaincreasingparatangencyindexmanyorbits}, close in $C^{r}$ to $(f_{\zeta})$, there exists a parametric family $(f^1_{\zeta})$ which, at $\zeta=0$, has a paratangency of index $(1,1)$ with the bifurcation stratum corresponding to $\frac{h}{2}$ corank-$1$ quadratic homoclinic tangencies $\Gamma^1_1,\ldots,\Gamma^1_{\frac{h}{2}}$ between $W^s(L)$ and $W^u(L)$. Applying Lemma \ref{lemmaincreasingparatangencyindexmanyorbits} again, we find a parametric family $(f^2_{\zeta})$ which is close in $C^{r}$ to $(f^1_{\zeta})$ and, at $\zeta=0$, has a paratangency of index $(1,2)$ with the bifurcation stratum corresponding to $\frac{h}{4}$ corank-$1$ quadratic homoclinic tangencies $\Gamma^2_1,\ldots,\Gamma^2_{\frac{h}{4}}$ between $W^s(L)$ and $W^u(L)$. Repeating this procedure $p$ times, we obtain a parametric family $(g_{\zeta})\equiv(f^p_{\zeta})$ which is close in $C^{r}$ to $(f_{\zeta})$ and, at $\zeta=0$, has a paratangency of index $(m,1)$ with the bifurcation stratum corresponding to $\frac{h}{2^p}=1$ corank-$1$ quadratic homoclinic tangency between $W^s(L)$ and $W^u(L)$.
\end{demo}

Recall that, by Lemma \ref{lemma1highordertangencies}, any bifurcation stratum corresponding to the simultaneous existence of a pair of corank-1 homoclinic tangencies, one quadratic and one of order $n$, is accumulated by a sequence of bifurcation strata corresponding to a corank-1 homoclinic tangency of order $n+1$. This, in an obvious way, leads to
\begin{lemma}\label{lemmakeepingparatangencyorder}
Let $f_0:\mathcal M^d\rightarrow \mathcal M^d$ be a $C^r$ diffeomorphism, where $d\ge 2$ and $r=\infty,\omega$. Let $f_0$ have a saddle periodic orbit $L$ and two orbits $\Gamma,\widetilde\Gamma$ of corank-1 homoclinic tangency between $W^s(L)$ and $W^u(L)$. Let $\Gamma$ be of order $n$ and $\widetilde\Gamma$ be quadratic. Let $f_0$ be embedded into an $\ell$-parameter $C^r$ family $(f_{\zeta})$ which has a paratangency of order $m$ with a bifurcation stratum $\mathcal Y(\Gamma,\widetilde\Gamma)$, where $\ell\ge 1,m\ge0$. Then, $C^{r}$-close to $(f_{\zeta})$, there exists a parametric family $(g_{\zeta})$ which, at $\zeta=0$, has a paratangency of order $m$ with the bifurcation stratum corresponding to corank-1 homoclinic tangency of order $n+1$ between $W^s(L)$ and $W^u(L)$. 
\end{lemma}

\begin{demo}
Employing Lemma \ref{lemmaperturbationofparatangencyfunctionals}, we find a $C^r$ family $(f_{\zeta,\varepsilon})$ such that $(f_{\zeta,0})=(f_{\zeta})$ and condition \eqref{lemma1highordertangencieseq0} holds. By this condition, we can choose parameters $\varepsilon$ such that 
\begin{equation}\label{lemmakeepingparatangencyordereq1}
\begin{aligned} 
\eta(\zeta,\varepsilon)&=\varepsilon+\eta(\zeta,0).
\end{aligned}
\end{equation}
By Lemma \ref{lemma1highordertangencies}, there exists a small $\varepsilon_{\alpha}(\zeta)$ corresponding to the existence of a corank-1 homoclinic tangency $\Gamma'$ of order $(n+1)$. The relation between the splitting coefficients $\eta$ for the homoclinic tangencies $\Gamma,\widetilde\Gamma$ and the splitting coefficients $\nu$ for the tangency $\Gamma'$ is given by Remark \ref{remarkhighordertangencies}. We define the sought $C^r$-small perturbation $(g_{\zeta})$ of the family $(f_{\zeta})$ as $(g_{\zeta})\equiv (f_{\zeta,Q_{\alpha}(\zeta,0)})$, where $Q_{\alpha}$ is the function from \eqref{remarkhighordertangencieseq1}. Thus, letting $\varepsilon=Q_{\alpha}(\zeta,0)$ in \eqref{lemmakeepingparatangencyordereq1}, we obtain from \eqref{remarkhighordertangencieseq1} that
\begin{equation}\label{lemmakeepingparatangencyordereq2}
\begin{aligned} 
\eta(\zeta,0)&=Q_{\alpha}(\zeta,\nu)-Q_{\alpha}(\zeta,0),
\end{aligned}
\end{equation}
where $\nu$ are the splitting coefficients for the tangency $\Gamma'$ in the family $(g_{\zeta})$. Since $\left.\det\frac{\partial Q_{\alpha}}{\partial\nu}\right|_{(\zeta,\nu)=0}\not=0$, formula \eqref{lemmakeepingparatangencyordereq2} entails that the Taylor expansion of $\nu$ by $\zeta$ starts with the terms of order not lower than the expansion of $\eta(\zeta,0)$ does.
\end{demo}

Lemmas \ref{lemma2highordertangencies} and \ref{lemmakeepingparatangencyorder} entail
\begin{lemma}\label{lemmahighordertangenciesandparatangencies}
Let $f_0:\mathcal M^d\rightarrow \mathcal M^d$ be a $C^r$ diffeomorphism, where $d\ge 2$ and $r=\infty,\omega$. Let $f_0$ have a saddle periodic orbit $L$ and $n$ orbits $\Gamma_1,\ldots,\Gamma_n$ of corank-1 quadratic homoclinic tangency between $W^s(L)$ and $W^u(L)$. Let $f_0$ be embedded into an $\ell$-parameter $C^r$ family $(f_{\zeta})$ which has a paratangency of order $m$ with a bifurcation stratum $\mathcal Y(\Gamma_1,\ldots,\Gamma_n)$, where $\ell\ge 1,m\ge0$. Then, $C^{r}$-close to $(f_{\zeta})$, there exists a parametric family $(g_{\zeta})$ which, at $\zeta=0$, has a paratangency of order $m$ with the bifurcation stratum corresponding to corank-1 homoclinic tangency of order $n$ between $W^s(L)$ and $W^u(L)$. 
\end{lemma}

\begin{remark}
Lemmas \ref{lemmaperturbationofparatangencyfunctionals}--\ref{lemmahighordertangenciesandparatangencies} also hold true for sufficiently large finite $r$.
\end{remark}

\begin{demo5}
Let $(f_{\zeta})$ be an $\ell$-parameter $C^r$ family such that $(f_{\zeta})\in\mathcal{PN}(L)$ with $f_{\zeta^*}\in\mathcal N(L)$. Without loss of generality, we assume $\zeta^*=0$. Also, we assume that $(f_{\zeta})$ is real-analytic, since otherwise it can always be $C^r$-approximated by a real-analytic family $(f_{\zeta})^{new}\in\mathcal{PN}(L)$.

By Lemma \ref{lemmamanytangencies}, $C^{\omega}$-close to $(f_{\zeta})$, there exists a parametric family $(\hat f_{\zeta})$ such that $\hat f_0$ has $h =n\cdot 2^{\raisebox{0.5ex}{$\binom{\ell+m}{m} - 1$}}$ orbits $\widehat\Gamma_1,\ldots,\widehat\Gamma_h$ of corank-$1$ quadratic homoclinic tangency between $W^s(L)$ and $W^u(L)$. Consecutively applying Lemma \ref{lemmaincreasingparatangencyorder} to collections $\{\widehat\Gamma_1,\ldots,\widehat\Gamma_{\frac{h}{n}}\},\ldots,\{\widehat\Gamma_{\frac{h(n-1)}{n}+1},\ldots,\widehat\Gamma_{h}\}$, $C^{\omega}$-close to $(\hat f_{\zeta})$, we find a parametric family $(\overline f_{\zeta})$ such that $\overline f_{0}$ has $n$ orbits $\overline\Gamma_1,\ldots,\overline\Gamma_n$ of corank-1 quadratic homoclinic tangency and $(\overline f_{\zeta})$, at $\zeta=0$, has a paratangency of order $m$ with a bifurcation stratum $\mathcal Y(\overline\Gamma_1,\ldots,\overline\Gamma_n)$. Making use of Lemma \ref{lemmahighordertangenciesandparatangencies}, $C^{\omega}$-close to $(\overline f_{\zeta})$, we obtain a parametric family $(g_{\zeta})$ which, at $\zeta=0$, has a paratangency of order $m$ with the bifurcation stratum corresponding to corank-1 homoclinic tangency of order $n$ between $W^s(L)$ and $W^u(L)$. 
\end{demo5}

\section{Universal dynamics}\label{sectionuniversaldynamics}

In this section, we derive Theorems \ref{theorem6} and \ref{theorem7} from Theorem \ref{theorem4} (see Subsections \ref{subsectionfirstreturnmaps} and \ref{subsectiondynamicsnearhomoclinicdisks}, respectively).

\subsection{Formulas for the global map}\label{subsectionformulasforglobalmap}
Let $f_0:\mathcal M^d\rightarrow \mathcal M^d$ be a $C^r$ diffeomorphism, where $d\ge 2$ and $r=2,.\ldots,\infty,\omega$. Let $f$ have a saddle periodic orbit $L$ whose stable and unstable manifolds $W^s(L)$ and $W^u(L)$ intersect over an orbit $\Gamma$ of corank-$c$ homoclinic tangency. Let $L$ satisfy condition \textbf{\hyperlink{C1.}{C1}}, and $\Gamma$ satisfy condition \textbf{\hyperlink{C2.}{C2}}. 

Suppose $f_0$ is embedded into a finite-parameter $C^r$ family $(f_{\varepsilon})$. We choose a point $O\in L$ and fix a small neighborhood $U(O)$, where the local map $T_{O}$ is well-defined for all small $\varepsilon$. We bring $T_O$ to the main normal form. As in Section \ref{subsectionglobalmap}, for all small $\varepsilon$, we define the global map $T_{\Gamma}:\Pi^-\to\Pi^+$, where $\Pi^+$ and $\Pi^-$ are small neighborhoods of the homoclinic points $M^{+},M^{-}\in\Gamma$ such that $M^{+}\in W^s_{loc}(O)$ and $M^{-}\in W^u_{loc}(O)$. 

Further, we apply a series of coordinate transformations to reduce the global map $T_{\Gamma}$ to the specific form detailed in Lemma \ref{lemmaglobalmap} below. All these transformations are linear, so they keep the local map $T_O$ in the main normal form. 

\textit{Step 1.} The manifolds $W^s_{loc}(O)$ and $W^u_{loc}(O)$ are straightened, so the homoclinic points $M^+$ and $M^-$ are given by $M^+=(x^+,u^+,0,0)$ and $M^-=(0,0,\tilde y^-,\tilde v^-)$. Since $T_{\Gamma}(M^-)=M^+$ at $\varepsilon=0$, we can write the global map $T_{\Gamma}$  for all small $\varepsilon$ in the form 
\begin{equation}\label{globalmapeq1}
\begin{aligned} 
\overline x -x^+&=P_1\left(\tilde x,\tilde u,\tilde y-\tilde y^-,\tilde v-\tilde v^-,\varepsilon\right), \quad & \overline y &= P_3\left(\tilde x,\tilde u,\tilde y-\tilde y^-,\tilde v-\tilde v^-,\varepsilon\right),\\
\overline u-u^+&=P_2\left(\tilde x,\tilde u,\tilde y-\tilde y^-,\tilde v-\tilde v^-,\varepsilon\right), \quad &\overline v&= P_4\left(\tilde x,\tilde u,\tilde y-\tilde y^-,\tilde v-\tilde v^-,\varepsilon\right),\\
\end{aligned}
\end{equation}
where $P$ are $C^r$-smooth functions. 

\textit{Step 2.} By condition \textbf{\hyperlink{C2.2.}{C2.2}}, the manifold $W^{se}_{loc}(O)$ and the leaf $T_{\Gamma}(l^{uu})$ are transverse at the point $M^-$ at $\varepsilon=0$. Note that $l^{uu}$ is given by $\{\tilde y=\tilde y^-,(\tilde x,\tilde u)=0\}$ (see \eqref{foliationseq3}), and the tangent space $\mathcal T_{M^+}(W^{se}_{loc}(O))$ is given by $\{\Delta\overline v=0\}$ (see \eqref{foliationseq4}). Therefore, the transversality condition \textbf{\hyperlink{C2.2.}{C2.2}} implies 
\begin{equation*}
\begin{aligned} 
&\det\left.\frac{\partial P_4}{\partial \tilde v}\right|_{M^-}\not=0.
\end{aligned}
\end{equation*}
Thus, we rewrite \eqref{globalmapeq1} for all small $\varepsilon$ as
\begin{equation}\label{globalmapeq2}
\begin{aligned} 
\overline x -x^+&= Q_1\left(\tilde x,\tilde u,\tilde y-\tilde y^-,\overline v,\varepsilon\right), \quad &\overline y &= Q_3\left(\tilde x,\tilde u,\tilde y-\tilde y^-,\overline v,\varepsilon\right),\\
\overline u-u^+&= Q_2\left(\tilde x,\tilde u,\tilde y-\tilde y^-,\overline v,\varepsilon\right), \quad &\tilde v-\tilde v^-&= Q_4\left(\tilde x,\tilde u,\tilde y-\tilde y^-,\overline v,\varepsilon\right).\\
\end{aligned}
\end{equation}
Here, $Q$ are $C^r$-smooth functions and 
\begin{equation*}
\begin{aligned} 
&\det\left.\frac{\partial Q_4}{\partial \overline v}\right|_{\widehat M}\not=0,
\end{aligned}
\end{equation*}
where we denote $\widehat M=(\tilde x=0,\tilde u=0,\tilde y=\tilde y^-,\overline v=0)$.

\textit{Step 3.} Recall that $\Gamma$ is an orbit of corank-$c$ tangency, i.e., the manifolds $W^s_{loc}(O)$ and $T_\Gamma(W^u(O))$ have $c$ linearly independent common tangent vectors at the point $M^+$ at $\varepsilon=0$. Since $\mathcal T_{M^+}(W^s_{loc}(O))$ is given by $\{(\Delta\overline y, \Delta\overline v)=0\}$ and $\mathcal T_{M^-}(W^u_{loc}(O))$ is given by $\{(\Delta x,\Delta u)=0\}$, this is equivalent to
\begin{equation}\label{globalmapeq2.5}
\begin{aligned} 
&\text{corank}\, \left.\frac{\partial Q_3}{\partial \tilde y}\right|_{\widehat M}=c \quad \text{at}\quad \varepsilon=0.
\end{aligned}
\end{equation}
Let $\alpha$ be the zero eigenvalue of the matrix $\left.\frac{\partial Q_3}{\partial \tilde y}\right|_{\widehat M}$. The equality \eqref{globalmapeq2.5} means that $\text{geom}(\alpha)=c$. We require the following genericity condition to hold: 

\begin{doublespace}
\textbf{\hypertarget{C5.}{C5.}} $\text{alg}(\alpha)=\text{geom}(\alpha)=c$ at $\varepsilon=0$.
\end{doublespace}

Note that for the main cases $(1,1), (2,1), (1,2), (2,2)$ and $c=1$ this condition can be achieved by the choice of the homoclinic point $M^+$, see Lemma 8 in \cite{Gonchenko3}.

Let us perform a linear change of $y$-coordinates to bring the matrix $\left.\frac{\partial Q_3}{\partial \tilde y}\right|_{\widehat M}$ to the real Jordan form. By condition \textbf{\hyperlink{C5.}{C5}}, the system \eqref{globalmapeq2} for all small $\varepsilon$ is recast as 
\begin{equation}\label{globalmapeq3}
\begin{aligned} 
\overline x -x^+&=R_1\left(\tilde x,\tilde u,\tilde y-\tilde y^-,\overline v,\varepsilon\right), \quad &\overline y_1&= R_{31}\left(\tilde x,\tilde u,\tilde y-\tilde y^-,\overline v,\varepsilon\right),\\
&&\overline y_2&= R_{32}\left(\tilde x,\tilde u,\tilde y-\tilde y^-,\overline v,\varepsilon\right),\\
\overline u-u^+&=R_2\left(\tilde x,\tilde u,\tilde y-\tilde y^-,\overline v,\varepsilon\right), \quad &\tilde v-\tilde v^-&= R_4\left(\tilde x,\tilde u,\tilde y-\tilde y^-,\overline v,\varepsilon\right),\\
\end{aligned}
\end{equation}
where $\tilde y_1\in\mathbb R^c$, $\tilde y_2\in\mathbb R^{m_u-c}$, and
\begin{equation*}
\begin{aligned} 
\left.\frac{\partial R_{31}}{\partial\tilde y_1}\right|_{\widehat M}&=0 \quad \text{at}\quad \varepsilon=0,\\
\end{aligned}
\end{equation*}
\begin{equation*}
\begin{aligned} 
\left.\frac{\partial R_{31}}{\partial\tilde y_2}\right|_{\widehat M}&=0, \quad \left.\frac{\partial R_{32}}{\partial \tilde y_1}\right|_{\widehat M}=0, \quad \det\left.\frac{\partial R_{32}}{\partial \tilde y_2}\right|_{\widehat M}\not=0, \quad \det\left.\frac{\partial R_4}{\partial \overline v}\right|_{\widehat M}\not=0 \quad \text{for all small} \; \varepsilon.
\end{aligned}
\end{equation*}
By the Implicit Function Theorem, we rewrite \eqref{globalmapeq3} as 
\begin{equation*}
\begin{aligned} 
\overline x -x^+&= S_1\left(\tilde x,\tilde u,\tilde y_1-\tilde y_1^-,\overline y_2,\overline v,\varepsilon\right), \quad &\overline y_1 &= S_{31}\left(\tilde x,\tilde u,\tilde y_1-\tilde y_1^-,\overline y_2,\overline v,\varepsilon\right),\\
&&\tilde y_2-\tilde y_2^-&= S_{32}\left(\tilde x,\tilde u,\tilde y_1-\tilde y_1^-,\overline y_2,\overline v,\varepsilon\right),\\
\overline u-u^+&= S_2\left(\tilde x,\tilde u,\tilde y_1-\tilde y_1^-,\overline y_2,\overline v,\varepsilon\right), \quad &\tilde v-\tilde v^-&= S_4\left(\tilde x,\tilde u,\tilde y_1-\tilde y_1^-,\overline y_2,\overline v,\varepsilon\right).\\
\end{aligned}
\end{equation*}
The functions $S$ are $C^r$-smooth and satisfy
\begin{equation*}
\begin{aligned} 
\left.\frac{\partial S_{31}}{\partial\tilde y_1}\right|_{M}&=0, \quad \text{at}\quad \varepsilon=0,\\
\end{aligned}
\end{equation*}
\begin{equation*}
\begin{aligned} 
\left.\frac{\partial S_{31}}{\partial\overline y_2}\right|_{M}&=0, \quad \left.\frac{\partial S_{32}}{\partial \tilde y_1}\right|_{M}=0, \quad \det\left.\frac{\partial S_{32}}{\partial \overline y_2}\right|_{M}\not=0, \quad \det\left.\frac{\partial S_4}{\partial \overline v}\right|_{M}\not=0 \quad \text{for all small} \; \varepsilon,
\end{aligned}
\end{equation*}
where $M=(\tilde x=0,\tilde u=0,\tilde y_1=\tilde y_1^-,\overline y_2=0,\overline v=0)$. 

\textit{Step 4.} By condition \textbf{\hyperlink{C2.1.}{C2.1}}, the manifold $T_\Gamma(W^{ue}_{loc}(O))$ and the leaf $l^{ss}$ are transverse at $M^+$ at $\varepsilon=0$. The leaf $l^{ss}$ is given by $\{\overline x=x^+, (\overline y,\overline v)=0\}$ (see \eqref{foliationseq3}) and the tangent space $\mathcal T_{M^-}(W^{ue}_{loc}(O))$ is given by $\{\Delta \tilde u=0\}$ (see \eqref{foliationseq4}). Thus, the transversality entails
\begin{equation*}
\begin{aligned} 
&\det \begin{pmatrix} \left.\frac{\partial S_1}{\partial\tilde x}\right|_{M}& \left.\frac{\partial S_1}{\partial\tilde y_1}\right|_{M}\\ \left.\frac{\partial S_{31}}{\partial\tilde x}\right|_{M}& 0 \end{pmatrix}\not=0,
\end{aligned}
\end{equation*}
which immediately implies that the matrices $\left.\frac{\partial S_1}{\partial\tilde y_1}\right|_{M}$ and $\left.\frac{\partial S_{31}}{\partial\tilde x}\right|_{M}$ have full rank:
\begin{equation*}
\begin{aligned} 
&\text{rank}\,  \left.\frac{\partial S_1}{\partial\tilde y_1}\right|_{M}=c \quad \text{and} \quad \text{rank}\, \left.\frac{\partial S_{31}}{\partial\tilde x}\right|_{M}=c \quad \text{for all small } \varepsilon.
\end{aligned}
\end{equation*}

\textit{Step 5.} Next, we apply a linear change of $x$-coordinates to ensure that the minor determinant in the first $c$ rows of the matrix $\left.\frac{\partial S_1}{\partial\tilde y_1}\right|_{M}$ is non-zero and the remaining $m_s-c$ rows are nullified. 

Thus, we have proved
\begin{lemma}\label{lemmaglobalmap}
Under conditions \textbf{\hyperlink{C1.}{C1}},\textbf{\hyperlink{C2.}{C2}} and \textbf{\hyperlink{C5.}{C5}}, the system of coordinates in $U(O)$ can be chosen in such a way that, for all small $\varepsilon$, the local map $T_{O}$ is in the main normal form, and the global map $T_{\Gamma}:(\tilde x,\tilde u,\tilde y_1,\overline y_2,\overline v)\mapsto(\overline x,\overline u,\overline y_1,\tilde y_2,\tilde v)$ satisfies: 
\begin{equation}\label{lemmaglobalmapeq1}
\begin{aligned} 
\overline x_1 -x_1^{+}&=\theta_{11}(\tilde x,\tilde u,\tilde y_1-\tilde y_1^-,\overline y_2,\overline v,\varepsilon), \\
\overline y_1 &=F(\tilde y_1-\tilde y_1^-,\varepsilon)+\theta_{31}(\tilde x,\tilde u,\tilde y_1-\tilde y_1^-,\overline y_2,\overline v,\varepsilon),\\
\overline x_2 -x_{2}^+&=\theta_{12}(\tilde x,\tilde u,\tilde y_1-\tilde y_1^-,\overline y_2,\overline v,\varepsilon), \quad \tilde y_2-\tilde y_2^-=\theta_{32}(\tilde x,\tilde u,\tilde y_1-\tilde y_1^-,\overline y_2,\overline v,\varepsilon),\\
\overline u-u^+&= \theta_2(\tilde x,\tilde u,\tilde y_1-\tilde y_1^-,\overline y_2,\overline v,\varepsilon), \quad\;\;\; \tilde v-\tilde v^-= \theta_4(\tilde x,\tilde u,\tilde y_1-\tilde y_1^-,\overline y_2,\overline v,\varepsilon),\\
\end{aligned}
\end{equation}
where $(\tilde x_1,\tilde x_2,\tilde u,\tilde y_1,\tilde y_2,\tilde v)\in\mathbb R^{c}\times\mathbb R^{m_s-c}\times\mathbb R^{d_s-m_s}\times\mathbb R^{c}\times\mathbb R^{m_u-c}\times\mathbb R^{d_u-m_u}$. At $\varepsilon=0$, one has
\begin{equation*}
\begin{aligned} 
F(\tilde y_1=\tilde y_1^-,\varepsilon=0)&=0, &\left.\frac{\partial F}{\partial \tilde y_1}\right|_{(\tilde y_1=\tilde y_1^-,\varepsilon=0)}&=0.
\end{aligned}
\end{equation*}
For all small $\varepsilon$, one has 
\begin{equation}\label{lemmaglobalmapeq4}
\begin{aligned} 
\theta(0,0,\tilde y_1=\tilde y_1^-,0,0,\varepsilon)&=0, &\theta_{31}(0,0,\tilde y_1-\tilde y_1^-,0,0,\varepsilon)&=0,
\end{aligned}
\end{equation}
and
\begin{equation}\label{lemmaglobalmapeq5}
\begin{aligned} 
\left.\frac{\partial \theta_{12}}{\partial \tilde y_1}\right|_{M}&=0, &\left.\frac{\partial \theta_{31}}{\partial\overline y_2}\right|_{M}&=0, &\left.\frac{\partial \theta_{32}}{\partial \tilde y_1}\right|_{M}&=0,\\ 
\left.\frac{\partial \theta_{11}}{\partial\tilde y_1}\right|_{M}&\not=0, &\text{rank}\, \left.\frac{\partial \theta_{31}}{\partial\tilde x}\right|_{M}&=c, &\det\left.\frac{\partial \theta_{32}}{\partial \overline y_2}\right|_{M}&\not=0, &\det\left.\frac{\partial \theta_4}{\partial \overline v}\right|_{M}&\not=0,
\end{aligned}
\end{equation}
where $M=(\tilde x=0,\tilde u=0,\tilde y_1=\tilde y_1^-,\overline y_2=0,\overline v=0)$.
\end{lemma}

\begin{lemma}\label{lemmaglobalmaporder}
Let $\varepsilon=0$. The tangency $\Gamma$ is of order $n\le r-1$ if and only if 
\begin{equation}\label{lemmaglobalmapordereq1}
\begin{aligned} 
&\left.\frac{\partial^i F}{\partial\tilde y^i_1}\right|_{(\tilde y_1=\tilde y_1^-,\varepsilon=0)}=0 \quad \text{for} \quad 0\le |i|\le n, \quad \text{and} \quad \left.\frac{\partial^{i} F}{\partial\tilde y^{i}_1}\right|_{(\tilde y_1=\tilde y_1^-,\varepsilon=0)}\not=0 \;\; \text{for some} \;\; i \;\; \text{with} \;\; |i|=n+1.
\end{aligned}
\end{equation}
\end{lemma}

\begin{demo}
By setting $\tilde x=0$ and $\tilde u=0$ in \eqref{lemmaglobalmapeq1}, we obtain the following formula for $\mathcal W^u_{\Gamma}(O)\equiv T_{\Gamma}(W^u_{loc}(O))\cap\Pi^+$:
\begin{equation}\label{unstablemanifoldneartangencyeq1}
\begin{aligned} 
\overline x_1 -x_1^{+}&=\theta_{11}(0,0,\tilde y_1-\tilde y_1^-,\overline y_2,\overline v,\varepsilon), \\
\overline y_1 &=F(\tilde y_1-\tilde y_1^-,\varepsilon)+\theta_{31}(0,0,\tilde y_1-\tilde y_1^-,\overline y_2,\overline v,\varepsilon),\\
\overline x_2 -x_{2}^+&=\theta_{12}(0,0,\tilde y_1-\tilde y_1^-,\overline y_2,\overline v,\varepsilon), \quad \tilde y_2-\tilde y_2^-= \theta_{32}(0,0,\tilde y_1-\tilde y_1^-,\overline y_2,\overline v,\varepsilon),\\
\overline u-u^+&= \theta_2(0,0,\tilde y_1-\tilde y_1^-,\overline y_2,\overline v,\varepsilon), \quad\;\;\; \tilde v-\tilde v^-= \theta_4(0,0,\tilde y_1-\tilde y_1^-,\overline y_2,\overline v,\varepsilon).\\
\end{aligned}
\end{equation}
Resolving the first equation with respect to $(\tilde y_1-\tilde y_1^-)$ yields
\begin{equation}\label{unstablemanifoldneartangencyeq2}
\begin{aligned} 
\overline y_1 &=\widehat F(\overline x_1- x_1^+,\varepsilon)+\hat\theta_1(\overline x_1- x_1^+,\overline y_2,\overline v,\varepsilon),\\
\overline x_2 -x_2^+ &=\hat\theta_2(\overline x_1- x_1^+,\overline y_2,\overline v,\varepsilon),\\
\overline u-u^+&= \hat\theta_3(\overline x_1- x_1^+,\overline y_2,\overline v,\varepsilon),\\
\end{aligned}
\end{equation}
where
\begin{equation*}
\begin{aligned} 
&\widehat F(\overline x_1=x_1^+,\varepsilon=0)=0, \quad \left.\frac{\partial \hat F}{\partial\overline x_1}\right|_{(\overline x_1= x_1^+,\varepsilon=0)}=0, \quad \hat\theta_{1}(\overline x_1-x_1^+,0,0,\varepsilon)=0,\\
\end{aligned}
\end{equation*}
and the functions $\hat\theta$ vanish at $\overline x_1=x_1^+, \overline y_2=0, \overline v=0$.

The local stable manifold $W^{s}_{loc}(O)$ is straightened and given by $(\overline y,\overline v)=(0,0)$. According to Definition \ref{coordinatedefinitionorder}, the tangency $\Gamma$ is of order $n\le r-1$ at $\varepsilon=0$ if and only if the functions $\widehat F$ in \eqref{unstablemanifoldneartangencyeq2} satisfy
\begin{equation*}
\begin{aligned} 
&\left.\frac{\partial^i \widehat F}{\partial\overline x^i_1}\right|_{(\overline x_1= x_1^+,\varepsilon=0)}=0 \quad \text{for} \quad 0\le |i|\le n, \quad \text{and} \quad \left.\frac{\partial^{i} \widehat F}{\partial\overline x^{i}_1}\right|_{(\overline x_1= x_1^+,\varepsilon=0)}\not=0 \;\; \text{for some} \;\; i \;\; \text{with} \;\; |i|=n+1.
\end{aligned}
\end{equation*}
The relation between the coordinates $\tilde y_1$ and $\overline x_1$ given by the first equation of \eqref{unstablemanifoldneartangencyeq1} is a $C^r$ diffeomorphism. Since a $C^r$ diffeomorphism preserves the order of the first non-vanishing term in a Taylor expansion, these conditions on $\widehat F$ are equivalent to \eqref{lemmaglobalmapordereq1}. 
\end{demo}

\begin{remark}
Let $\varepsilon=0$. The tangency $\Gamma$ belongs to the homoclinic disk if and only if $F\equiv0$ in some neighborhood of $\tilde y_1=\tilde y_1^-$.
\end{remark}

Suppose the tangency $\Gamma$ is of order $n$. For sufficiently small $\varepsilon$, the function $F$ of \eqref{lemmaglobalmapeq1} can be expanded as 
\begin{equation*}
\begin{aligned} 
F(\tilde y_1-\tilde y_1^-,\varepsilon)&=\sum\limits_{|i|=0}^{n} \mu_{i} \cdot \left(\tilde y_1-\tilde y_1^-\right)^{i}+O(||\tilde y_1-\tilde y_1^-||^{n+1}),\\
\end{aligned}
\end{equation*}
where the coefficients $\mu_i$ are smooth functions of $\varepsilon$ that vanish at $\varepsilon=0$ by Lemma \ref{lemmaglobalmaporder}. Denote by $\mathbf M_{n}$ the vector of all coefficients $\mu_{i}$ from this expansion.
\begin{lemma}\label{lemmaglobalmapsplitting}
The tangency $\Gamma$ splits freely up to order $n$ if and only if
\begin{equation}\label{lemmaglobalmapsplittingeq1} 
\left.\text{rank}\;\frac{\partial\mathbf M_{n}}{\partial\varepsilon}\right |_{\varepsilon=0}=\dim\mathbf M_{n}. 
\end{equation}
\end{lemma}

\begin{demo}
By Definition \ref{definitionfreeperturbationorbit}, the tangency $\Gamma$ splits freely up to order $n$ if and only if $\left.\text{rank}\;\frac{\partial\mathbf H_{n}}{\partial\varepsilon}\right |_{\varepsilon=0}=\dim\mathbf H_{n}$, where $\mathbf H_{n}$ is the vector of all Taylor coefficients $\eta_i$ of $\widehat F$ (see \eqref{unstablemanifoldneartangencyeq2}) up to order $n$:
\begin{equation*}
\begin{aligned} 
\widehat F(\overline x_1- x_1^+,\varepsilon)&=\sum\limits_{|i|=0}^{n} \eta_{i} \cdot (\overline x_1- x_1^+)^{i}+O(||\overline x_1- x_1^+||^{n+1}).
\end{aligned}
\end{equation*}
Comparing \eqref{unstablemanifoldneartangencyeq1} and \eqref{unstablemanifoldneartangencyeq2}, we see that $\eta_i$ are related to $\mu_i$ by
\begin{equation}\label{globalmapsplittingeq3}
\begin{aligned} 
\eta_0&=\mu_0, &\eta_1&=\mu_1, &\eta_i&=\mu_i+\mathcal F_i(\mu_1,\ldots,\mu_{i-1}) \quad \text{for} \quad |i|\in\{2,\ldots,n\},
\end{aligned}
\end{equation}
with smooth functions $\mathcal F$. The triangular form of \eqref{globalmapsplittingeq3} implies that it can be inverted as 
\begin{equation}\label{globalmapsplittingeq4}
\begin{aligned} 
\mu_0&=\eta_0, &\mu_1&=\eta_1, &\mu_i&=\eta_i+\mathcal {\widetilde F}_i(\eta_1,\ldots,\eta_{i-1}) \quad \text{for} \quad |i|\in\{2,\ldots,n\},
\end{aligned}
\end{equation}
with smooth functions $\mathcal {\widetilde F}$. By \eqref{globalmapsplittingeq3} and \eqref{globalmapsplittingeq4}, we get
\begin{equation*}
\begin{aligned} 
\left.\text{rank}\;\frac{\partial\mathbf H_{n}}{\partial\varepsilon}\right |_{\varepsilon=0}=\left.\text{rank}\;\frac{\partial\mathbf M_{n}}{\partial\varepsilon}\right |_{\varepsilon=0}. 
\end{aligned}
\end{equation*}
\end{demo}

\subsection{First-return maps}\label{subsectionfirstreturnmaps}
Let $f:\mathcal M^d\rightarrow \mathcal M^d$ be a $C^r$ diffeomorphism, where $d\ge 2$ and $r=2,\ldots,\infty,\omega$. Let $f$ have a saddle periodic orbit $L$ whose stable and unstable manifolds intersect over an orbit $\Gamma$ of corank-$c$ homoclinic tangency of order $n$. Suppose $f$ is embedded into a finite-parameter $C^r$ family $(f_{\varepsilon})$, where $f_0\equiv f$. As before, we introduce the local map $T_{O}$, see Section \ref{subsectionlocalmap}. Also, we define  the global map $T_{\Gamma}$ acting between small neighborhoods $\Pi^-$ and $\Pi^+$ of some homoclinic points, see Section \ref{subsectionglobalmap}. For all sufficiently large $k$, we consider the \textit{first-return maps} to $\Pi^+$:
\begin{equation}\label{definitionoffirstreturnmap}
\begin{aligned} 
\mathcal T_{k}\equiv T_{\Gamma}\circ T^k_{O},
\end{aligned}
\end{equation}
which are defined on balls $\sigma_k\equiv\Pi^+\cap T_O^{-k}(\Pi^-)$.

Let $L$ satisfy conditions \textbf{\hyperlink{C1.}{C1}},\textbf{\hyperlink{C4.}{C4}}, and $\Gamma$ satisfy conditions \textbf{\hyperlink{C2.}{C2}},\textbf{\hyperlink{C5.}{C5}}. We fix coordinates as in Lemmas \ref{lemmaglobalmap},\ref{lemmaglobalmapsplitting}, and assume that $\Gamma$ splits freely up to order $n$ in these coordinates. For simplicity, we further assume that the number of parameters $\varepsilon$ equals to the number of the splitting coefficients up to order $n$. 
\begin{lemma}\label{lemmafirstreturnmap}
There exist a sequence $\langle k_\alpha\rangle_{\alpha\in\mathbb N}$ of natural numbers with $k_\alpha\xrightarrow{}+\infty$ and a $\langle k_{\alpha}\rangle$-tame sequence $\langle\mathfrak{D}_{\alpha}\rangle_{\alpha\in\mathbb N}$ such that for all $\varepsilon\in\mathfrak{D}_{\alpha}$ there are $C^{r}$-smooth coordinates $(X,U,Y,V)$ on $\sigma_{k_{\alpha}}$ and a $C^{r-\max(2,n)}$-smooth parameter transformation $\varepsilon\mapsto E$ for which the first-return map $\mathcal T_{k_{\alpha}}$ takes the form
\begin{equation}\label{lemmafirstreturnmapeq0}
\begin{aligned} 
\overline Y_1&=\sum\limits_{|i|=0}^{n} E_{i} \cdot Y_1^{i}+\vartheta_{1,k_{\alpha}}(X,U,Y_1,\overline Y_2,\overline V,E),\\
(\overline X,\overline U,Y_2,V)&=\vartheta_{2,k_{\alpha}}(X,U,Y_1,\overline Y_2,\overline V,E).\\
\end{aligned}
\end{equation}
The range of $(X,U,Y,V)$ and $E$ covers a ball centered at zero whose radius tends to infinity as $\alpha\rightarrow+\infty$. The functions $\vartheta_{k_{\alpha}}$, along with their derivatives with respect to $(X,U,Y_1,\overline Y_2,\overline V)$ up to order $r-2$ and with respect to $E$ up to order $r-\max(2,n)$, tend to zero uniformly on any compact set as $\alpha\rightarrow+\infty$.
\end{lemma}

\begin{demo}
By Lemmas \ref{lemmaglobalmap} and \ref{lemmaglobalmaporder}, for all small $\varepsilon$, the local map $T_{O}$ is in the main normal form and the global map $T_{\Gamma}$ is given by
\begin{equation}\label{firstreturnmapeq1}
\begin{aligned} 
&\overline y_1 =\sum\limits_{|i|=0}^{n} \mu_{i} \cdot (\tilde y_1-\tilde y_1^-)^{i}+O(||\tilde y_1-\tilde y_1^-||^{n+1})+\theta_{31}(\tilde x,\tilde u,\tilde y_1-\tilde y_1^-,\overline y_2,\overline v,\varepsilon),\\
&(\overline x -x^{+},\overline u-u^+,\tilde y_2-\tilde y_2^-,\tilde v-\tilde v^-)=(\theta_{1},\theta_{2},\theta_{32},\theta_{4})(\tilde x,\tilde u,\tilde y_1-\tilde y_1^-,\overline y_2,\overline v,\varepsilon),
\end{aligned}
\end{equation}
where the functions $\theta$ satisfy \eqref{lemmaglobalmapeq4} and \eqref{lemmaglobalmapeq5}. Since the tangency $\Gamma$ splits freely up to order $n$, Lemma \ref{lemmaglobalmapsplitting} implies that the coefficients $\mu_0,\ldots,\mu_n$ satisfy condition \eqref{lemmaglobalmapsplittingeq1}. Therefore, we can take these coefficients as the new independent parameters, and assume $\varepsilon=(\mu_0,\ldots,\mu_n)$ without loss of generality. Note that we lose two degrees of parameter smoothness when bringing the local map $T_O$ to the main normal form (see Section \ref{subsectionlocalmap}), and $n$ degrees of parameter smoothness when defining the coefficients $\mu_0,\ldots,\mu_n$. So, the total loss of the smoothness with respect to the parameters is $\max(2,n)$.

\textit{Step 1.} By Corollary \ref{corollarylocalmap}, there exists a sequence $\langle k_\alpha\rangle_{\alpha\in\mathbb N}$ such that $k_\alpha\xrightarrow{}+\infty$ and for any $\langle k_{\alpha}\rangle$-tame sequence $\langle\mathfrak{D}_{\alpha}\rangle_{\alpha\in\mathbb N}$ and all $\varepsilon\in\mathfrak{D}_{\alpha}$, the iterations $T_O^{k_{\alpha}}$ of the local map are given by \eqref{eqlocal3tamesequence}. Combining formulas \eqref{eqlocal3tamesequence} and \eqref{firstreturnmapeq1}, we write the first-return map $\mathcal T_{k_{\alpha}}$ for such $\varepsilon$ in the form
\begin{equation}\label{firstreturnmapeq2}
\begin{aligned} 
\gamma^{-k_{\alpha}}\cdot(\overline Y_1+\tilde y_1^-)&=\sum\limits_{|i|=0}^{n} \mu_{i} \cdot Y_1^{i}+\Theta_{1}(X,U,Y,V, \overline X,\overline U,\overline Y,\overline V,\varepsilon),\\
(\overline X,\overline U,Y_2,V)&=\Theta_{2}(X,U,Y,V, \overline X,\overline U,\overline Y,\overline V,\varepsilon),\\
\end{aligned}
\end{equation}
where
\begin{equation*}
\begin{aligned} 
&X=x-x^+, \qquad U=u-u^+, \qquad Y=\tilde y-\tilde y^-, \qquad V=\tilde v-\tilde v^-,
\end{aligned}
\end{equation*}
and $\Theta$ are some functions satisfying 
\begin{equation*}
\begin{aligned} 
\Theta_1&=O(||Y_1||^{n+1})+o(\gamma^{-k_{\alpha}}), &\Theta_2&=O(Y_1)+O(\gamma^{-k_{\alpha}}).
\end{aligned}
\end{equation*}

\textit{Step 2.} We nullify the constant term in $\Theta_2$ by making the shift
\begin{equation}\label{firstreturnmapeq3}
\begin{aligned} 
(X,U,Y_2,V)&=(X,U,Y_2,V)^{new}+O(\gamma^{-k_{\alpha}}).
\end{aligned}
\end{equation}
As a result, we rewrite \eqref{firstreturnmapeq2} as
\begin{equation}\label{firstreturnmapeq4}
\begin{aligned} 
\gamma^{-k_{\alpha}}\cdot\overline Y_1&=\sum\limits_{|i|=0}^{n} \mu^*_{i} \cdot Y_1^{i}+\Theta^*_{1}(X,U,Y,V, \overline X,\overline U,\overline Y,\overline V,\varepsilon^*),\\
(\overline X,\overline U,Y_2,V)&=\Theta^*_{2}(X,U,Y,V, \overline X,\overline U,\overline Y,\overline V,\varepsilon^*),\\
\end{aligned}
\end{equation}
where
\begin{equation}\label{firstreturnmapeq4.5}
\begin{aligned} 
\mu^*_0&=\mu_0-\gamma^{-k_{\alpha}}\cdot\tilde y_1^-+o(\gamma^{-k_{\alpha}}),  \quad &\mu^*_i&=\mu_i \quad \text{for} \quad |i|\in\{1,\ldots,n\}, 
\end{aligned}
\end{equation}
and the functions $\Theta^*$ satisfy 
\begin{equation}\label{firstreturnmapeq4.75}
\begin{aligned} 
\Theta^*_1&=O(||Y_1||^{n+1})+o(\gamma^{-k_{\alpha}}), &\Theta^*_2&=O(Y_1)+\gamma^{-k_{\alpha}}\cdot O(X,U,Y,V).
\end{aligned}
\end{equation}

\textit{Step 3.} We rescale the coordinates as
\begin{equation}\label{firstreturnmapeq5}
\begin{aligned} 
Y_1&=\varrho_{k_{\alpha}}\cdot\gamma^{-\frac{k_{\alpha}}{n}}\cdot Y_1^{new}, \quad\quad (X,U,Y_2,V)=\gamma^{-\frac{k_{\alpha}}{n}}\cdot(X,U,Y_2,V)^{new},
\end{aligned}
\end{equation}
where $\varrho_{k_{\alpha}}$ tends sufficiently slowly to zero as $\alpha\rightarrow+\infty$. Moreover, we rescale the parameters as
\begin{equation}\label{firstreturnmapeq6}
\begin{aligned} 
\mu^*_i&=\frac{1}{\varrho^{|i|-1}_{k_{\alpha}}}\cdot\gamma^{-k_{\alpha}\cdot\left(1-\frac{|i|-1}{n}\right)}\cdot E_i \quad \text{for} \quad |i|\in\{0,\ldots,n\}.
\end{aligned}
\end{equation}
By \eqref{firstreturnmapeq4.5} and \eqref{firstreturnmapeq6}, the parameters $\varepsilon=(\mu_0,\ldots,\mu_n)$ approach zero exponentially fast, so there exists a $\langle k_{\alpha}\rangle$-tame sequence $\langle\mathfrak{D}_{\alpha}\rangle_{\alpha\in\mathbb N}$ such that \eqref{eqlocal3tamesequence} and, hence, \eqref{firstreturnmapeq4} holds for any bounded set of $E$. After the rescaling, the first-return map $\mathcal T_{k_{\alpha}}$ takes the form
\begin{equation}\label{firstreturnmapeq7}
\begin{aligned} 
\overline Y_1&=\sum\limits_{|i|=0}^{n} E_{i} \cdot Y_1^{i}+\vartheta^*_{1,k_{\alpha}}(X,U,Y,V, \overline X,\overline U,\overline Y,\overline V,E),\\
(\overline X,\overline U,Y_2,V)&=\vartheta^*_{2,k_{\alpha}}(X,U,Y,V, \overline X,\overline U,\overline Y,\overline V,E),\\
\end{aligned}
\end{equation}
where the functions $\vartheta^*_{k_{\alpha}}$ and their derivatives tend to zero as $\alpha\to+\infty$. 

\textit{Step 4.} Expressing $\overline X,\overline U,\overline Y_1,Y_2,V$ in \eqref{firstreturnmapeq7} as functions of $X,U,Y_1,\overline Y_2,\overline V,E$, we reduce \eqref{firstreturnmapeq7} to the required form \eqref{lemmafirstreturnmapeq0}.
\end{demo}

Lemma \ref{lemmafirstreturnmap} implies 
\begin{lemma}\label{lemmauniversalmaps}
Let $\mathcal Q$ be an open region in the space of $C^r$ diffeomorphisms of a $d$-dimensional manifold $\mathcal M^d$, where $d\ge 2$ and $r=1,\ldots,\infty,\omega$. Suppose there exists a $C^{r}$-dense subset $\widehat{\mathcal Q}\subset\mathcal Q$ such that every diffeomorphism $\hat f\in\widehat{\mathcal Q}$ has orbits of corank-$c$ homoclinic tangency of every order. Then a $C^{r}$-generic diffeomorphism $f\in\mathcal Q$ is cross $(c,r)$-universal.
\end{lemma}

\begin{demo}
Let $\mathbb D^d$ be the closed unit ball in $\mathbb R^d$. We take any $\delta>0$ and any $C^r$ map $\beta:\mathbb D^d\to\mathbb R^d$ of the form
\begin{equation}\label{statementuniversaldynamicseq1}
\begin{aligned} 
(\overline X_1,\ldots,\overline X_{c})&=\Phi(X_{1},\ldots,X_c),\\
(\overline X_{c+1},\ldots,\overline X_q,X_{q+1},\ldots,X_d)&=0,
\end{aligned}
\end{equation}
where $\Phi: \mathbb D^c\to\mathbb R^c$ is an arbitrary $C^r$ map. Let $\mathcal H(\beta,\delta)$ denote the set of all $C^r$ diffeomorphisms of $\mathcal M^d$ which possess a $q$-crossed renormalized iteration that lies within $C^{r}$ distance less than $\delta$ of $\beta$. By definition, $\mathcal H(\beta,\delta)$ is open. We shall show that this set is $C^{r}$-dense in $\mathcal Q$.

Consider any $f\in\mathcal Q$. Since $\widehat{\mathcal Q}$ is $C^{r}$-dense in $\mathcal Q$, there exists a diffeomorphism $\hat f\in\widehat{\mathcal Q}$ that is $C^{r}$-close to $f$. By assumption, $\hat f$ has an orbit $\Gamma$ of corank-$c$ homoclinic tangency of order $n$, where $n$ can be chosen to be as large as we want. Employing Lemma \ref{lemmasplittingoftangency}, we embed $\hat f$ into a finite-parameter family $(\hat f_{\varepsilon})$ $($with $\hat f_0\equiv \hat f)$ such that the tangency $\Gamma$ splits freely up to order $n$. According to Lemma \ref{lemmafirstreturnmap}, there exist rescaled coordinates such that for all sufficiently small $\varepsilon$, the first-return map can be approximated by the polynomial expression
\begin{equation}\label{statementuniversaldynamicseq2}
\begin{aligned} 
\overline Y_1&=\sum\limits_{|i|=0}^{n} E_{i} \cdot Y_1^{i}, &(\overline X,\overline U,Y_2,V)&=0.
\end{aligned}
\end{equation}
Given that $n$ can be arbitrarily large, we may choose $E_i$ such that the polynomial expression \eqref{statementuniversaldynamicseq2} approximates $\beta$ arbitrarily well. Thus, $C^{r}$-close to $f$, we found a diffeomorphism having a $q$-crossed renormalized iteration that lies within $C^{r}$ distance less than $\delta$ of $\beta$. This implies the $C^r$ density of $\mathcal H(\beta,\delta)$ in $\mathcal Q$. 

Let:
\begin{enumerate}
\item $\langle\beta_i\rangle_{i\in\mathbb N}$ be a sequence of maps that are dense in the space of $C^r$ maps of the form \eqref{statementuniversaldynamicseq1};

\item $\langle\delta_j\rangle_{j\in\mathbb N}$ be a sequence of positive real numbers converging to zero.

\end{enumerate}
The intersection $\bigcap\limits_{i,j}\mathcal H(\beta_i,\delta_j)$ is a residual set in $\mathcal Q$ and consists of cross $(c,r)$-universal diffeomorphisms.
\end{demo}

\vspace{2mm}

\begin{demo6}
Theorem \ref{theorem6} follows from Theorem \ref{theorem4} and Lemma \ref{lemmauniversalmaps} immediately. In the smooth case, just note that the cross-universality leads to universality by Lemmas \ref{lemmauniversalityreductiontodisk} and \ref{lemmauniversalityinftysmooth}. 
\end{demo6}

\subsection{Dynamics near homoclinic disks}\label{subsectiondynamicsnearhomoclinicdisks}

Let $f:\mathcal M^d\rightarrow \mathcal M^d$ be a $C^r$ diffeomorphism, where $d\ge 2$ and $r=2,.\ldots,\infty$. Let $f$ have a saddle periodic orbit $L$ whose stable and unstable manifolds intersect over an orbit $\Gamma$ of corank-$c$ homoclinic tangency.

We select a point $O\in L$ and fix a small neighborhood $U(O)$, where the local map $T_{O,f}$ is well-defined. We define the global map $T_{\Gamma,f}:\Pi^-\to\Pi^+$, where $\Pi^+$ and $\Pi^-$ are small neighborhoods of the homoclinic points $M^{+},M^{-}\in\Gamma$ such that $M^{+}\in W^s_{loc}(O)$ and $M^{-}\in W^u_{loc}(O)$. 

In this section, we will focus on the case where $\Gamma$ belongs to a continuous family of homoclinic orbits of corank-$c$ tangency, i.e., the homoclinic point $M^+$ lies in a $c$-dimensional homoclinic disk. But, before that, we prove Lemmas \ref{lemmaperturbationforglobalmap} and \ref{lemmaproducinghomoclinicdisk} without this assumption.

Let $L$ satisfy conditions \textbf{\hyperlink{C1.}{C1}},\textbf{\hyperlink{C4.}{C4}}, and $\Gamma$ satisfies conditions \textbf{\hyperlink{C2.}{C2}},\textbf{\hyperlink{C5.}{C5}}. Similar to Lemma \ref{lemmasplittingoftangency}, we have the following result that allows us to add arbitrary $C^r$-small perturbations to the global map without perturbing the local map.
\begin{lemma}\label{lemmaperturbationforglobalmap}
Let $\langle\omega_{\alpha}\rangle_{\alpha\in\mathbb N}$ be any sequence of $C^r$ functions from $\mathbb R^d$ to $\mathbb R^d$ such that $\omega_{\alpha}\to 0$ in $C^r$ as $\alpha\to+\infty$. There exists a sequence $\langle g_{\alpha}\rangle_{\alpha\in\mathbb N}$ of $C^r$ diffeomorphisms such that 
\begin{equation}\label{lemmaperturbationforglobalmapeq1}
\begin{aligned} 
&g_{\alpha}\to f \;\, \text{in} \;\, C^r \; \text{as} \;\, \alpha\to+\infty, \quad\quad T_{O,g_{\alpha}}\equiv T_{O,f}, \quad\quad T_{\Gamma, g_{\alpha}}=T_{\Gamma, f}+\omega_{\alpha}.
\end{aligned} 
\end{equation} 
\end{lemma}

\begin{demo}
We choose sufficiently small $\rho>0$ and denote by $B^{\rho}_{M^+}$ a ball of radius $\rho$ centered at $M^+$. Making use of a bump function, we define a sequence $\langle\mathcal F_{\alpha}\rangle_{\alpha\in\mathbb N}$ of $C^{r}$ diffeomorphisms such that $\mathcal F_{\alpha}$ is an identity map outside the ball $B^{2\rho}_{M^+}$ and inside the ball $B^{\rho}_{M^+}$ acts as
\begin{equation*}
\begin{aligned} 
&z\mapsto z+\omega_{\alpha}.
\end{aligned} 
\end{equation*} 
Consider a sequence $\langle g_{\alpha}\rangle_{\alpha\in\mathbb N}$ of $C^r$ diffeomorphisms such that
\begin{equation*}
\begin{aligned} 
&g_{\alpha}\equiv\mathcal F_{\alpha}\circ f.
\end{aligned} 
\end{equation*} 
Clearly, $g_{\alpha}\to f$ in $C^r$ as $\alpha\to+\infty$. Considering that $\rho$ is sufficiently small, we obtain that the local map does not change, i.e., $T_{O,g_{\alpha}}\equiv T_{O,f}$. Since $g_{\alpha}$ coincides with $f$ outside a small neighborhood of $f^{-1}(M^+)$ (for all $\alpha$), the global map $T_{\Gamma, g_{\alpha}}$ acting from $\Pi^-$ to $\Pi^+$ is given by 
\begin{equation*}
\begin{aligned} 
&T_{\Gamma, g_{\alpha}}=\mathcal F_{\alpha}\circ T_{\Gamma, f}=T_{\Gamma, f}+\omega_{\alpha}.
\end{aligned} 
\end{equation*} 
\end{demo}

Lemma \ref{lemmaperturbationforglobalmap} entails the following 
\begin{lemma}\label{lemmaproducinghomoclinicdisk}
Let the tangency $\Gamma$ be of order $n$. Then, there exists a sequence $\langle g_{\alpha}\rangle_{\alpha\in\mathbb N}$ of $C^r$ diffeomorphisms such that $g_{\alpha}\to f$ in $C^n$ as $\alpha\to+\infty$, and $g_{\alpha}$ has a $c$-dimensional homoclinic disk $\mathbb B\owns M^+$ in the intersection of $W^s(L)$ and $W^u(L)$.
\end{lemma}

\begin{demo}
By Lemmas \ref{lemmaglobalmap} and \ref{lemmaglobalmaporder}, the global map $T_{\Gamma,f}$ is given by
\begin{equation*}
\begin{aligned} 
\overline y_1&=F(\tilde y_1-\tilde y_1^-)+\theta_{31}(\tilde x,\tilde u,\tilde y_1-\tilde y_1^-,\overline y_2,\overline v,\varepsilon),\\
(\overline x -x^{+},\overline u-u^+,\tilde y_2-\tilde y_2^-,\tilde v-\tilde v^-)&=(\theta_{1},\theta_{2},\theta_{32},\theta_{4})(\tilde x,\tilde u,\tilde y_1-\tilde y_1^-,\overline y_2,\overline v,\varepsilon),
\end{aligned}
\end{equation*}
where the function $F$ satisfies \eqref{lemmaglobalmapordereq1} and the functions $\theta$ satisfy \eqref{lemmaglobalmapeq4},\eqref{lemmaglobalmapeq5}. 

By Lemma \ref{lemmaperturbationforglobalmap}, there exists a sequence $\langle g_{\alpha}\rangle_{\alpha\in\mathbb N}$ of $C^r$ diffeomorphisms such that \eqref{lemmaperturbationforglobalmapeq1} holds with $C^r$ functions $\langle\omega_{\alpha}\rangle_{\alpha\in\mathbb N}$ such that 
\begin{equation}\label{lemmaproducinghomoclinicdiskeq2}
\begin{aligned} 
\omega_{31,\alpha}(\tilde x,\tilde u,\tilde y_1-\tilde y_1^-,\overline y_2,\overline v)&=-F(\tilde y_1-\tilde y_1^-),\\
(\omega_{1,\alpha},\omega_{2,\alpha},\omega_{32,\alpha},\omega_{4,\alpha})(\tilde x,\tilde u,\tilde y_1-\tilde y_1^-,\overline y_2,\overline v)&=0.
\end{aligned}
\end{equation}
The tangency $\Gamma$ is of order $n$, so the $C^n$ norm of the map \eqref{lemmaproducinghomoclinicdiskeq2} vanishes at the point $(0,0,\tilde y_1=\tilde y_1^-,0,0)$. By choosing a sequence $\langle\rho_{\alpha}\rangle_{\alpha\in\mathbb N}$ such that $\rho_{\alpha}\to 0$ as $\alpha\to+\infty$ (see proof of Lemma \ref{lemmaperturbationforglobalmap}), we make the functions $\omega_{\alpha}$, along with their derivatives up to order $n$, decay to zero as $\alpha\to+\infty$. We get that the global map $T_{\Gamma,g_{\alpha}}$ is given by  
\begin{equation*}
\begin{aligned} 
\overline y_1&=\theta_{31}(\tilde x,\tilde u,\tilde y_1-\tilde y_1^-,\overline y_2,\overline v),\\
(\overline x -x^{+},\overline u-u^+,\tilde y_2-\tilde y_2^-,\tilde v-\tilde v^-)&=(\theta_{1},\theta_{2},\theta_{32},\theta_{4})(\tilde x,\tilde u,\tilde y_1-\tilde y_1^-,\overline y_2,\overline v),
\end{aligned}
\end{equation*}
which proves the lemma.
\end{demo}

Now we assume that $M^+\in\Gamma$ lies inside a $c$-dimensional homoclinic disk $\mathbb B$. The following result, similar to Lemma \ref{lemmafirstreturnmap}, describes the first-return maps $\mathcal T_{k}$ on small balls $\sigma_k\equiv\Pi^+\cap T_O^{-k}(\Pi^-)$, see \eqref{definitionoffirstreturnmap}. Let $\langle k_\alpha\rangle_{\alpha\in\mathbb N}$ be a sequence of natural numbers given by Corollary \ref{corollarylocalmap}, so $k_\alpha\to+\infty$ and the equations \eqref{eqlocal3tamesequence} are fulfilled.

\begin{lemma}\label{lemmafirstreturnmaphomoclinicdisk}
For every $\alpha\in\mathbb N$ there exist functions $\vartheta_{1,k_{\alpha}}:\mathbb R^{2d}\to\mathbb R^{c}$ and $\vartheta_{2,k_{\alpha}}:\mathbb R^{2d}\to\mathbb R^{d-c}$, which tend to zero in $C^{r-2}$ on any compact set as $\alpha\rightarrow+\infty$, such that the following holds. Given any sequence of $C^r$ functions $\Phi_{\alpha}: \mathbb R^c\to\mathbb R^c$, uniformly bounded in $C^r$ for all $\alpha$, and any sequence of $C^r$ functions $\Omega_{\alpha}: \mathbb R^c\to\mathbb R^{d-c}$ such that $\Omega_{\alpha}\to 0$ in $C^r$ as $\alpha\to+\infty$, for every $\alpha\in\mathbb N$ there exist a $C^r$ diffeomorphism $g_{\alpha}$ $(g_{\alpha}\to f$ in $C^r$ as $\alpha\to+\infty)$ and $C^{r}$-smooth coordinates $\left(X,U,Y,V\right)$ on $\sigma_{k_{\alpha}}$ for which the first-return map $\mathcal T_{k_{\alpha},g_{\alpha}}$ takes the form
\begin{equation}\label{lemmafirstreturnmaphomoclinicdiskeq0}
\begin{aligned} 
\overline Y_1&=\Phi_{\alpha}(Y_1)+\vartheta_{1,k_{\alpha}}(Y_1,X,U,\overline Y_2,\overline V,\overline Y_1,\overline X,\overline U,Y_2,V),\\
(\overline X,\overline U,Y_2,V)&=\Omega_{\alpha}(Y_1)+\vartheta_{2,k_{\alpha}}(Y_1,X,U,\overline Y_2,\overline V,\overline Y_1,\overline X,\overline U,Y_2,V).\\
\end{aligned}
\end{equation}
The range of $\left(X,U,Y,V\right)$ covers a centered at zero ball whose radius tends to infinity as $\alpha\to+\infty$. 
\end{lemma}

\begin{demo}
Applying Lemma \ref{lemmaglobalmap}, we choose coordinates such that the local map $T_{O,f}$ is in the main normal form and the global map $T_{\Gamma,f}$ has the form
\begin{equation*}
\begin{aligned} 
\overline y_1&=\theta_{31}(\tilde x,\tilde u,\tilde y_1-\tilde y_1^-,\overline y_2,\overline v),\\
(\overline x-x^{+},\overline u-u^+,\tilde y_2-\tilde y_2^-,\tilde v-\tilde v^-)&=(\theta_{1},\theta_{2},\theta_{32},\theta_{4})(\tilde x,\tilde u,\tilde y_1-\tilde y_1^-,\overline y_2,\overline v),
\end{aligned}
\end{equation*}
where the functions $\theta$ satisfy \eqref{lemmaglobalmapeq4} and \eqref{lemmaglobalmapeq5}. Note that the function $F$ of \eqref{lemmaglobalmapeq1} vanishes, since $M^+$ lies inside the homoclinic disk $\mathbb B$. 

Let $\langle\omega_{\alpha}\rangle_{\alpha\in\mathbb N}$ be a sequence of $C^r$ functions such that 
\begin{equation*}
\begin{aligned} 
\omega_{31,\alpha}(\tilde x,\tilde u,\tilde y_1-\tilde y_1^-,\overline y_2,\overline v)&=\gamma^{-2k_{\alpha}}\cdot(\Phi_{\alpha}(\gamma^{k_{\alpha}}\cdot(\tilde y_1-\tilde y_1^-))+\tilde y^-),\\
(\omega_{1,\alpha},\omega_{2,\alpha},\omega_{32,\alpha},\omega_{4,\alpha})(\tilde x,\tilde u,\tilde y_1-\tilde y_1^-,\overline y_2,\overline v)&=\frac{\gamma^{-k_{\alpha}}}{\varrho_{k_{\alpha}}}\cdot\Omega_{\alpha}(\gamma^{k_{\alpha}}\cdot(\tilde y_1-\tilde y_1^-)),
\end{aligned}
\end{equation*}
where $\varrho_{k_{\alpha}}$ tends sufficiently slowly to zero as $\alpha\rightarrow+\infty$. Applying Lemma \ref{lemmaperturbationforglobalmap}, we construct a sequence $\langle g_{\alpha}\rangle_{\alpha\in\mathbb N}$ of $C^r$ diffeomorphisms such that the relations \eqref{lemmaperturbationforglobalmapeq1} hold. Thus, the global map $T_{\Gamma,g_{\alpha}}$ is given by  
\begin{equation}\label{lemmafirstreturnmaphomoclinicdiskeq3}
\begin{aligned} 
\overline y_1&=\gamma^{-2k_{\alpha}}\cdot(\Phi_{\alpha}(\gamma^{k_{\alpha}}\cdot(\tilde y_1-\tilde y_1^-))+\tilde y^-)+\theta_{31}(\tilde x,\tilde u,\tilde y_1-\tilde y_1^-,\overline y_2,\overline v),\\
(\overline x -x^{+},\overline u-u^+,\tilde y_2-\tilde y_2^-,\tilde v-\tilde v^-)&=\frac{\gamma^{-k_{\alpha}}}{\varrho_{k_{\alpha}}}\cdot\Omega_{\alpha}(\gamma^{k_{\alpha}}\cdot(\tilde y_1-\tilde y_1^-))+(\theta_{1},\theta_{2},\theta_{32},\theta_{4})(\tilde x,\tilde u,\tilde y_1-\tilde y_1^-,\overline y_2,\overline v).
\end{aligned}
\end{equation}

Substituting formulas \eqref{eqlocal3tamesequence} for the local map $T_{O,g_{\alpha}}$ into \eqref{lemmafirstreturnmaphomoclinicdiskeq3}, we write the first-return maps $\mathcal T_{k_{\alpha},g_{\alpha}}$ as
\begin{equation}\label{lemmafirstreturnmaphomoclinicdiskeq4}
\begin{aligned} 
\gamma^{-k_{\alpha}}\cdot(\overline Y_1+\tilde y_1^-)&=\gamma^{-2k_{\alpha}}\cdot(\Phi_{\alpha}(\gamma^{k_{\alpha}}\cdot Y_1)+\tilde y^-)+\Theta_{1}(Y_1,X,U,\overline Y_2,\overline V,\overline Y_1,\overline X,\overline U,Y_2,V),\\
(\overline X,\overline U,Y_2,V)&=\frac{\gamma^{-k_{\alpha}}}{\varrho_{k_{\alpha}}}\cdot\Omega_{\alpha}(\gamma^{k_{\alpha}}\cdot Y_1)+\Theta_{2}(Y_1,X,U,\overline Y_2,\overline V,\overline Y_1,\overline X,\overline U,Y_2,V),\\
\end{aligned}
\end{equation}
where $X,U,Y,V$ are given by \eqref{lemma1highordertangencieseq2.5}, and $\Theta$ are functions satisfying 
\begin{equation*}
\begin{aligned} 
\Theta_1&=o(\gamma^{-k_{\alpha}}), &\Theta_2&=O(Y_1)+O(\gamma^{-k_{\alpha}}).
\end{aligned}
\end{equation*}
We nullify the constant term in $\Theta_2$ by making the shift of the form \eqref{firstreturnmapeq3}, and rescale the coordinates as
\begin{equation*}
\begin{aligned} 
Y_1&=\gamma^{-k_{\alpha}}\cdot Y_1^{new}, \quad\quad (X,U,Y_2,V)=\frac{\gamma^{-k_{\alpha}}}{\varrho_{k_{\alpha}}}\cdot(X,U,Y_2,V)^{new}.
\end{aligned}
\end{equation*}
As a result, the system \eqref{lemmafirstreturnmaphomoclinicdiskeq4} takes the form \eqref{lemmafirstreturnmaphomoclinicdiskeq0}.
\end{demo}

Before proving the main result of this Section (Lemma \ref{lemmachoosingsequences}) , we prove auxiliary
\begin{lemma}\label{lemmainvariantdisk}
Let $\phi: \mathbb R^{c}\to\mathbb R^{c}$ be any diffeomorphism such that $\phi$ and $\phi^{-1}$ are uniformly bounded in $C^r$ for some finite $r$. Let $\hat\vartheta:\mathbb R^{2d}\to\mathbb R^{d-c}$ be $C^r$-small. Let $\widehat {\mathcal T}$ be a (possibly multi-valued) map between subsets of $\mathbb R^{d}$ defined by 
\begin{equation}\label{lemmainvariantdiskeq1}
\begin{aligned} 
\overline Y_1&=\phi(Y_1),\\
(\overline X,\overline U,Y_2,V)&=\hat\vartheta(Y_1,X,U,\overline Y_2,\overline V,\overline Y_1,\overline X,\overline U,Y_2,V),
\end{aligned}
\end{equation} 
meaning that the point $\left(\overline X,\overline U,\overline Y,\overline V\right)$ is an image of $(X,U,Y,V)$ if and only if the relation \eqref{lemmainvariantdiskeq1} holds. Then $\widehat {\mathcal T}$ has a normally-hyperbolic invariant $c$-dimensional $C^r$-smooth manifold $\widehat {\mathbb M}^c$ given by 
\begin{equation}\label{lemmainvariantdiskeq2}
\begin{aligned} 
(X,U,Y_2,V)&=\hat h(Y_1),
\end{aligned}
\end{equation}
where the function $\hat h$ is $C^r$-small.
\end{lemma}

\begin{demo}
By the second part of Proposition \ref{propositionnormallycontractingexpanding} (where we take $u=(Y_1,X,U)$ and $v=(Y_2,V)$), the map $\widehat {\mathcal T}$ has a normally-hyperbolic repelling invariant $c$-dimensional $C^r$-smooth manifold $\widetilde {\mathbb M}^c$ given by 
\begin{equation*}
\begin{aligned} 
(Y_2,V)&=\tilde h(Y_1,X,U),
\end{aligned}
\end{equation*}
where the function $\tilde h$ is $C^r$-small if $\hat \vartheta$ is $C^r$-small. Applying the first part of Proposition  \ref{propositionnormallycontractingexpanding} to the restriction $\widehat {\mathcal T}_{\widetilde {\mathbb M}^c}$ (we take $u=(X,U)$ and $v=Y_1$), we obtain that $\widehat {\mathcal T}_{\widetilde {\mathbb M}^c}$ has a normally-hyperbolic attracting invariant $c$-dimensional $C^r$-smooth manifold $\widehat {\mathbb M}^c$ given by \eqref{lemmainvariantdiskeq2}. Clearly, this manifold is normally-hyperbolic and invariant for the map $\widehat {\mathcal T}$.
\end{demo}

\begin{remark}
The invariance of the manifold $\widehat {\mathbb M}^c$ under $\widehat {\mathcal T}$ means that for each point $\left(Y_1,X,U,\overline Y_2,\overline V\right)\in\widehat {\mathbb M}^c$, its image $(\overline Y_1,\overline X,\overline U,Y_2,V)\in\widehat {\mathbb M}^c$, i.e.,
\begin{equation*}
\begin{aligned} 
(\overline X,\overline U,Y_2,V)&=\hat h(\overline Y_1),
\end{aligned}
\end{equation*}
which along with \eqref{lemmainvariantdiskeq1} gives
\begin{equation}\label{corollaryinvariantdiskeq1}
\begin{aligned} 
\hat h(\phi(Y_1))&=\hat\vartheta(Y_1,\hat h(Y_1),\phi(Y_1),\hat h(\phi(Y_1)).
\end{aligned}
\end{equation}
\end{remark}

\begin{lemma}\label{lemmachoosingsequences}
Let $\phi: \mathbb R^{c}\to\mathbb R^{c}$ be any $C^r$ diffeomorphism. The sequences $\langle\Phi_{\alpha}\rangle_{\alpha\in\mathbb N},\langle\Omega_{\alpha}\rangle_{\alpha\in\mathbb N}$ and $\langle g_{\alpha}\rangle_{\alpha\in\mathbb N}$ in Lemma \ref{lemmafirstreturnmaphomoclinicdisk} can be chosen such that the first-return map $\mathcal T_{k_{\alpha},g_{\alpha}}$ of \eqref{lemmafirstreturnmaphomoclinicdiskeq0} has, for every $\alpha$, a normally-hyperbolic invariant $c$-dimensional $C^{\infty}$-smooth disk $D^c_\alpha$ given by
\begin{equation*}
\begin{aligned} 
(X,U,Y_2,V)&=h_{\alpha}(Y_1),
\end{aligned}
\end{equation*}
and the restriction of $\mathcal T_{k_{\alpha},g_{\alpha}}$ to $D^c_{\alpha}$ is
\begin{equation}\label{lemmachoosingsequenceseq2}
\begin{aligned} 
\overline Y_1&=\phi(Y_1).
\end{aligned}
\end{equation}
\end{lemma}

\begin{demo}
Let $s=r$ if $r$ is finite, or $s$ be sufficiently large if $r=\infty$. Let $\langle D^c_{\alpha}\rangle_{\alpha\in\mathbb N}$ be a sequence of $c$-dimensional $C^{\infty}$ disks given by
\begin{equation*}
\begin{aligned} 
(X,U,Y_2,V)&=h_{\alpha}(Y_1),
\end{aligned}
\end{equation*}
such that $D^c_{\alpha}\to\widehat {\mathbb M}^c$ in $C^s$ as $\alpha\to+\infty$, where the manifold $\widehat {\mathbb M}^c$ is taken from Lemma \ref{lemmainvariantdisk}. Making use of Lemma \ref{lemmafirstreturnmaphomoclinicdisk}, we choose the sequence $\langle g_{\alpha}\rangle_{\alpha\in\mathbb N}$ such that the first-return maps $\mathcal T_{k_{\alpha},g_{\alpha}}$ are given by \eqref{lemmafirstreturnmaphomoclinicdiskeq0}, where
\begin{equation*}
\begin{aligned} 
\Phi_{\alpha}(Y_1)&\equiv\phi(Y_1)-\vartheta_{1,k_{\alpha}}(Y_1,h_{\alpha}(Y_1),\phi(Y_1),h_{\alpha}(\phi(Y_1))),\\
\Omega_{\alpha}(Y_1)&\equiv h_{\alpha}(\phi(Y_1))-\vartheta_{2,k_{\alpha}}(Y_1,h_{\alpha}(Y_1),\phi(Y_1),h_{\alpha}(\phi(Y_1))).
\end{aligned}
\end{equation*}
Thus, the first-return map $\mathcal T_{k_{\alpha},g_{\alpha}}$ takes the form
\begin{equation*}
\begin{aligned} 
\overline Y_1&=\phi(Y_1)-\vartheta_{1,k_{\alpha}}(Y_1,h_{\alpha}(Y_1),\phi(Y_1),h_{\alpha}(\phi(Y_1)))+\\
&+\vartheta_{1,k_{\alpha}}(Y_1,X,U,\overline Y_2,\overline V,\overline Y_1,\overline X,\overline U,Y_2,V),\\
(\overline X,\overline U,Y_2,V)&=h_{\alpha}(\phi(Y_1))-\vartheta_{2,k_{\alpha}}(Y_1,h_{\alpha}(Y_1),\phi(Y_1),h_{\alpha}(\phi(Y_1)))+\\
&+\vartheta_{2,k_{\alpha}}(Y_1,X,U,\overline Y_2,\overline V,\overline Y_1,\overline X,\overline U,Y_2,V).
\end{aligned}
\end{equation*}
Clearly, the disk $D^c_{\alpha}$ is invariant under $\mathcal T_{k_{\alpha},g_{\alpha}}$ and the formula \eqref{lemmachoosingsequenceseq2} holds. Moreover, by \eqref{corollaryinvariantdiskeq1}, $\Omega_{\alpha}\to 0$ in $C^s$ as $\alpha\to+\infty$ since $h_{\alpha}\to\hat h$ in $C^s$.
\end{demo}

Theorem \ref{theorem7} follows from Theorem \ref{theorem4} and Lemma \ref{lemmainvariantmanifolddensesubset} immediately.

\begin{lemma}\label{lemmainvariantmanifolddensesubset}
Let $\mathcal Q$ be an open region in the space of $C^r$ diffeomorphisms of a $d$-dimensional manifold $\mathcal M^d$, where $d\ge 2$ and $r=1,\ldots,\infty$. Suppose there exists a $C^r$-dense subset $\widehat{\mathcal Q}\subset\mathcal Q$ such that every $\hat f\in\widehat{\mathcal Q}$ has orbits of corank-$c$ homoclinic tangency of every order. Then, for every $C^r$ diffeomorphism $\phi:\mathbb D^c\to \mathbb D^c$, there exists a dense subset $\mathcal A_{\phi}$ of a region $\mathcal Q$ such that for every $\overline f\in\mathcal A_{\phi}$, its $c$-restricted renormalized iteration $\overline f_{k,\psi}$, for some $k>0$ and $\psi\in C^{\infty}$, equals to $\phi$.
\end{lemma}

\begin{demo}
Take any $C^r$ diffeomorphism $\phi:\mathbb D^c\to \mathbb D^c$ and any diffeomorhism $f\in\mathcal Q$. We choose a number $n\in\mathbb N$ as follows: 
\begin{enumerate}

\item if $r$ is finite, then $n=r$; 

\item if $r=\infty$, then $n$ is sufficiently large.
\end{enumerate}
Let $\hat f\in\widehat{\mathcal Q}$ be $C^r$-close to $f$ and let $\Gamma$ be an orbit of corank-$c$ homoclinic tangency of order $n$ of diffeomorphism $\hat f$. Employing Lemma \ref{lemmaproducinghomoclinicdisk}, $C^n$-close to $\hat f$, we find a diffeomorphism $\hat f^{new}$ which has a $c$-dimensional homoclinic disk. Applying Lemma \ref{lemmachoosingsequences}, $C^r$-close to $\hat f$, we find a diffeomorphism $\overline f$ such that its $c$-restricted renormalized iteration $\overline f_{k,\psi}$, for some $k>0$ and $\psi\in C^{\infty}$, equals to $\phi$.
\end{demo}

\section*{Acknowledgments}

We would like to thank M. Asaoka, P.G. Barrientos, P. Berger, A. Khodaeian Karim, D. Li, and S. van Strien for useful discussions. Dmitrii Mints was supported by the EPSRC Centre for Doctoral Training in Mathematics
of Random Systems: Analysis, Modelling and Simulation (EP/S023925/1).

\appendix

\section{Appendix. Cross-universal dynamics}\label{sectionappendix}

\counterwithin{definition}{section} 
\renewcommand{\thedefinition}{\thesection.\arabic{definition}}

\counterwithin{lemma}{section} 
\renewcommand{\thelemma}{\thesection.\arabic{lemma}}

\counterwithin{remark}{section} 
\renewcommand{\theremark}{\thesection.\arabic{remark}}

\counterwithin{proposition}{section} 
\renewcommand{\theproposition}{\thesection.\arabic{proposition}}

\begin{definition}\label{definitioncrossform}
Let $\mathbb U,\mathbb V,\overline {\mathbb U},\overline {\mathbb V}$ be arbitrary sets. Let a map $g$ take a subset of $\mathbb U\times \mathbb V$ to $\overline {\mathbb U}\times \overline {\mathbb V}$. We say that a map $g$ is defined in a \textit{cross form} by the pair of maps $g^{\times}_1:\mathbb U\times \overline {\mathbb V}\to \overline {\mathbb U}$ and $g^{\times}_2:\mathbb U\times \overline {\mathbb V}\to \mathbb V$, if, for every $(u,v)\in \mathbb U\times \mathbb V$ and $(\overline u,\overline v)\in \overline {\mathbb U}\times \overline {\mathbb V}$, the following equivalence holds:
\begin{equation}\label{crossformeq1}
\begin{aligned} 
(\overline u,\overline v)&=g(u,v) \iff \begin{cases} \overline u=g^{\times}_1(u,\overline v),\\ v=g^{\times}_2(u,\overline v).\end{cases}
\end{aligned}
\end{equation}
\end{definition}
As we see, the cross form $g^{\times}=(g^{\times}_1,g^{\times}_2)$ is a single-valued map $\mathbb U\times \overline {\mathbb V}\to\overline {\mathbb U}\times \mathbb V$, and $g$ itself is a map (possibly multi-valued) of $\{(u,v):u\in \mathbb U,v\in g_2^{\times}(u,\overline {\mathbb V})\}$ to $\{(\overline u,\overline v):\overline u\in g_1^{\times}(\mathbb U,\overline v), \overline v\in \overline {\mathbb V}\}$. 
\begin{definition}\label{definitionnormallycontracting}
Let $g$ be a map represented in a cross form \eqref{crossformeq1} for some convex closed subsets $\mathbb U,\overline {\mathbb U},\mathbb V,\overline {\mathbb V}$ of $\mathbb R^d$ and some functions $g^{\times}_1,g^{\times}_2$ which are uniformly bounded in $C^r$. Let $s$ be an integer such that $1\le s\le r$.

We say that $g$ is normally $C^s$-contracting if $\mathbb U=\overline {\mathbb U}$ and $g^{\times}_1,g^{\times}_2$ satisfy
\begin{equation*}
\begin{aligned} 
\sqrt[s+1]{\sup\limits_{(u,\overline v)\in \mathbb U\times \overline {\mathbb V}}\left(\left|\left|\frac{\partial g_1^{\times}}{\partial u}\right|\right|\cdot \left|\left|\frac{\partial g^{\times}_2}{\partial  \overline v}\right|\right|^s\right)}+\sqrt{\sup_{(u,\overline v)\in \mathbb U\times \overline {\mathbb V}}\left|\left|\frac{\partial g_1^{\times}}{\partial  \overline v}\right|\right|\cdot \sup_{(u,\overline v)\in \mathbb U\times \overline {\mathbb V}}\left|\left|\frac{\partial g^{\times}_2}{\partial u}\right|\right|}<1,\\
\sup\limits_{(u,\overline v)\in \mathbb U\times \overline {\mathbb V}}\left|\left|\frac{\partial g_1^{\times}}{\partial u}\right|\right|+\sqrt{\sup\limits_{(u,\overline v)\in \mathbb U\times \overline {\mathbb V}}\left|\left|\frac{\partial g_1^{\times}}{\partial \overline v}\right|\right|\cdot\sup\limits_{(u,\overline v)\in \mathbb U\times \overline {\mathbb V}}\left|\left|\frac{\partial g_2^{\times}}{\partial u}\right|\right|}<1.
\end{aligned}
\end{equation*}
We say that $g$ is normally $C^s$-expanding if $\mathbb V=\overline {\mathbb V}$ and $g^{\times}_1,g^{\times}_2$ satisfy
\begin{equation*}
\begin{aligned} 
\sqrt[s+1]{\sup\limits_{(u,\overline v)\in \mathbb U\times \overline {\mathbb V}}\left(\left|\left|\frac{\partial g_2^{\times}}{\partial \overline v}\right|\right|\cdot \left|\left|\frac{\partial g^{\times}_1}{\partial u}\right|\right|^s\right)}+\sqrt{\sup_{(u,\overline v)\in \mathbb U\times \overline {\mathbb V}}\left|\left|\frac{\partial g_2^{\times}}{\partial u}\right|\right|\cdot \sup_{(u,\overline v)\in \mathbb U\times \overline {\mathbb V}}\left|\left|\frac{\partial g^{\times}_1}{\partial \overline v}\right|\right|}<1,\\
\sup\limits_{(u,\overline v)\in \mathbb U\times \overline {\mathbb V}}\left|\left|\frac{\partial g_2^{\times}}{\partial \overline v}\right|\right|+\sqrt{\sup\limits_{(u,\overline v)\in \mathbb U\times \overline {\mathbb V}}\left|\left|\frac{\partial g_2^{\times}}{\partial u}\right|\right|\cdot\sup\limits_{(u,\overline v)\in \mathbb U\times \overline {\mathbb V}}\left|\left|\frac{\partial g_1^{\times}}{\partial \overline v}\right|\right|}<1.
\end{aligned}
\end{equation*}
\end{definition}
It is directly verified that a normally $C^s$-contracting map is normally $C^i$-contracting for each $i\in\{1,\ldots,s-1\}$. Similarly, a normally $C^s$-expanding map is normally $C^i$-explanding for each $i\in\{1,\ldots,s-1\}$.

\begin{proposition}\label{propositionnormallycontractingexpanding}(\cite{Shilnikov1}, Theorem 4.3 and Theorem 4.4). The following holds: 
\begin{enumerate}
\item When $\mathbb V=\overline {\mathbb V}$, a normally $C^s$-contracting map $g$ has a normally-hyperbolic (attracting) invariant manifold $\mathbb M$, which is a graph of a $C^s$-smooth function $u(v)$.

\item When $\mathbb U=\overline {\mathbb U}$, a normally $C^s$-expanding map $g$ has a normally-hyperbolic (repelling) invariant manifold $\mathbb M$, which is a graph of a $C^s$-smooth function $v(u)$.
\end{enumerate}
\end{proposition}

Note that we do not require the map $g$ to be single-valued, or everywhere defined, or smooth. Therefore, Proposition \ref{propositionnormallycontractingexpanding} needs a clarification. In the normally $C^s$-contracting case, for each point of $\mathbb M$, its full image by $g$ lies in $\mathbb M$, and $g^{-1}|_{\mathbb M}$ has a single-valued $C^s$-smooth branch $b_a$ defined everywhere on $\mathbb M$. In the normally $C^s$-expanding case, for each point of $\mathbb M$, its full image by $g^{-1}$ lies in $\mathbb M$, and $g|_{\mathbb M}$ has a single-valued $C^s$-smooth branch $b$ defined everywhere on $\mathbb M$. In both cases, the manifold $\mathbb M$ and the map $b$ depend on $g^{\times}$ continuously in $C^s$.

Let $f:\mathcal M^d\to\mathcal M^d$ be a $C^r$ diffeomorphism, where $d\ge 2$ and $r=1,\ldots,\infty,\omega$. Take a positive integer $q\le d$, and let $\mathbb D^{q}, \mathbb D^{d-q}$ be closed unit disks in $\mathbb R^{q}, \mathbb R^{d-q}$, respectively. Take any $C^r$ embedding $\psi: \mathbb R^d\hookrightarrow \mathcal M^d$. Suppose that the map $\psi^{-1}\circ f^k\circ \psi$ is normally $C^1$ expanding with $\mathbb U=\mathbb D^q$, $\overline {\mathbb U}=\mathbb R^q$, $\mathbb V=\overline {\mathbb V}=\mathbb D^{d-q}$; in particular, it has a cross form representation 
\begin{equation*}
\begin{aligned} 
f^{\times}_{k,\psi}\equiv(\psi^{-1}\circ f^k\circ \psi)^{\times}.
\end{aligned}
\end{equation*}
We call the map $f^{\times}_{k,\psi}: \mathbb D^q\times \mathbb D^{d-q}\to\mathbb R^q\times \mathbb D^{d-q}$ a \textit{$q$-crossed renormalized iteration} of $f$. The set of all possible $q$-crossed renormalized iterations is called the \textit{$q$-crossed renormalization class} of $f$. 
\begin{definition}\label{definitioncrossuniversaldynamics}
Let $f:\mathcal M^d\rightarrow\mathcal M^d$ be a $C^r$ diffeomorphism and let $1\le p< d$. We say that $f$ is cross $(p,r)$-universal if there exists $q$ such that $p<q\le d$ and the $C^r$-closure of the $q$-crossed renormalization class of $f$ contains all maps $\mathbb D^q\times \mathbb D^{d-q}\to\mathbb R^d$ of the form
\begin{equation}\label{definitioncrossuniversaldynamicseq1}
X\mapsto \overline X: \quad
\begin{aligned} 
(\overline X_1,\ldots,\overline X_{p})&= \Phi(X_{1},\ldots,X_p), \\
(\overline X_{p+1},\ldots,\overline X_q,X_{q+1},\ldots,X_d)&=0,
\end{aligned}
\end{equation}
where $\Phi$ is an arbitrary $C^r$ map of the closed unit ball $\mathbb D^p$ into $\mathbb R^p$. 
\end{definition}

\begin{remark}\label{remarkcrossuniversaldynamics}
In the case $q=d$ $($where there are no uniformly expanding directions in the partially hyperbolic decomposition of $f)$, Definitions \ref{definitionuniversaldynamics} and \ref{definitioncrossuniversaldynamics} coincide. 
\end{remark}

In the case $q\not=d$, Definition \ref{definitioncrossuniversaldynamics} means that for any $C^r$ map $\Phi: \mathbb D^p\to\mathbb R^p$ there exists $k\in\mathbb N$ and a $C^r$ embedding $\psi: \mathbb R^d\hookrightarrow \mathcal M^d$ such that a renormalized iteration $f_{k,\psi}$ can be written in a cross form \eqref{crossformeq1} as
\begin{equation}\label{renormalizediterationcrossform}
\begin{aligned} 
(\overline X_1,\ldots,\overline X_{p})&= \Phi(X_{1},\ldots,X_p)+o(1), &(\overline X_{p+1},\ldots,\overline X_q)&=o(1),\\
(X_{q+1},\ldots,X_d)&=o(1),
\end{aligned}
\end{equation}
where $(X_{1},\ldots,X_q)=u$ and $(X_{q+1},\ldots,X_d)=v$. Obviously, for a finite $r$, this map is $C^r$-expanding. We extend the right-hand side of \eqref{renormalizediterationcrossform} for all $X_1,\ldots,X_q$ (this includes extension of $\Phi$ to a $C^r$ map of $\mathbb R^p$ to $\mathbb R^p$) so that the extended map remains $C^r$-expanding. Then, since $\mathbb U=\overline {\mathbb U}=\mathbb R^q$, Proposition \ref{propositionnormallycontractingexpanding} gives an existence of a repelling invariant  $q$-dimensional $C^r$-smooth manifold $\mathbb M^q$. Define $B^q\equiv\psi\left(\mathbb M^q\cap\{||(X_1,\ldots,X_q)||\le 1\}\right)$. Since $B^q\cup f^k(B^q)$ lies in $\psi(\mathbb M^q)$, this union is a $C^r$ submanifold of $\mathcal M^d$. Thus, $B^q$ is a disk from Definition \ref{definitionuniversaldynamics}. The restriction of $f^k$ to this disk, in the coordinates defined by $\psi$, is a $q$-restricted renormalized iteration of Definition \ref{definitionuniversaldynamics} which is $C^r$-close to
\begin{equation*}
\begin{aligned} 
(\overline X_1,\ldots,\overline X_{p})&= \Phi(X_{1},\ldots,X_p), &(\overline X_{p+1},\ldots,\overline X_q)&=0.\end{aligned}
\end{equation*}
Therefore, we have the following
\begin{lemma}\label{lemmauniversalityreductiontodisk}
For a finite $r$, universality in a sense of Definition \ref{definitioncrossuniversaldynamics} implies universality in a sense of Definition \ref{definitionuniversaldynamics}. 
\end{lemma}

In the $C^{\infty}$ case, similarly to arguments of Lemma \ref{lemmachoosingsequences} (approximating finitely smooth disks $\mathbb M^q$ by $C^{\infty}$ disks and adding small perturbations that make these disks invariant), one proves
\begin{lemma}\label{lemmauniversalityinftysmooth}
In any $C^{\infty}$ neighborhood of a $(p,\infty)$-universal diffeomorphism in a sense of Definition \ref{definitioncrossuniversaldynamics} there exists a $(p,\infty)$-universal diffeomorphism in a sense of Definition \ref{definitionuniversaldynamics}.  
\end{lemma}

In a similar way, one defines the cross-universality for parametric families of diffeomorphisms.

\renewcommand\refname{References}


\begin{thebibliography}{99} 

\bibitem{Asaoka2008} M. Asaoka. Hyperbolic sets exhibiting $C^1$-persistent homoclinic tangency for higher dimensions. \textit{Proceedings of the American Mathematical Society}. \textbf{136}(2) (2008), 677--686.

\bibitem{Asaoka} M. Asaoka. Stable intersection of Cantor sets in higher dimension and robust homoclinic tangency of the largest codimension. \textit{Trans. Amer. Math. Soc.} \textbf{375}(2) (2022), 873--908.

\bibitem{Asaokadegeneratebehavior} M. Asaoka, K. Shinohara and D. Turaev. Degenerate behavior in non-hyperbolic semigroup actions on the interval: fast growth of periodic points and universal dynamics. \textit{Math. Ann.} \textbf{368} (2017), 1277--1309.

\bibitem{Asaokafastgrowth} M. Asaoka, K. Shinohara and D. Turaev. Fast growth of the number of periodic points arising from heterodimensional connections. \textit{Compositio Mathematica}. \textbf{157}(9) (2021), 1899--1963.

\bibitem{AvilaLyubichZhang} A. Avila, M. Lyubich, Z. Zhang. To appear.

\bibitem{Barrientos1} P.G Barrientos and A. Raibekas. Robust tangencies of large codimension. \textit{Nonlinearity}. \textbf{30} (2017), 4369--4409.

\bibitem{Barrientos3} P.G Barrientos and A. Raibekas, Robust degenerate unfoldings of cycles and tangencies. \textit{J. Dynam. Differential Equations}. \textbf{33}(1) (2021), 177--209.

\bibitem{BenedicksMartensPalmisano2018} M. Benedicks, M. Martens and L. Palmisano. Newhouse Laminations. \textit{arXiv:1811.00617} (2018).

\bibitem{Berger2016} P. Berger. Generic family with robustly infinitely many sinks. \textit{Inventiones mathematicae}. \textbf{205}(1) (2016), 121--172.

\bibitem{Bergeremergence2017} P. Berger. Emergence and non-typicality of the finiteness of the attractors in many topologies. \textit{Proceedings of the Steklov Institute of Mathematics}. \textbf{297}(1) (2017), 1--27.

\bibitem{Berger2026} P. Berger. Wild dynamics on manifolds. \textit{International Congress of Mathematicians 2026}, 24--41.

\bibitem{BergerCrovisierPujals2022} P. Berger, S. Crovisier and E. Pujals. Germ-typicality of the coexistence of infinitely many sinks. \textit{Advances in Mathematics.} \textbf{407} (2022), 108528.

\bibitem{Berger1} P. Berger, A. Florio and D. Peralta-Salas. Steady Euler flows on $\mathbb R^3$ with wild and universal dynamics. \textit{Comm. Math. Phys.} \textbf{401} (2023), 1--47.

\bibitem{Bergereverydiffeomorphism} P. Berger, N. Gourmelon and M. Helfter. Every diffeomorphism is a total renormalization of a close to identity map. \textit{Inventiones mathematicae}. \textbf{239}(2) (2025), 431--468.

\bibitem{Berger2025} P. Berger and D. Turaev. On Kolmogorov-typical properties of symplectic dynamics. \textit{arXiv:2507.11375} (2025).

\bibitem{Bonatti2003} C. Bonatti and L. D{\'\i}az. On maximal transitive sets of generic diffeomorphisms. \textit{Publ. Math. Inst. Hautes {\'E}tudes Sci.} \textbf{96} (2003), 171--197.

\bibitem{Bonatti2012} C. Bonatti and L. D{\'\i}az. Abundance of $C^1$-robust homoclinic tangencies. \textit{Transactions of the American Mathematical Society}. \textbf{364}(10) (2012), 5111--5148.

\bibitem{Buzzard1997} G.T. Buzzard. Infinitely many periodic attractors for holomorphic maps of 2 variables. \textit{Annals of mathematics}. \textbf{145}(2) (1997), 389--417.

\bibitem{Callis2022} K.G. Callis. Absolutely periodic billiard orbits of arbitrarily high order. \textit{arXiv:2209.11721} (2022).

\bibitem{Colli1998} E. Colli. Infinitely many coexisting strange attractors. \textit{Annales de l'Institut Henri Poincar{\'e}.} \textbf{15}(5) (1998), 539--579.

\bibitem{Domitrz} W. Domitrz, P. Mormul and P. Pragacz. Order of tangency between manifolds. \textit{In Schubert Calculus and Its Applications in Combinatorics and Representation Theory}. (2017), 27--42.

\bibitem{Downarowicz} T. Downarowicz and S. Newhouse. Symbolic extensions and smooth dynamical systems. \textit{Inventiones mathematicae}. \textbf{160}(3) (2005), 453--499.

\bibitem{Duarte1994} P. Duarte. Plenty of elliptic islands for the standard family of area preserving maps. \textit{Annales de l'Institut Henri Poincar{\'e}}. \textbf{11}(4) (1994), 359--409.

\bibitem{Duarte1999} P. Duarte. Abundance of elliptic isles at conservative bifurcations. \textit{Dynamics and Stability of Systems.} \textbf{14}(4) (1999), 339--356.

\bibitem{Duarte2008} P. Duarte. Elliptic isles in families of area-preserving maps. \textit{Ergodic Theory and Dynamical Systems.} \textbf{28}(6) (2008), 1781--1813.

\bibitem{Garrido2026} M. Garrido and P. Mart{\'\i}n. Universal Dynamics in a family of Celestial Mechanics models. \textit{arXiv:2610.09966} (2026). 

\bibitem{Garrido2025} M. Garrido, P. Mart{\'\i}n and J. Paradela. Parabolic Saddles and Newhouse Domains in Celestial Mechanics. \textit{Communications in mathematical physics.} \textbf{406}(7) (2025), 173.

\bibitem{Gelfreich2010} V. Gelfreich and D. Turaev. Universal dynamics in a neighborhood of a generic elliptic periodic point. \textit{Regular and Chaotic Dynamics}. \textbf{15}(2) (2010), 159--164.

\bibitem{GonchenkoOvsyannikov2010} S.V. Gonchenko and I.I. Ovsyannikov. On bifurcations of three-dimensional diffeomorphisms with a non-transversal heteroclinic cycle containing saddle-foci. \textit{Nelin. Dinam.}. \textbf{6}(1) (2010), 61--77 (in Russian).

\bibitem{GonchenkoShilnikov1986} S.V. Gonchenko and L.P. Shilnikov. Dynamical systems with structurally unstable homoclinic curves. \textit{Dokl. Akad. Nauk SSSR.} \textbf{286}(5) (1986), 1049--1053.

\bibitem{GonchenkoShilnikov1987} S.V. Gonchenko and L.P. Shilnikov. Arithmetic properties of topological invariants of systems with nonstructurally-stable homoclinic trajectories. \textit{Ukrainian Mathematical Journal}. \textbf{39}, (1987) 15--21.

\bibitem{Gonchenkoonmoduli1993} S.V. Gonchenko and L.P. Shilnikov. On moduli of systems with a structurally unstable homoclinic Poincare curve. \textit{Russian Acad. Sci. Izv. Math.} \textbf{41}(3) (1993), 417--445.

\bibitem{Gonchenkoonmodelswith} S.V. Gonchenko, L.P. Shilnikov and D.V. Turaev. On models with non-rough Poincar{\'e} homoclinic curves. \textit{Physica D: Nonlinear Phenomena}. \textbf{62} (1993), 1--14.

\bibitem{Gonchenko3} S.V. Gonchenko, L.P. Shilnikov and D.V. Turaev. On dynamical properties of multidimensional diffeomorphisms from Newhouse regions: I. \textit{Nonlinearity}. \textbf{21}(5) (2008), 923--972.

\bibitem{GonchenkoShilnikovTuraev2009} S.V. Gonchenko, L.P. Shilnikov and D.V. Turaev. On global bifurcations in three-dimensional diffeomorphisms leading to wild Lorenz-like attractors. \textit{Regular and Chaotic Dynamics.} \textbf{14}(1) (2009), 137--147.

\bibitem{Gonchenkoonthreetypes} S.V. Gonchenko and D.V. Turaev. On three types of dynamics and the notion of attractor. \textit{Proceedings of the Steklov Institute of Mathematics}. \textbf{297}(1) (2017), 116--137.

\bibitem{Gonchenkocomplexity} S.V. Gonchenko, D.V. Turaev, P. Gaspard and G. Nicolis. Complexity in the bifurcation structure of homoclinic loops to a saddle-focus. \textit{Nonlinearity}. \textbf{10}(2) (1997), 409--423.

\bibitem{Gonchenko1} S.V. Gonchenko, D.V. Turaev and L.P. Shilnikov. On the existence of Newhouse domains in a neighborhood of systems with a structurally unstable Poincar{\'e} homoclinic curve (the higher-dimensional case). \textit{Russian Acad. Sci. Dokl. Math.} \textbf{47}(2) (1993), 268--273.  

\bibitem{dynamicalphenomena} S.V. Gonchenko, D.V. Turaev and L.P. Shilnikov. Dynamical phenomena in multidimensional systems with a structurally unstable homoclinic Poincar{\'e} curve. \textit{Dokl. Math.} \textbf{47}(3) (1993), 410--415.

\bibitem{GonchenkoTuraevShilnikov1997} S.V. Gonchenko, D.V. Turaev and L. P. Shilnikov. On Newhouse domains of two-dimensional diffeomorphisms that are close to a diffeomorphism with a structurally unstable heteroclinic contour. \textit{Proc. Steklov Inst. Math.}. \textbf{216} (1997), 70--118.

\bibitem{arbitraryorder} S.V. Gonchenko, D.V. Turaev and L.P. Shilnikov. Homoclinic tangencies of an arbitrary order in Newhouse domains. \textit{Itogi Nauki i Tekhniki, Seriya Sovremennaya Matematika iEe Prilozheniya (Moscow: VINITI)} \textbf{67} (1999), 69--128.
\\
S.V. Gonchenko, D.V. Turaev and L.P. Shilnikov. \textit{J. Math. Sci. (N.Y.)} \textbf{105}(1) (2001), 1738--1778 (English transl.).

\bibitem{arbitraryorderconservative} S. Gonchenko, D. Turaev and L. Shilnikov. Homoclinic tangencies of arbitrarily high orders in conservative and dissipative two-dimensional maps. \textit{Nonlinearity}. \textbf{20}(2) (2007), 241--275.

\bibitem{GorodetskiKaloshin2007} A. Gorodetski and V. Kaloshin. How often surface diffeomorphisms have infinitely many sinks and hyperbolicity of periodic points near a homoclinic tangency. \textit{Advances in Mathematics}. \textbf{208}(2) (2007), 710--797.

\bibitem{GorodetskiKaloshin2026} A. Gorodetski and V. Kaloshin. On the Hausdorff dimension of oscillatory motions in three-body problems. \textit{arXiv:2610.25750} (2026).

\bibitem{Gourmelon} N. Gourmelon, D. Li, D. Mints and D. Turaev. Stabilization of Homoclinic Tangencies. \textit{Work in progress}. 

\bibitem{Hunt1992} B.R. Hunt, T. Sauer and J.A. Yorke. Prevalence: a translation-invariant "almost every" on infinite-dimensional spaces. \textit{Bulletin of the American mathematical society}. \textbf{27}(2) (1992), 217--238.

\bibitem{Kaloshin1996} V.Y. Kaloshin. Some prevalent properties of smooth dynamical systems. \textit{Proc. Steklov Inst. Math.} \textbf{213} (1996), 115--140.

\bibitem{Kaloshin1} V.Yu. Kaloshin. Generic diffeomorphisms with superexponential growth of number of periodic orbits. \textit{Comm. Math. Phys.} \textbf{211} (2000), 253--271.

\bibitem{Leal2008} B. Leal. High dimension diffeomorphisms exhibiting infinitely many strange attractors. \textit{Annales de l'Institut Henri Poincar{\'e}.} \textbf{25}(3) (2008), 587--607. 

\bibitem{Li2024} D. Li. $C^1$-robust homoclinic tangencies. \textit{arXiv:2406.12500} (2024).

\bibitem{Lirobustheterodimensionalcycles} D. Li, X. Li, K. Shinohara, D. Turaev. Robust heterodimensional cycles in two-parameter unfolding of homoclinic tangencies. \textit{arXiv:2203.14075} (2024).

\bibitem{MartensPalmisano2026} M. Martens and L. Palmisano. Non-Hyperbolic Chaotic Dynamics: Renormalization and More. \textit{arXiv:2511.12926v4} (2026).

\bibitem{Mints} D. Mints. High-order homoclinic tangencies of corank 2. \textit{Ergodic Theory and Dynamical Systems.} \textbf{46}(1) (2026), 218--262.

\bibitem{Newhouse2} S.E. Newhouse. Non-density of axiom A on $S^2$. \textit{Proc. Symp. Pure Math.}, Amer. Math. Soc. \textbf{14} (1970), 191--202.

\bibitem{Newhouse1974} S.E. Newhouse. Diffeomorphisms with infinitely many sinks. \textit{Topology}. \textbf{13}(1) (1974), 9--18.

\bibitem{Newhouse3} S.E. Newhouse. The abundance of wild hyperbolic sets and non-smooth stable sets for diffeomorphisms. \textit{Publ. Math. Inst. Hautes {\'E}tudes Sci.} \textbf{50} (1979), 101--151.

\bibitem{bifurcationsandstability} S.E. Newhouse, J. Palis and F. Takens. Bifurcations and stability of families of diffeomorphisms. \textit{Publ. Math. Inst. Hautes {\'E}tudes Sci.} \textbf{57} (1983), 5--71.

\bibitem{Newhouseoccurrence} S. Newhouse, D. Ruelle and F. Takens. Occurrence of strange axiom $A$ attractors near quasi periodic flows on $T^m$, $m\ge 3$. \textit{Communications in Mathematical Physics}. \textbf{64}(1) (1978), 35--40.

\bibitem{Palis2000} J. Palis. A global view of dynamics and a conjecture on the denseness of finitude of attractors. \textit{Ast{\'e}risque.} \textbf{261} (2000), 335--347.

\bibitem{Palis1} J. Palis and M. Viana. High dimension diffeomorphisms displaying infinitely many periodic attractors. \textit{Ann. of Math. (2)} \textbf{140}(1) (1994), 207--250.

\bibitem{PughShub1996} C. Pugh and M. Shub. Stable ergodicity and partial hyperbolicity. \textit{In International Conference on Dynamical Systems (Montevideo, 1995).} \textbf{362} (1996), 182--187.

\bibitem{Romero1} N. Romero. Persistence of homoclinic tangencies in higher dimensions. \textit{Ergodic Theory Dynam. Systems}. \textbf{15}(4) (1995), 735--757.

\bibitem{Ruelleonthenature} D. Ruelle and F. Takens. On the nature of turbulence. \textit{Commun. Math. Phys.} \textbf{20} (1971), 167--192.

\bibitem{Shilnikov2} L.P. Shilnikov. On a Poincar{\'e}--Birkhoff problem. \textit{Math. USSR-Sb.} \textbf{3}(3) (1967), 353--371.

\bibitem{Shilnikov1} L.P. Shilnikov, A.L. Shilnikov, D.V. Turaev and L.O. Chua. \textit{Methods of Qualitative Theory in Nonlinear Dynamics, part I}. World Scientific, 1998.

\bibitem{Simon} C. P. Simon. Instability in $Diff^r(T^3)$ and the nongenericity of rational zeta functions. \textit{Transactions of the American Mathematical Society}. \textbf{174} (1972), 217--242.

\bibitem{Smalediffeomorphismswith} S. Smale. Diffeomorphisms with many periodic points. \textit{Differential and Combinatorial Topology}. (1965), 63--80.

\bibitem{Smale} S. Smale. Differentiable dynamical systems. \textit{Bull. Amer. Math. Soc.} \textbf{73}(6) (1967), 747--817.

\bibitem{TedeschiniLalliYorke1986} L. Tedeschini-Lalli and J.A. Yorke. How often do simple dynamical processes have infinitely many coexisting sinks? \textit{Communications in mathematical physics}. \textbf{106}(4) (1986), 635--657.

\bibitem{Turaev1} D. Turaev. On dimension of non-local bifurcational problems. \textit{Internat. J. Bifur. Chaos Appl. Sci. Engrg.} \textbf{6}(5) (1996), 919--948.

\bibitem{Turaev2} D. Turaev. Polynomial approximations of symplectic dynamics and richness of chaos in non-hyperbolic area-preserving maps. \textit{Nonlinearity}. \textbf{16}(1) (2002), 123--135.

\bibitem{TuraevEquadiff} D. Turaev. On the genericity of the Newhouse phenomenon. In: EQUADIFF 2003. World
Sci. Publ., Hackensack (2005).

\bibitem{Turaevintro1} D. Turaev. Richness of chaos in the absolute Newhouse domain. \textit{Proceedings of the International Congress of Mathematicians 2010 (ICM 2010)}, 1804--1815.

\bibitem{Turaevintro2} D. Turaev. Maps close to identity and universal maps in the Newhouse domain. \textit{Comm. Math. Phys.} \textbf{335} (2015), 1235--1277.
     
\end{thebibliography}
\end{document}